\documentclass[11pt]{article}

\usepackage{amsmath,amssymb,amsfonts,amsthm}
\usepackage{latexsym}
\usepackage[english]{babel}
\usepackage{showidx}
\usepackage{graphicx}
\usepackage{graphics}
\usepackage{babel,amsmath,amssymb,amsbsy,amsfonts,latexsym,amsthm}
\usepackage{color}

\def\E{{\mathbb E}}

\def\<{\langle}
\def\>{\rangle}

\newtheorem{theorem}{Theorem}[section]

\newtheorem{corollary}[theorem]{Corollary}

\newtheorem{lemma}[theorem]{Lemma}

\newtheorem{proposition}[theorem]{Proposition}
\newtheorem{remark}[theorem]{Remark}

\numberwithin{equation}{section}
\definecolor{mycolorred}{rgb}{1, 0, 0}

\begin{document}
	
\pagenumbering{roman}

\title{Convergence in distribution norms in the CLT\\
for non identical distributed random variables}
\author{\textsc{Vlad Bally}\thanks{%
Universit\'e Paris-Est, LAMA (UMR CNRS, UPEMLV, UPEC), MathRisk INRIA,
F-77454 Marne-la-Vall\'{e}e, France. Email: \texttt{bally@univ-mlv.fr} }%
\smallskip \\
\textsc{Lucia Caramellino}\thanks{%
Dipartimento di Matematica, Universit\`a di Roma ``Tor Vergata'', and
INDAM-GNAMPA, Via della Ricerca Scientifica 1, I-00133 Roma, Italy. Email: 
\texttt{caramell@mat.uniroma2.it}. Corresponding author.}\smallskip\\
\textsc{Guillaume Poly}\thanks{%
IRMAR, Universit\'e de Rennes 1, 263 avenue du G\'en\'eral Leclerc, CS 74205
35042 Rennes, France. Email: \texttt{guillaume.poly@univ-rennes1.fr}
}}
\maketitle

\begin{abstract}
We study the convergence in distribution norms in the Central Limit Theorem
for non identical distributed random variables that is 
\begin{equation*}
\varepsilon _{n}(f):={\mathbb{E}}\Big(f\Big(\frac 1{\sqrt
n}\sum_{i=1}^{n}Z_{i}\Big)\Big)-{\mathbb{E}}\big(f(G)\big)\rightarrow 0
\end{equation*}%
where $Z_{i}$ are centred independent random variables and $G$ is a Gaussian
random variable. We also consider local developments (Edgeworth expansion).
This kind of results is well understood in the case of smooth test functions 
$f$. If one deals with measurable and bounded test functions (convergence in
total variation distance), a well known theorem due to Prohorov shows that
some regularity condition for the law of the random variables $Z_{i}$, $i\in 
{\mathbb{N}}$, on hand is needed. Essentially, one needs that the law of $%
Z_{i}$ is locally lower bounded by the Lebesgue measure (Doeblin's
condition). This topic is also widely discussed in the literature. Our main
contribution is to discuss convergence in distribution norms, that is to
replace the test function $f$ by some derivative $\partial_{\alpha }f$ and
to obtain upper bounds for $\varepsilon _{n}(\partial_{\alpha }f)$ in terms
of the infinite norm of $f$. Some applications are also discussed: an
invariance principle for the occupation time for random walks, small balls
estimates and expected value of the number of roots of trigonometric
polynomials with random coefficients.
\end{abstract}

\medskip


\medskip


\clearpage

\tableofcontents

\clearpage

\parindent 0pt

\pagestyle{plain}

\pagenumbering{arabic}

\section{Introduction}

$\quad$ \textbf{The framework.} We consider $n$ independent (but not
necessarily identically distributed) random variables $Y_{k}$, $k=1,\ldots,n$%
, taking values in ${\mathbb{R}}^m$, which are centered and with identity
covariance matrix. Moreover, we consider $n$ matrices $C_{n,k}\in\mathrm{Mat}%
(d\times m)$ and we look to 
\begin{equation}
S_{n}(Y)=\frac 1{\sqrt{n}}\sum_{k=1}^{n}C_{n,k}Y_{k}.  \label{MR1}
\end{equation}%
Our aim is to obtain a Central Limit Theorem (CLT) as well as Edgeworth
developments in this framework. The basic hypotheses are the following. We
assume the normalization condition%
\begin{equation}
\frac 1n\sum_{k=1}^{n}C_{n,k}C_{n,k}^*=\mathrm{Id}_{d},  \label{MR2'-intro}
\end{equation}%
where $*$ denotes transposition and $\mathrm{Id}_{d}\in \mathrm{Mat}(d\times
d)$ is the identity matrix. Moreover we assume that for each $p\in {\mathbb{N%
}}$ there exists a constant ${\mathrm{C}}_{p}(Y)\geq 1$ such that%
\begin{equation}
\max_{1\leq k\leq n}{\mathbb{E}}(|C_{n,k}Y_k| ^{p})\leq {\mathrm{C}}_{p}(Y).
\label{MR2}
\end{equation}

$\quad $ \textbf{The case of smooth test functions.} Let $\left\Vert
f\right\Vert _{k,\infty }$ denote the norm in $W^{k,\infty },$ that is the
uniform norm of $f$ and of all its derivatives of order less or equal to $k.$
First, we want to prove that%
\begin{equation}
\Big|{\mathbb{E}}(f(S_{n}(Y))-\int_{{\mathbb{R}}^{d}}f(x)\gamma _{d}(x)dx%
\Big|\leq \frac{{\mathrm{C}}_{0}}{n^{\frac{1}{2}}}\left\Vert f\right\Vert
_{3,\infty }  \label{MR5}
\end{equation}%
where $\gamma _{d}(x)=(2\pi )^{-d/2}\exp (-\frac{1}{2}\left\vert
x\right\vert ^{2})$ is the density of the standard normal law. This
corresponds to the Central Limit Theorem (hereafter CLT). Moreover we look
for some polynomials $\psi _{n,k}:{\mathbb{R}}^{d}\rightarrow {\mathbb{R}}$
such that for $N\in {\mathbb{N}}$ and for every $f\in C_{b}^{\widehat{N}}({%
\mathbb{R}}^{d})$, with $\widehat{N}=N(2\lfloor N/2\rfloor +N+5)$, 
\begin{equation}
\Big|{\mathbb{E}}(f(S_{n}(Z))-\int_{{\mathbb{R}}^{d}}f(x)\Big(%
1+\sum_{k=1}^{N}\frac{1}{n^{k/2}}\psi _{n,k}(x)\Big)\gamma _{d}(x)dx\Big|%
\leq \frac{{\mathrm{C}}_{N}}{n^{\frac{1}{2}(N+1)}}\left\Vert f\right\Vert _{%
\widehat{N},\infty }.  \label{MR3}
\end{equation}

This is Theorem \ref{Smooth}, giving the Edgeworth development of order $N$.
In the case of smooth test functions $f$ (as it is the case in (\ref{MR3})),
this topic has been widely discussed and well understood. One should mention
the seminal paper by Essen \cite{[E]} the books of Gnedenko and Kolmogorov 
\cite{[GK]}, Petrov \cite{[Pe]}, Battacharaya and Rao \cite{[BR]} and
Zolotarev \cite{[Z]}. Such development has been obtained by Sirazhdinov and
Mamatov \cite{[SM]} in the case of identically distributed random variables
and then by G\"{o}tze and Hipp \cite{GH} in the non identically distributed
case. A complete presentation of this topic may be found in the recent
review paper by Bobkov \cite{[B1]}. The coefficients $\psi _{n,k}$ in the
development (\ref{MR3}) are linear combinations of Hermite polynomials. An
explicit expression, in the one dimensional case, is given in \cite{[B1]}.
Ourself we give the explicit formula of these coefficients in the multi -
dimensional case. This is important because, in the working paper \cite%
{[BCP]}, the development of order three, in the bi-dimensional case, is used
in order to study invariance principles for the variance of trigonometric
polynomials.

It is worth to mention that the classical approach is based on Fourier
analysis. In our paper we use a different approach based on the Lindeberg
method for Markov semigroups (this is inspired from works concerning the
parametrix method for Markov semigroups in \cite{[BH]}, see also Chatterjee 
\cite{[C]}). This alternative approach is convenient for the proof of our
main result concerning \textquotedblleft distribution
norms\textquotedblright (see below).

$\quad$ \textbf{The case of general test functions.} A second problem is to
obtain the estimate (\ref{MR3}) for test functions $f$ which are not
regular, in particular to replace $\left\Vert f\right\Vert _{\widehat{N}%
,\infty }$ by $\left\Vert f\right\Vert _{\infty }.$ This amounts to estimate
the error in total variation distance. In the case of identically
distributed random variables, and for $N=0$ (so at the level of the standard
CLT), this problem has been widely studied. First of all, one may prove the
convergence in Kolmogorov distance, that is for $f=1_{D}$ where $D$ is a
rectangle. Many refinements of this type of result has been obtained by
Battacharaya and Rao and they are presented in \cite{[BR]}. But it turns out
that one may not prove such a result for a general measurable set $D$
without assuming more regularity on the law of $Y_{k}$, $k\in {\mathbb{N}}.$

Indeed, consider the standard CLT, so take $m=d$, $C_{n,k}=\mathrm{Id}_{d}$
and $Y_{k}$, $k=1,\ldots ,n$, i.i.d. In his seminal paper \cite{[P]}
Prohorov proved that the convergence in total variation distance is
equivalent to the fact that there exists $r$ such that the law of $%
Y_{1}+\cdots +Y_{r}$ has an absolutely continuous component. This is
\textquotedblleft essentially\textquotedblright\ equivalent to the Doeblin's
condition that we present now (see Remark \ref{PPP}): we assume that there
exists $r,\varepsilon >0$ and there exists $y_{k}\in {\mathbb{R}}^{m}$ such
that for every measurable set $A\subset B_{r}(y_{k})$ 
\begin{equation}
{\mathbb{P}}(Y_{k}\in A)\geq \varepsilon \lambda (A)  \label{MR4}
\end{equation}%
where $\lambda $ is the Lebesgue measure. Under (\ref{MR4}) we are able to
obtain (\ref{MR3}) in total variation distance. 

Let us finally mention another line of research which has been strongly
developed in the last years: it consists in estimating the convergence in
the CLT in entropy distance. This starts with the papers of Barron \cite{[B]}
and Johnson and Barron \cite{[BJ]}. In these papers the case of identically
distributed random variables is considered, but recently,
Bobkov, Chistyakov and G\"{o}tze  \cite{[BCG]} have obtained the estimate in entropy distance
for the case of random variables which are no more identically distributed
as well. We recall that the convergence in entropy distance implies the
convergence in total variation distance, so such results are stronger.
However, in order to work in entropy distance one has to assume that the law
of $Z_{n,k}=C_{n,k}Y_{k}$ is absolutely continuous with respect to the
Lebesgue measure and have finite entropy and this is more limiting than (\ref%
{MR4}). So the hypotheses and the results are slightly different. Finally, 
other types of distances ($W_{p}$-transport distances) have been recently
studied in \cite{[B2], [Bo],[R]}.

\smallskip

$\quad $ \textbf{Convergence in distribution norms.}  Consider
first the particular case when $Z_{n,k}=C_{n,k}Y_{k}$ are identically
distributed and have a density which is one time differentiable with
derivative belonging to $L^{1}.$ Then the law of $S_{n}(Y)$ is absolutely
continuous with $C^{n}$ density and then, in Proposition \ref{D}, we prove
that for every $k\in {\mathbb{N}}$ and every multiindex $\alpha $%
\begin{equation*}
\sup_{x}(1+\left\vert x\right\vert ^{2})^{k}\left\vert \partial _{\alpha
}p_{S_{n}}(x)-\partial _{\alpha }\gamma _{d}(x)\right\vert \leq \frac{C}{%
\sqrt{n}}
\end{equation*}%
which is the standard convergence in distribution norms. Notice also that
here we are at the level of the CLT and we are not able to deal with
Edgeworth expansions.

Unfortunately we fail to obtain such a result in the general framework
(which is the interesting case): this is moral because we do not assume that
the laws of $C_{n,k}Y_{k}$, $k=1,...,n$ are absolutely continuous, and then
the law of $S_{n}(Y)$ may have atoms. However we obtain a similar result,
but we have to keep a \textquotedblleft small error\textquotedblright . Let
us give a precise statement of our result. For a function $f\in C_{p}^{q}({%
\mathbb{R}}^{d})$ ($q$ times differentiable with polynomial growth) we
define $L_{q}(f)$ and $l_{q}(f)$ to be two constants such that 
\begin{equation}
\sum_{0\leq \left\vert \alpha \right\vert \leq q}\left\vert \partial
_{\alpha }f(x)\right\vert \leq L_{q}(f)(1+\left\vert x\right\vert
)^{l_{q}(f).}.  \label{Norm1}
\end{equation}%
Our main result is given in Theorem \ref{MAIN} and says the following: for a
fixed $q\in {\mathbb{N}}$, there exist some constants ${\mathrm{C}}_{N}\geq
1\geq {\mathrm{c}}_{N}>0$ (depending on $r,\varepsilon $ from (\ref{MR4})
and on ${\mathrm{C}}_{p}(Y)$ from (\ref{MR2})) such that for every
multiindex $\gamma $ with $\left\vert \gamma \right\vert =q$ and for every $%
f\in C_{p}^{q}({\mathbb{R}}^{d})$%
\begin{equation}
\begin{array}{rcl}
&  & \displaystyle\Big|{\mathbb{E}}\Big(\partial _{\gamma }f(S_{n}(Z))\Big)%
-\int_{{\mathbb{R}}^{d}}\partial _{\gamma }f(x)\Big(1+\sum_{k=1}^{N}\frac{1}{%
n^{k/2}}\psi _{n,k}(x)\Big)\gamma _{d}(x)dx\Big|\smallskip \\ 
&  & \displaystyle\leq C_{N}\Big(L_{q}(f)e^{-{\mathrm{c}}_{N}\times n}+\frac{%
1}{n^{\frac{1}{2}(N+1)}}L_{0}(f)\Big).%
\end{array}
\label{M3''}
\end{equation}%
However we fail to get convergence in distribution norms because $%
L_{q}(f)e^{-{\mathrm{c}}_{N}\times n}$ appears in the upper bound of the
error and $L_{q}(f)$ depends on the derivatives of $f$. But we are close to
such a result: notice first that if $f_{n}=f\ast \phi _{\delta _{n}}$ is a
regularization by convolution with $\delta _{n}=\exp (-\frac{{\mathrm{c}}_{N}%
}{2q}\times n)$ then (\ref{M3''}) gives 
\begin{equation}
\Big|{\mathbb{E}}\Big(\partial _{\gamma }f_{n}(S_{n}(Z))\Big)-\int_{{\mathbb{%
R}}^{d}}\partial _{\gamma }f_{n}(x)\Big(1+\sum_{k=1}^{N}\frac{1}{n^{k/2}}%
\psi _{n,k}(x)\Big)\gamma _{d}(x)dx\Big|\leq \frac{{\mathrm{C}}_{N}}{n^{%
\frac{1}{2}(N+1)}}L_{0}(f).  \label{M3'}
\end{equation}

\smallskip

We discuss now three applications.

\smallskip

$\quad$ \textbf{Application 1: an invariance principle related to the local
time.} Let 
\begin{equation*}
S_n(k,Y)=\frac 1{\sqrt n}\sum_{i=1}^k Y_i,
\end{equation*}
where $Y_1,\ldots,Y_n$ are independent and identically distributed random
variables. We set $\varepsilon _{n}=n^{-\frac{1}{2}(1-\rho )}$ with $\rho
\in (0,1)$ and in Theorem \ref{OT} we prove that, for every $\rho ^{\prime
}<\rho$, 
\begin{equation*}
\Big| \frac 1n\sum_{k=1}^n{\mathbb{E}}\Big(\frac 1{2\varepsilon_n}
1_{\{\vert S_n(k,Y)\vert \leq \varepsilon _{n}\}}\Big)- {\mathbb{E}}\Big(%
\int_{0}^{1}\frac 1{2\varepsilon _n}1_{\{\left\vert W_s\right\vert \leq
\varepsilon _{n}\}}ds\Big)\Big| \leq \frac{{\mathrm{C}}}{n^{\frac{%
1+\rho^{\prime }}{2}}}.
\end{equation*}%
with $W_{s}$ a Brownian motion (we recall that $ \int_{0}^{1}\frac 1{2\varepsilon _n}1_{\{\left\vert W_s\right\vert \leq \varepsilon _{n}\}}ds$ 
converges to the local time of $W).$ Here the test function is $f_{n}(x)=%
\frac{1}{2\varepsilon _{n}}1_{|x|<\varepsilon _{n}}$ and this converges to
the Dirac function. This example shows that (\ref{M3''}) is an appropriate
estimate in order to deal with some singular problems.

\smallskip

$\quad$ \textbf{Application 2: small ball probabilities.} We consider the
case in which the matrices $C_{n,k}$ can depend on a parameter $u\in {%
\mathbb{R}}^\ell$, that is, 
\begin{equation*}
S_{n}(u,Y)=\frac{1}{\sqrt{n}}\sum_{k=1}^{n}C_{n,k}(u)Y_{k},\quad u\in{%
\mathbb{R}}^\ell .
\end{equation*}
We assume that $u\mapsto C_{n,k}(u)\in\mathrm{Mat}(d\times d)$ is twice
differentiable with bounded derivatives up to order two and that the
covariance matrix field of $S_{n}(u,Y)$ is the identity matrix, that is, $%
\Sigma_{n}(u)=\frac 1n\sum_{k=1}^nC_{n,k}(u)C_{n,k}^*(u)=\mathrm{Id}_d$.
Then in Theorem \ref{SB} we prove the following estimate: if $d>\ell$, $%
a\geq 0$ and $\theta >\frac{a\ell }{d-\ell }$ then, for every $\varepsilon
>0 $, 
\begin{equation}  \label{SB-intr}
{\mathbb{P}}\Big(\inf_{\left\vert u\right\vert \leq n^{a}}\left\vert
S_{n}(u,Y)\right\vert \leq \frac{1}{n^{\theta }}\Big)\leq \frac{{\mathrm{C}}%
}{n^{\theta (d-\ell )-a\ell -\varepsilon }}.
\end{equation}
This is done by applying (\ref{M3'}) to the multiindex $\gamma=(1,\ldots,d)$
and the function $f=f_n$, with 
\begin{equation*}
f_{n }(x)=n ^{-\theta d}\int_{-\infty }^{x_{1}}dx_{2}...\int_{-\infty
}^{x_{d-1}}dx_{d}1_{\{|x|<n^{-\theta}\}}(x).
\end{equation*}
Then (\ref{M3'}) allows one to replace $S_{n}(u,Y)$ with a Gaussian random
variable, and in this case we have a nice estimate of the error. We
emphasize that, contrarily to the case of supremas of random processes, much
less is known regarding infimas. As such, the last result can be seen as a
preliminary step enabling one to switch to the Gaussian case for which more
accurate tools are available.

\smallskip

$\quad$ \textbf{Application 3: an invariance principle for the expected
roots of trigonometric polynomials.} Let $N_n(Y)$ be the number of roots in $%
(0,\pi)$ of the polynomial 
\begin{equation*}
P_{n}(t,Y)=\sum_{k=1}^{n}\big(Y_{k}^{1}\cos (kt)+Y_{k}^{2}\sin (kt)\big).
\end{equation*}
It is known, see e.g. \cite{D}, that if the $Y_k$'s are replaced by independent
and identically distributed standard normal random variables $G_k$'s then 
\begin{equation*}
\lim_{n}\frac{1}{n}{\mathbb{E}}(N_{n}(G))=\frac{1}{\sqrt 3}.
\end{equation*}%
Note that the aforementioned asymptotic still holds when the Gaussian
coefficients display some strong form of dependence \cite{ADP17}. In the
recent paper \cite{F}, the above result has been proved for general
independent and identically distributed random variables $Y_k$, $k\in {%
\mathbb{N}}$, which are centered and with variance one. In Theorem \ref{roots}
we drop the assumption of being identically distributed: we prove that the
same limit holds for $N_n(Y)$ when the $Y_k$'s are independent and fulfill
the Doeblin's condition. We stress that it is not completely clear whether
the strategy used in \cite{F} can be adapted to the setting of
non-identically distributed coefficients since it is explicitly used at
several moments in the proof that the characteristic function of each
coefficients behaves \textit{in a same way} near the origin, which is more
restrictive that our normalization condition (\ref{MR2'-intro}). Our main
result enters in the following way: thanks to the Kac-Rice formula, we have 
\begin{equation*}
N_n(Y)=\lim_{\delta\to 0} \int_{a}^{b}\left\vert \partial_t
P_n(t,Y)\right\vert 1_{\{\left\vert P_n(t,Y)\right\vert \leq \delta \}}\frac{%
dt}{2\delta },
\end{equation*}
so we apply (\ref{M3''}) to the pair $(\partial_t P_n(t,Y), P_n(t,Y))$.
Although the article only focuses on the expectation, we stress that this
methods paves the way to an investigation of higher moments (and hence
variance or CLT's) by using Kac-Rice formulas of higher order. This is
actually the main content of the forthcoming article \cite{[BCP]} which
follows the series \cite{GW,AL13,ADL16} of articles dedicated to this task
in the Gaussian case.

\section{Notation and main results}

\label{main-results}

We fix $n\in{\mathbb{N}}$ and we consider $n$ independent random variables $%
\{Y_k\}_{1\leq k\leq n}$, with $Y_{k}=(Y_{k}^{1},\ldots,Y_{k}^{m})\in {%
\mathbb{R}}^{m}$, which are centered and whose covariance matrix is the
identity. Let $\{C_{n,k}\}_{1\leq k\leq n}$ denote $n$ matrices in $\mathrm{%
Mat}(d\times m)$ and set 
\begin{equation*}
\sigma _{n,k}=C_{n,k}C^\ast_{n,k}\in \mathrm{Mat}(d\times d),
\end{equation*}
$\ast$ denoting transposition, so $\sigma_{n,k}$ is the covariance matrix of
the random variable $C_{n,k}Y_k$. We define 
\begin{equation}
S_{n}(Y)=\frac 1{\sqrt n}\sum_{k=1}^{n}C_{n,k}Y_k.  \label{BD1}
\end{equation}%
Sometimes, but not everywhere, we consider the normalizing condition 
\begin{equation}
\frac 1n\sum_{k=1}^{n}\sigma _{n,k}=\mathrm{Id}_{d},  \label{MR2'}
\end{equation}%
$\mathrm{Id}_{d}$ denoting the $d\times d$ identity matrix. Our aim is to
compare the law of $S_{n}(Y)$ with the law of $S_{n}(G)$ where $%
G=(G_{k})_{1\leq k\leq n}$ denote $n$ standard independent Gaussian random
variables. This is a CLT result (but we stress that it is not asymptotic)
and we will obtain an Edgeworth development as well.

We assume that $Y_{k}$ has finite moments of any order and more precisely, 
\begin{equation}
\max_{1\leq k\leq n}{\mathbb{E}}(\left\vert C_{n,k}Y_{k}\right\vert
^{p})\leq {\mathrm{C}}_{p}(Y),\quad \forall \ p.  \label{bd00}
\end{equation}%
Notice that by (\ref{MR2'}) $|\sigma _{n,k}^{i,j}|\leq 1$ so we may assume
without loss of generality that ${\mathbb{E}}(\left\vert
C_{n,k}G_{k}\right\vert ^{p})\leq {\mathrm{C}}_{p}(Y)$ for the standard
normal random variables as well.

\bigskip

\subsection{Doeblin's condition and Nummelin's splitting}

\label{doeb}

We say that the law of the random variable $Y\in {\mathbb{R}}^{m}$ is
locally lower bounded by the Lebesgue measure if there exists $y_{Y}\in {%
\mathbb{R}}^{d}$ and $\varepsilon ,r>0$ such that for every non negative and
measurable function $f:{\mathbb{R}}^{d}\rightarrow {\mathbb{R}}_{+}$ 
\begin{equation}
\begin{array}{ll}
\qquad & {\mathbb{E}}(f(Y))\geq \varepsilon \int
f(y-y_{Y})1_{B(0,2r)}(y-y_{Y})dy.%
\end{array}
\label{s1}
\end{equation}%
(\ref{s1}) is known as the Doeblin's condition. We denote by ${\mathfrak{D}}%
(r,\varepsilon )$ the class of the random variables which verify (\ref{s1}).
Given $r>0$\ we consider the functions $a_{r},\psi _{r}:{\mathbb{R}}%
\rightarrow {\mathbb{R}}_{+}$ defined by 
\begin{equation}
a_{r}(t)=1-\frac{1}{1-(\frac{t}{r}-1)^{2}}\qquad \psi
_{r}(t)=1_{\{\left\vert t\right\vert \leq r\}}+1_{\{r<\left\vert
t\right\vert \leq 2r\}}e^{a_{r}(\left\vert t\right\vert )}.  \label{s2}
\end{equation}%
If $Y\in {\mathfrak{D}}(r,\varepsilon )$ then%
\begin{equation*}
{\mathbb{E}}(f(Y))\geq \varepsilon \int f(y-y_{Y})\psi _{r}(\left\vert
y-y_{Y}\right\vert ^{2})dy.
\end{equation*}%
The advantage of $\psi _{r}(\left\vert y-y_{Y}\right\vert ^{2})$ is that it
is a smooth function (which replaces the indicator function of the ball) and
(it is easy to check) that for each $l\in {\mathbb{N}},p\geq 1$ there exists
a universal constant ${\mathrm{C}}_{l,p}\geq 1$ such that 
\begin{equation}
\psi _{r}(t)|a_{r}^{(l)}(|t|)|^{p}\leq \frac{{\mathrm{C}}_{l,p}}{r^{lp}}
\label{s3}
\end{equation}%
where $a_{r}^{(l)}$ denotes the derivative of order $l$ of $a_{r}.$ Moreover
one can check (see \cite{[CLT]}) that if $Y\in {\mathfrak{D}}(r,\varepsilon
) $ then it admits the following decomposition (the equality is understood
as identity of laws): 
\begin{equation}
Y=\chi V+(1-\chi )U  \label{s4}
\end{equation}%
where $\chi ,V,U$ are independent random variables with the following laws: 
\begin{equation}
\begin{array}{c}
\displaystyle{\mathbb{P}}(\chi =1)=\varepsilon {\mathfrak{m}}_{r}\quad %
\mbox{and}\quad {\mathbb{P}}(\chi =0)=1-\varepsilon {\mathfrak{m}}%
_{r},\smallskip \\ 
\displaystyle{\mathbb{P}}(V\in dy)=\frac{1}{{\mathfrak{m}}_{r}}\psi
_{r}(\left\vert y-y_{Y})\right\vert ^{2})dy\smallskip \\ 
\displaystyle{\mathbb{P}}(U\in dy)=\frac{1}{1-\varepsilon {\mathfrak{m}}_{r}}%
({\mathbb{P}}(Z\in dy)-\varepsilon \psi _{r}(\left\vert y-y_{Y}\right\vert
^{2})dy)%
\end{array}
\label{s5}
\end{equation}%
with%
\begin{equation}
{\mathfrak{m}}_{r}=\int \psi _{r}(\left\vert y-y_{Y}\right\vert ^{2})dy.
\label{s6}
\end{equation}%
The decomposition (\ref{s4}) is also known as the Nummelin's splitting. We
will see later on, specifically in next Section \ref{sect-mall}, that the
noise coming from the Nummelin's decomposition allows one to set-up a
Malliavin type calculus, which in turn will be our main tool in order to get
our CLT result in distribution norms.

\begin{remark}
\label{PPP}In his seminal paper \cite{[P]} Prohorov considers a sequence of $%
i.i.d$ random variables $X_{n}$ and proves that the convergence in the CLT
holds in total variation distance if and only if the following hypothesis
holds: there exists $n_{\ast }$ such that the law of $X_{1}+...+X_{n_{\ast
}} $ has an absolute continuous component, that is $X_{1}+...+X_{n_{\ast
}}\sim \mu (dx)+p(x)dx.$ Of course this is much weaker than Doeblin's
condition, but, as long as we want to prove the CLT in total variation
distance, we may proceed as follows: we denote $Y_{k}=X_{kn_{\ast
}+1}+...+X_{(k+1)n_{\ast }}$ and take $Z_{k}=Y_{2k}+Y_{2k+1}.$ Since the
convolution of two functions from $L^{1}$ is a continuous function, $p\ast p$
is continuous and consequently locally lower bonded by the Lebesgue measure.
So $Z_{k}$ verifies Doeblin's condition. We prove the CLT in total
variation for $Z_{k}$ and then it easily follows for $X_{n}$ (see Corollary %
\ref{Prohorov} below). So, as long as one is concerned with the CLT the
two conditions are (in the above sense) equivalent.
\end{remark}

\begin{remark}
We stress that in \cite{[CLT]} Proposition 2.4 there is fault: it is asserted
that, if $X\sim \mu (dx)+p(x)dx$ then $X$ satisfies the Doeblin's condition
-- and of course this is false if we do not ask $p$ to be lower
semicontinous. However, in Lemma $A.1$ from the appendix in the same paper,
the lower continuity hypothesis is mentioned.
\end{remark}

\bigskip

\subsection{Main results}

\label{main-res}

In order to give the expression of the terms which appear in the Edgeworth
development we need to introduce some notation.

We say that $\alpha$ is a multiindex if $\alpha\in\{1,\ldots,d\}^k$ for some 
$k\geq 1$, and we set $|\alpha|=k$ its length. We allow the case $k=0$,
giving the void multiindex $\alpha=\emptyset$.

Let $\alpha $ be a multiindex and set $k=|\alpha |$. For for $x\in {\mathbb{R%
}}^{d}$ and $f\,:\,{\mathbb{R}}^{d}\rightarrow {\mathbb{R}}$, we denote $%
x^{\alpha }=x_{\alpha _{1}}\cdots x_{\alpha _{k}}$ and $\partial _{\alpha
}f(x)=\partial _{x_{\alpha _{1}}}\cdots \partial _{x_{\alpha _{k}}}f(x)$,
the case $k=0$ giving $x^{\emptyset }=1$ and $\partial _{\emptyset }f=f$. In
the following, we denote with $C^{k}({\mathbb{R}}^{d})$ the set of the
functions $f$ such that $\partial _{\alpha }f$ exists and is continuous for
any $\alpha $ with $|\alpha |\leq k$. The set $C_{p}^{k}({\mathbb{R}}^{d})$,
resp. $C_{b}^{k}({\mathbb{R}}^{d})$, is the subset of $C^{k}({\mathbb{R}}%
^{d})$ such that $\partial _{\alpha }f$ has polynomial growth, resp. is
bounded, for any $\alpha $ with $|\alpha |\leq k$. $C^{\infty }({\mathbb{R}}%
^{d})$, resp. $C_{p}^{\infty }({\mathbb{R}}^{d})$ and $C_{b}^{\infty }({%
\mathbb{R}}^{d})$, denotes the intersection of $C^{k}({\mathbb{R}}^{d})$,
resp. of $C_{p}^{k}({\mathbb{R}}^{d})$ and of $C_{b}^{k}({\mathbb{R}}^{d})$,
for every $k$. For $f\in C_{p}^{k}({\mathbb{R}}^{d})$ we define $L_{k}(f)$
and $l_{k}(f)$ to be some constants such that 
\begin{equation}
\sum_{0\leq \left\vert \alpha \right\vert \leq k}\left\vert \partial
_{\alpha }f(x)\right\vert \leq L_{k}(f)(1+\left\vert x\right\vert
)^{l_{k}(f)}.  \label{Lklk}
\end{equation}%
Notice that if $f\in C_{b}^{\infty }({\mathbb{R}}^{d})$ then $l_{k}(f)=0$
and $L_{k}(f)=\sum_{0\leq \left\vert \alpha \right\vert \leq k}\left\Vert
\partial _{\alpha }f\right\Vert _{\infty }.$

Moreover, for a non negative definite matrix $\sigma \in \mathrm{Mat}(d\times d)$ we denote by $L_{\sigma }$ the Laplace operator associated to $%
\sigma $, i.e. 
\begin{equation}
L_{\sigma }=\sum_{i,j=1}^{d}\sigma ^{i,j}\partial _{z_{i}}\partial _{z_{j}}.
\label{bd0a}
\end{equation}%
For $r\geq 1$ and $l\geq 0$ we set 
\begin{equation}
\Delta _{n,r}(\alpha )={\mathbb{E}}((C_{n,r}Y_r)^{\alpha })-{\mathbb{E}}%
((C_{n,r}G_r)^{\alpha })\quad \mbox{and}\quad D_{n,r}^{(l)}=\sum_{\left\vert
\alpha \right\vert =l}\Delta _{n,r}(\alpha )\partial_{\alpha }.  \label{bd1}
\end{equation}

Notice that $D_{n,r}^{(l)}\equiv 0$ for $l=0,1,2$ and, by (\ref{bd00}), for $%
l\geq 3$ and $\left\vert \alpha \right\vert =l$ then 
\begin{equation}
\left\vert \Delta _{n,r}(\alpha )\right\vert \leq 2{\mathrm{C}}_{l}(Y),\quad
r=1,\ldots,n.  \label{M2'}
\end{equation}%
We construct now the coefficients of our development. Let $N$ be fixed: this
is the order of the development that we will obtain. Given $1\leq m\leq
k\leq N$ we define%
\begin{equation}  \label{M1}
\begin{array}{rcl}
\Lambda _{m} & = & \displaystyle \{((l_{1},l_{1}^{\prime
}),...,(l_{m},l_{m}^{\prime })):N+2\geq l_{i}\geq 3, \lfloor N/2\rfloor \geq
l_{i}^{\prime }\geq 0,i=1,...,m\}, \smallskip \\ 
\Lambda _{m,k} & = & \displaystyle \{((l_{1},l_{1}^{\prime
}),...,(l_{m},l_{m}^{\prime }))\in \Lambda
_{m}:\sum_{i=1}^{m}l_{i}+2\sum_{i=1}^{m}l_{i}^{\prime }=k+2m\}.%
\end{array}%
\end{equation}
Then, for $1\leq k\leq N,$ we define the differential operator%
\begin{equation}
\Gamma _{n,k}=\sum_{m=1}^{k}\sum_{((l_{1},l_{1}^{\prime
}),...,(l_{m},l_{m}^{\prime }))\in \Lambda _{m,k}}\frac 1{n^m}\sum_{1\leq
r_{1}<...<r_{m}\leq n}\prod_{i=1}^{m}\frac{1}{l_{i}!}D_{n,r_{i}}^{(l_{i})}%
\prod_{j=1}^{m}\frac{(-1)^{l_{j}^{\prime }}}{2^{l_{j}^{\prime
}}l_{j}^{\prime }!}L_{\sigma _{n,r_{j}}}^{l_{j}^{\prime }}.  \label{M2}
\end{equation}

By using (\ref{bd00}) and (\ref{M2'}), one easily gets the following
estimates: 
\begin{equation}
\left\vert \Gamma _{n,k}f(x)\right\vert \leq C\times {\mathrm{C}}%
_{3k}(Y)L_{3k}(f)(1+\left\vert x\right\vert )^{l_{3k}(f)},\quad f\in
C_{p}^{3k}({\mathbb{R}}^{d}),  \label{M2''LUCIA}
\end{equation}%
where $L_{3k}(f)$ and $l_{3k}(f)$ are given in (\ref{Lklk}) and $C>0$ is a
suitable constant which does not depend on $n$.

We introduce now the Hermite polynomials, we refer to Nualart \cite{bib:[N]}
for definitions and properties. 
The Hermite polynomial $H_{m}$ of order $m$ on ${\mathbb{R}}$ is defined as 
\begin{equation}
H_{m}(x)=(-1)^{m}e^{\frac{1}{2}x^{2}}\frac{d^{m}}{x^{m}}e^{-\frac{1}{2}%
x^{2}}.  \label{tv1}
\end{equation}%
For a multiindex $\alpha\in\{1,\ldots,d\}^l$ we denote $\beta _{i}(\alpha )={%
\mathrm{card}}\{j:\alpha _{j}=i\}$ and we define the Hermite polynomial on ${%
\mathbb{R}}^{d}$ corresponding to the multiindex $\alpha $ by%
\begin{equation}
H_{\alpha }(x)=\prod_{i=1}^{l}H_{\beta _{i}(\alpha )}(x_{i})\qquad \mbox{for}%
\qquad x=(x_{1},...,x_{d}).  \label{tv2}
\end{equation}%
Equivalently, the Hermite polynomial $H_{\alpha }$ on ${\mathbb{R}}^{d}$
associated to the multiindex $\alpha$ is defined by%
\begin{equation}
{\mathbb{E}}(\partial_{\alpha }f(W))={\mathbb{E}}(f(W)H_{\alpha }(W))\quad
\forall f\in C_{p}^{\infty }({\mathbb{R}}^{d})  \label{M3}
\end{equation}%
where $W$ is a standard normal random variable in ${\mathbb{R}}^{d}$.
Moreover for a differential operator $\Gamma =\sum_{\left\vert \alpha
\right\vert \leq k}a(\alpha )\partial_{\alpha },$ with $a(\alpha )\in {%
\mathbb{R}},$ we denote $H_{\Gamma }=\sum_{\left\vert \alpha \right\vert
\leq k}a(\alpha )H_{\alpha }$ so that 
\begin{equation}
{\mathbb{E}}(\Gamma f(W))={\mathbb{E}}(f(W)H_{\Gamma }(W)).  \label{M4}
\end{equation}%
Finally we define%
\begin{equation}
\Phi _{n,N}(x)=1+\sum_{k=1}^{N}\frac 1{n^{k/2}}H_{\Gamma _{n,k}}(x) 
\mbox{
with $\Gamma _{n,k}$ defined in (\ref{M2}).}  \label{M5}
\end{equation}%
The polynomial $\Phi _{n,N}$ gives the Edgeworth expansion of order $N$ in
the CLT, as stated in the following result, which represents the main result
of this paper.

\begin{theorem}
\label{MAIN} Assume that $Y_{k}\in {\mathfrak{D}}(r,\varepsilon ),\forall
k\in N$ for some $\varepsilon >0,r>0.$ Let the normalizing condition (\ref%
{MR2'}) and the moment bounds condition (\ref{bd00}) both hold. Let $N,q\in {%
\mathbb{N}}$ be fixed. We assume that $n$ is sufficiently large in order to
have 
\begin{equation*}
n^{\frac{1}{2}(N+1)}e^{-\frac{{\mathfrak{m}}_{r}^{2}n}{256}}\leq 1\quad %
\mbox{and}\quad n\geq 4(N+1){\mathrm{C}}_{2}(Y).
\end{equation*}%
There exists $C\geq 1$, depending on $N$ and $q$ only, such that for every
multiindex $\gamma $ with $\left\vert \gamma \right\vert =q$ and every $f\in
C_{p}^{q}({\mathbb{R}}^{d})$%
\begin{equation}
\left\vert {\mathbb{E}}(\partial _{\gamma }f(S_{n}(Y)))-{\mathbb{E}}%
(\partial _{\gamma }f(W)\Phi _{n,N}(W))\right\vert \leq C\times {\mathrm{C}}%
_{\ast }(Y)\Big(\frac{L_{0}(f)}{n^{\frac{1}{2}(N+1)}}+L_{q}(f)e^{-\frac{{%
\mathfrak{m}}_{r}^{2}}{32}\times n}\Big)  \label{s20}
\end{equation}%
where ${\mathrm{C}}_{\ast }(Y)$\ is a constant which depends on $%
q,l_{q}(f),N $ and ${\mathrm{C}}_{p}(Y)$ for $p=2(N+3)\vee 2l_{0}(f)..$
\end{theorem}

\begin{remark}
The precise value of ${\mathrm{C}}_{\ast }(Y)$ is given by 
\begin{eqnarray}
&&%
\begin{array}{l}
{\mathrm{C}}_{\ast }(Y)=\displaystyle\Big(1\vee \frac{8}{{\mathfrak{m}}_{r}}%
\Big)^{2dp_{2}}2^{(N+3)l_{0}(f)+l_{q}(f)}\,\hat{c}%
_{p_{1},l_{0}(f)}c_{l_{q}(f)\vee (l_{0}(f)+p_{1})}\smallskip \\ 
\displaystyle\qquad \qquad \times \frac{{\mathrm{C}}%
_{16dp_{2}}^{(4d+1)p_{2}}(Y)}{r^{p_{2}(p_{2}+1)}}{\mathrm{C}}%
_{2(N+3)}^{(\lfloor N/2\rfloor +2)(N+1)}(Y)\big(1+{\mathrm{C}}%
_{2l_{0}(f)}^{l_{0}(f)\vee (N+1)}(Y)\big)%
\end{array}
\label{s20'} \\
&&%
\begin{array}{l}
\mbox{with }\displaystyle p_{1}=q+(N+1)(N+3),\quad p_{2}=q+p_{1},\smallskip
\\ 
\displaystyle{\mathrm{c}}_{\rho }=\int \left\vert \phi (z)\right\vert
(1+\left\vert z\right\vert )^{\rho }dz,\quad \hat{c}_{p,l}=1\vee \max_{0\leq
|\alpha |\leq p}\int (1+|x|)^{l}|\partial _{\alpha }\phi (x)|dx%
\end{array}
\label{s20''}
\end{eqnarray}%
in which $\phi $ denotes a super kernel (see next (\ref{kk1}) and (\ref{kk2}%
)).
\end{remark}

Actually the coefficients $H_{\Gamma _{n,k}}(x)$ of the polynomial $\Phi
_{n,N}(x)$ are cumbersome. The following corollary, whose proof is postponed in Section \ref{coefficients}, gives a plain expansion
of order three:

\begin{corollary}
\label{coeff} Let the set-up of Theorem \ref{MAIN} holds. For a multiindex $%
\alpha $ and $i,j\in \{1,\ldots ,d\}$, set 
\begin{equation}
c_{n}(\alpha )=\frac{1}{n}\sum_{r=1}^{n}\Delta _{n,r}(\alpha )\quad %
\mbox{and}\quad \overline{c}_{n}(\alpha ,i,j)=\frac{1}{n}\sum_{r=1}^{n}%
\Delta _{n,r}(\alpha )\sigma _{n,r}^{ij}.  \label{coeff-1}
\end{equation}%
Then there exists $C\geq 1$, depending on $N$ and $q$ only, such that for
every multiindex $\gamma $ with $\left\vert \gamma \right\vert =q$ and every 
$f\in C_{p}^{q}({\mathbb{R}}^{d})$%
\begin{equation}
\Big\vert{\mathbb{E}}\Big(\partial _{\gamma }f(S_{n}(Y))\Big)-{\mathbb{E}}%
\Big(\partial _{\gamma }f(W)\Big(1+\sum_{k=1}^{3}\frac{1}{n^{k/2}}\mathcal{H}%
_{n,k}(W)\Big)\Big)\Big\vert\leq C{\mathrm{C}}_{\ast }(Y)\Big(\frac{L_{0}(f)%
}{n^{2}}+L_{q}(f)e^{-\frac{{\mathfrak{m}}_{r}^{2}}{32}\times n}\Big)
\label{s20-3}
\end{equation}%
where ${\mathrm{C}}_{\ast }(Y)$ is given in (\ref{s20''}) and 
\begin{align}
&\mathcal{H}_{n,1}(x)=\frac{1}{6}\sum_{\left\vert \alpha \right\vert
=3}c_{n}(\alpha )H_{\alpha }(x),  \label{H1} \\
&\mathcal{H}_{n,2}(x)=\frac{1}{24}\sum_{\left\vert \alpha \right\vert
=4}c_{n}(\alpha )\,H_{\alpha }(x)+\frac{1}{72}\sum_{\left\vert \alpha
\right\vert =3}\sum_{\left\vert \beta \right\vert =3}c_{n}(\alpha
)c_{n}(\beta )H_{(\alpha ,\beta )}(x),  \label{H2} \\
&%
\begin{array}{l}
\mathcal{H}_{n,3}(x)=\displaystyle-\frac{1}{12}\sum_{\left\vert \alpha
\right\vert =3}\sum_{i,j=1}^{2}\overline{c}_{n}(\alpha ,i,j)H_{(\alpha
,\beta )}(x)+\frac{1}{120}\sum_{\left\vert \alpha \right\vert
=5}c_{n}(\alpha )H_{\alpha }(x)\smallskip \\ 
\quad \displaystyle+\frac{1}{144}\sum_{\left\vert \alpha \right\vert
=3}\sum_{\left\vert \beta \right\vert =4}c_{n}(\alpha )c_{n}(\beta
)H_{(\alpha ,\beta )}+\frac{1}{1296}\sum_{\left\vert \alpha \right\vert
=3}\sum_{\left\vert \beta \right\vert =3}\sum_{\left\vert \gamma \right\vert
=3}c_{n}(\alpha )c_{n}(\beta )c_{n}(\gamma )H_{(\alpha ,\beta ,\gamma )}(x).%
\label{H3}
\end{array}
\end{align}
\end{corollary}

\begin{remark}
We stress that the coefficients of the Hermite polynomials appearing in $%
\mathcal{H}_{n,1}(x)$--$\mathcal{H}_{n,3}(x)$ 
depend on $n$ (this is because we work with $C_{n,k}Y_{k}$, $k=1,\ldots ,n$,
whose law depends on $n$) but in a bounded way. In fact, by the formula (\ref%
{coeff-1}) and by (\ref{M2'}), for $|\alpha |=l$ and $i,j\in \{1,\ldots ,d\}$%
, 
\begin{equation*}
|c_{n}(\alpha )|\leq 2{\mathrm{C}}_{l}(Y)\mbox{ and }|\overline{c}%
_{n}(\alpha ,i,j)|\leq 4{\mathrm{C}}_{l}(Y){\mathrm{C}}_{2}(Y),\quad %
\mbox{for every }n.
\end{equation*}
\end{remark}

\begin{remark}
In the one dimensional case Bobkov obtained in \cite{[B1]} (see Proposition
14.1 therein) the following development using Hermite polynomials:%
\begin{equation*}
\Phi _{n,N}(x)=1+\sum \frac{1}{k_{1}!...k_{N}!}\left( \frac{\gamma _{n,3}}{3!%
}\right) ^{k_{1}}...\left( \frac{\gamma _{n,N+2}}{(N+2)!}\right)
^{k_{N}}\times H_{k}(x)
\end{equation*}%
where $k=3k_{1}+...+(N+2)k_{N}$ and the summation is made over the non
negative integers $k_{1},...k_{N}$ such that $0<k_{1}+2k_{2}+...+Nk_{N}\leq
N.$ And $\gamma _{n,p}$ is the $p$-cumulant of $S_{n}(Y).$ This is an
alternative way to write the correctors which is ordered according to the
powers of the Hermite polynomials (and of course, the two expressions are
equivalent and one may pass from one to another).
\end{remark}

The proof of Theorem \ref{MAIN} is done by using a Malliavin type calculus
based on the random variables $V_{k}$'s coming from the Nummelin's splitting
associated to the $Y_k$'s. This differential calculus is developed in next
Section \ref{sect-mall}. The proof of Theorem \ref{MAIN} represents the main
effort in this paper, so we postpone it to Section \ref{sect-proof-MAIN}. As
for Corollary \ref{coeff}, the proof consists in heavy but straightforward
computations, so we postpone in Section \ref{coefficients}.

\smallskip

We give now two slight variants of Theorem \ref{MAIN} which will be used in
the following. First:

\begin{proposition}
\label{W} Let (\ref{MR2'}) and (\ref{bd00}) hold. Assume that for some $%
n_*<n $ one has $Y_{k}\in {\mathfrak{D}}(r,\varepsilon )$ for $k\leq n-n_*$
and $\frac 1n\sum_{k=1}^{n-n_*}\sigma _{n,k}\geq \frac{1}{2}\mathrm{Id}_d.$
Then (\ref{s20}) holds true.
\end{proposition}

The proof of Proposition \ref{W} mimics the one of Theorem \ref{MAIN} so we
postpone it as well, in next Section \ref{sect-proof-W}. This result will be
used in the proof of Corollary \ref{Prohorov} below.

Let us now show how to get the estimate in Theorem \ref{MAIN} without
assuming the normalization condition (\ref{MR2'}). We assume that $\Sigma
_{n}:=\frac{1}{n}\sum_{k=1}^{n}\sigma _{n,k},$ is invertible and we denote $%
\overline{C}_{n,k}=\Sigma _{n}^{-1/2}C_{n,k}.$ Then we construct $\Phi
_{n,N}^{\Sigma _{n}}$ as in (\ref{M2}) by using $\overline{\Delta }%
_{n,k}(\alpha )={\mathbb{E}}((\overline{C}_{n,k}Y_{k})^{\alpha })-{\mathbb{E}%
}((\overline{C}_{n,k}G_{k})^{\alpha })$.

\begin{proposition}
\label{Normalization} Assume that $Y_{k}\in {\mathfrak{D}}(r,\varepsilon
),\forall k\in N$ for some $\varepsilon >0,r>0$ and $\Sigma _{n}=\frac{1}{n}%
\sum_{k=1}^{n}\sigma _{n,k}$ is invertible and condition (\ref{bd00}) hold.
Let $N,q\in {\mathbb{N}}$ be fixed. Then Theorem \ref{MAIN} holds as well
and (\ref{s20}) reads: for a multiindex $\alpha $ with $\left\vert \alpha
\right\vert =q$, 
\begin{equation}
\left\vert {\mathbb{E}}(\partial _{\alpha }f(S_{n}(Y)))-{\mathbb{E}}%
(\partial _{\alpha }f(\Sigma _{n}^{1/2}W)\Phi _{n,N}^{\Sigma
_{n}}(W))\right\vert \leq C\underline{\lambda }_{n}^{-q}\times {\mathrm{C}}%
_{\ast }(Y)\Big(\frac{1}{n^{\frac{1}{2}(N+1)}}L_{0}(f)+L_{q}(f)e^{-\frac{{%
\mathfrak{m}}_{\overline{r}_{n}}^{2}}{16}\times n}\Big)  \label{s20a}
\end{equation}%
where $W$ is a standard Gaussian random variable,${\mathrm{C}}_{\ast }(Y)$
is given in (\ref{s20'}) and $\underline{\lambda }_{n}$ is the lower
eigenvalue of $\Sigma _{n}$.
\end{proposition}

\textbf{Proof}. For an invertible matrix $\sigma \in \mathrm{Mat}(d\times d)$
and for $f\,:\:{\mathbb{R}}^d\to{\mathbb{R}}$, let $f_{\sigma}(x)=f(\sigma
x).$ A simple computation shows that 
\begin{equation*}
(\partial _{\alpha }f)(\sigma x)=\sum_{\left\vert \beta \right\vert
=\left\vert \alpha \right\vert }(\sigma^{-1})^{\alpha ,\beta }\partial
_{\beta }f_{\sigma }(x),
\end{equation*}%
where, for any two multiindexes $\alpha$ and $\beta$ with $%
|\alpha|=q=|\beta| $, 
\begin{equation*}
(\sigma^{-1})^{\alpha ,\beta }=\prod_{i=1}^{q}(\sigma^{-1}) ^{\alpha
_{i},\beta _{i}}.
\end{equation*}
We denote now $\overline{S}_{n}(Y)=\frac 1{\sqrt{n}}\sum_{k=1}^n\overline{C}%
_{n,k}Y_k =\Sigma _{n}^{-1/2}S_{n}(Y)$ verifies the normalization condition (%
\ref{MR2'}). So using (\ref{s20}) for $\overline{S}_{n}(Y)$ we obtain%
\begin{eqnarray*}
{\mathbb{E}}(\partial _{\alpha }f(S_{n}(Y))) &=&{\mathbb{E}}(\partial
_{\alpha }f(\Sigma _{n}^{1/2}\overline{S}_{n}(Y)))=\sum_{\left\vert \beta
\right\vert =q}(\Sigma _{n}^{-1/2})^{\alpha ,\beta }{\mathbb{E}}(\partial
_{\beta }f_{\Sigma _{n}^{1/2}}(\overline{S}_{n}(Y))) \\
&=&\sum_{\left\vert \beta \right\vert =q}(\Sigma _{n}^{-1/2})^{\alpha ,\beta
}\big({\mathbb{E}}(\partial _{\beta }f_{\Sigma _{n}^{1/2}}(W)\Phi
_{N}^{\Sigma _{n}}(W))+R_{N}^{\beta }(n)\big) \\
&=&{\mathbb{E}}(\partial _{\alpha }f(\Sigma _{n}^{1/2}W)\Phi _{N}^{\Sigma
_{n}}(W))+\sum_{\left\vert \beta \right\vert =q}(\Sigma _{n}^{-1/2})^{\alpha
,\beta }R_{N}^{\beta }(n).
\end{eqnarray*}%
The estimate of $R_{N}(n)$ follows from the fact that $L_{q}(f_{\Sigma
_{n}^{1/2}})\leq $ $\overline{\lambda }_{n}^{q}L_{q}(f)$ and $%
\sum_{\left\vert \beta \right\vert =q}(\Sigma _{n}^{-1/2})^{\alpha ,\beta
}\leq C\underline{\lambda }_{n}^{-q}\overline{\lambda}_{n}^{dq}$. $\square $

\medskip

Another immediate consequence of Theorem \ref{MAIN} is given by the
following estimate for an ``approximative density'' of the law of $S_{n}(Y)$:

\begin{proposition}
\label{cor-main} Assume that $Y_{k}\in {\mathfrak{D}}(r,\varepsilon )$ for
some $\varepsilon >0,r>0$ and let (\ref{MR2'}) and (\ref{bd00}) hold.
Suppose that $n^{\frac{1}{2}(N+1)}e^{-\frac{{\mathfrak{m}}^2_r n}{256}}\leq
1 $ and $n\geq 4(N+1){\mathrm{C}}_2(Y)$. Let $\delta _{n}$ be such that 
\begin{equation*}
n^{(N+1)/2d}e^{-\frac{{\mathfrak{m}}^2_r }{32d}\times n}\leq \delta _{n}\leq 
\frac{1}{n^{\frac{1}{2}(N+1)}}.
\end{equation*}%
Then%
\begin{equation}
\left\vert {\mathbb{E}}\Big(\frac{1}{\delta _{n}^{d}}1_{\{\left\vert
S_{n}(Y)-a\right\vert \leq \delta _{n}\}}\Big)-\gamma _{d}(a)\Phi
_{n,N}(a)\right\vert \leq \frac{C}{n^{\frac{1}{2}(N+1)}},  \label{P0}
\end{equation}
where $\gamma _{d}$ denotes the density of the standard normal law in ${%
\mathbb{R}}^{d}.$
\end{proposition}

\textbf{Proof}. Let $h(x)=\int_{-\infty }^{x_{1}}dx_{1}...\int_{-\infty
}^{x_{d-1}}\frac{1}{\delta _{n}^{d}}1_{\{\left\vert x-a\right\vert \leq
\delta _{n}\}}dx_{d}$ so that $\frac{1}{\delta _{n}^{d}}1_{\{\left\vert
x-a\right\vert \leq \delta _{n}\}}=\partial _{x_{1}}...\partial _{xd}h(x).$
Using Theorem \ref{MAIN}%
\begin{eqnarray*}
{\mathbb{E}}\Big(\frac{1}{\delta _{n}^{d}}1_{\{\left\vert
S_{n}(Y)-a\right\vert \leq \delta _{n}\}}\Big) &=&{\mathbb{E}}(\partial
_{x_{1}}...\partial _{xd}h(S_{n}(Y)))={\mathbb{E}}(\partial
_{x_{1}}...\partial _{xd}h(W)\Phi _{n,N}(W))+R_{N}(n) \\
&=&{\mathbb{E}}\Big(\frac{1}{\delta _{n}^{d}}1_{\{\left\vert W-a\right\vert
\leq \delta _{n}\}}\Phi _{n,N}(W)\Big)+R_{N}(n)
\end{eqnarray*}%
with 
\begin{equation*}
\left\vert R_{N}(n)\right\vert \leq C\Big(\frac{1}{n^{\frac{1}{2}(N+1)}}+%
\frac{1}{\delta _{n}^{d}}e^{-\frac{{\mathfrak{m}}^2_r }{32}\times n}\Big)%
\leq \frac{C}{n^{\frac{1}{2}(N+1)}}
\end{equation*}%
the last inequality being true by our choice of $\delta _{n}.$ Moreover%
\begin{eqnarray*}
{\mathbb{E}}\Big(\frac{1}{\delta _{n}^{d}}1_{\{\left\vert W-a\right\vert
\leq \delta _{n}\}}\Phi _{n,N}(W)\Big) &=&\int_{{\mathbb{R}}^{d}}\frac{1}{%
\delta _{n}^{d}}1_{\{\left\vert y-a\right\vert \leq \delta _{n}\}})\Phi
_{n,N}(y)\gamma _{d}(y)dy \\
&=&\Phi _{n,N}(a)\gamma _{d}(a)+R^{\prime }(n)
\end{eqnarray*}%
with $\left\vert R^{\prime }(n)\right\vert \leq \frac{C}{n^{\frac{1}{2}(N+1)}%
}$, as a further consequence of the choice of $\delta _{n}$. $\square $

\medskip

We now prove a stronger version of Prohorov's theorem. We consider a
sequence of identical distributed, centered random variables $X_{k}\in {%
\mathbb{R}}^{d}$ which have finite moments of any order and we look to 
\begin{equation*}
S_{n}(X)=\frac{1}{\sqrt{n}}\sum_{k=1}^{n}X_{k}.
\end{equation*}%
Following Prohorov we assume that there exist $n_\ast\in {\mathbb{N}}$ such
that%
\begin{equation}
{\mathbb{P}}(X_{1}+\cdots+X_{n_\ast}\in dx)=\mu (dx)+\psi (x)dx  \label{P1}
\end{equation}%
for some measurable non negative function $\psi .$

\begin{corollary}
\label{Prohorov} We assume that (\ref{P1}) holds. We fix $q,N\in {\mathbb{N}}%
.$ There exist two constants $0<{\mathrm{c}}_{\ast }\leq 1\leq {\mathrm{C}}%
_{\ast }$, depending on $N$ and $q$, such that the following holds: if 
\begin{equation*}
n^{\frac{1}{2}(N+1)}e^{-{\mathrm{c}}_{\ast }n}\leq 1
\end{equation*}%
then, for every multiindex $\gamma $ with $\left\vert \gamma \right\vert
\leq q$ and for every $f\in C_{p}^{q}({\mathbb{R}}^{d})$ one has%
\begin{equation}
\left\vert {\mathbb{E}}(\partial _{\gamma }f(S_{n}(X)))-{\mathbb{E}}%
(\partial _{\gamma }f(W)\Phi _{n,N}(W))\right\vert \leq {\mathrm{C}}_{\ast }%
\Big(\frac{1}{n^{\frac{1}{2}(N+1)}}L_{0}(f)+L_{q}(f)e^{-{\mathrm{c}}_{\ast
}\times n}\Big).  \label{P2}
\end{equation}
\end{corollary}

\textbf{Proof}. We denote%
\begin{equation*}
Y_{k}=\sum_{i=2kn_\ast+1}^{2(k+1)n_\ast}X_{i}\quad \mbox{and}\quad Z_{k}=%
\frac{1}{\sqrt{n}}Y_{k}.
\end{equation*}%
Notice that we may take $\psi $ in (\ref{P1}) to be bounded with compact
support. Then $\psi \ast \psi $ is continuous and so we may find some $%
r>0,\varepsilon >0$ \ and $y\in {\mathbb{R}}^{d}$\ such that $\psi \ast \psi
\geq \varepsilon 1_{B_{2r}(y)}.$ It follows that $Y_{k}\in {\mathfrak{D}}%
(r,\varepsilon )$ and we may use Theorem \ref{MAIN} in order to obtain (\ref%
{P2}) for $n=2n_\ast\times n^{\prime }$ with $n^{\prime }\in {\mathbb{N}}.$
But this is not satisfactory because we claim that (\ref{P2}) holds for
every $n\in {\mathbb{N}}.$ This does not follow directly but needs to come
back to the proof of Theorem \ref{MAIN} and to adapt it in the following
way. Suppose that $2n_\ast n^{\prime }\leq n<2n_\ast(n^{\prime }+1).$ Then%
\begin{equation*}
S_{n}(X)=S_{2n_\ast n^{\prime }}(X)+\frac{1}{\sqrt{n}}\sum_{k=2n_\ast
n^{\prime }+1}^{n}X_{k}=\frac{1}{\sqrt{n}}\sum_{k=1}^{n^{\prime }}Y_{k}+%
\frac{1}{\sqrt{n}}\sum_{k=2n_\ast n^{\prime }+1}^{n}X_{k}.
\end{equation*}%
Since $X_{k}$, $2n_\ast n^{\prime }+1\leq k\leq n$, have no regularity
property, we may not use them in the regularization arguments employed in
the proof of Theorem \ref{MAIN}. But $Y_{k}$, $1\leq k\leq n^{\prime }$
contain sufficient noise in order to achieve the proof (see the proof of
Proposition \ref{W} in next Section \ref{sect-proof-W}). $\square $

\subsection{Convergence in distribution norms}

\label{conv-distr}

In this section we prove that, under some supplementary regularity
assumptions on the laws of $Y_{k}$, $k\in {\mathbb{N}}$, Theorem \ref{MAIN}
implies that the density of the law of $S_{n}(Y)$ converges in distribution
norms to the Gaussian density. We consider the case $C_{n,k}\equiv C_{k}$,
that is, 
\begin{equation*}
S_{n}(Y)=\frac{1}{\sqrt{n}}\sum_{k=1}^{n}C_{k}Y_{k},
\end{equation*}%
and we denote $\sigma _{k}=C_{k}C_{k}^{\ast }.$ We assume that 
\begin{equation}
0<\underline{\sigma }\leq \inf_{k}\sigma _{k}\leq \sup_{k}\sigma _{k}\leq 
\overline{\sigma }<\infty \quad \mbox{and}\quad \sup_{k}\left\Vert
Y_{k}\right\Vert _{p}^{p}<\infty .  \label{D1}
\end{equation}%
In particular each $\sigma _{k}$ is invertible. We denote $\gamma
_{k}=\sigma _{k}^{-1}.$ For a function $f\in C^{1}({\mathbb{R}}^{d})$ and
for $k\in {\mathbb{N}}$\ we denote%
\begin{equation*}
m_{1,k}(f)=\int_{{\mathbb{R}}^{d}}(1+\left\vert x\right\vert )^{k}\left\vert
\nabla f(x)\right\vert dx.
\end{equation*}

\begin{proposition}
\label{D}We fix $q\in {\mathbb{N}}$ and we also fix a polynomial $P.$
Suppose that $Y_{k}\in {\mathfrak{D}}(r,\varepsilon )$, $k\in {\mathbb{N}}$,
and (\ref{D1}) holds. Suppose moreover that 
\begin{equation}
{\mathbb{P}}(Y_{k}\in dy)=p_{Y_{k}}(y)dy\quad \mbox{with}\quad p_{Y_{k}}\in
C^{1}({\mathbb{R}}^{d})\quad\mbox{for every for $i=1,...,q$}.  \label{D5}
\end{equation}

\textbf{A}. There exist some constants ${\mathrm{c}}\in (0,1)$ (depending on 
$r$ and on $\varepsilon )$ and ${\mathrm{D}}_{q}(P)\geq 1$ (depending on $q,%
\underline{\sigma},\overline{\sigma}$ and on $P)$ such that, if $%
n^{(q+1)/2}e^{-{\mathrm{c}} n}\leq 1,$ then for every $f\in C_{p}^{q}({%
\mathbb{R}}^{d})$ and every multiindex $\alpha $ with $\left\vert \alpha
\right\vert \leq q$, 
\begin{equation}
\left\vert {\mathbb{E}}(P(S_{n}(Y))\partial _{\alpha }f(S_{n}(Z))-{\mathbb{E}%
}(P(S_{n}(G))\partial _{\alpha }f(S_{n}(G))\right\vert \leq \frac{{\mathrm{D}%
}_{q}(P)}{\sqrt{n}}\prod_{i=1}^{q}m_{1,l_{0}(f)+l_{0}(P)}(p_{Y_{i}})\times
L_{0}(f).  \label{D2}
\end{equation}%
\textbf{B}. Moreover, if $p_{S_{n}}$ is the density of the law of $S_{n}(Y)$
then, if $n^{(d+q+1)/2}e^{-{\mathrm{c}} n}\leq 1,$ we have 
\begin{equation}
\sup_{x\in {\mathbb{R}}^{d}}\left\vert P(x)(\partial _{\alpha
}p_{S_{n}}(x)-\partial _{\alpha }\gamma_d (x))\right\vert \leq \frac{{%
\mathrm{D}}_{q+d}(P)}{\sqrt{n}}%
\prod_{i=1}^{q+d}m_{1,l_{0}(f)+l_{0}(P)}(p_{Y_{i}})  \label{D3}
\end{equation}%
where $\gamma_d $ is the density of the standard normal law in ${\mathbb{R}}%
^{d}$.
\end{proposition}

\textbf{Proof A.} We proceed by recurrence on the degree $i $ of the
polynomial $P$. First we assume that $i =0$ (so that $P$ is a constant$)$
and we prove (\ref{D2}) for every $q\in {\mathbb{N}}.$ We write%
\begin{equation*}
S_{n}(Y)=\frac{1}{\sqrt{n}}\sum_{k =1}^{n}C_{k }Y_{k }=\frac{1}{\sqrt{n}}%
\sum_{k =1}^{q}C_{k }Y_{k }+S_{n}^{(q)}(Y).
\end{equation*}%
with%
\begin{equation*}
S_{n}^{(q)}(Z)=\frac{1}{\sqrt{n}}\sum_{k =q+1}^{n}C_{k }Y_{k }.
\end{equation*}%
Then we define 
\begin{equation*}
g(x)={\mathbb{E}}\Big(f\Big(\frac{1}{\sqrt{n}}\sum_{k =1}^{q}C_{k }Y_{k }+x%
\Big)\Big)
\end{equation*}%
and we have 
\begin{equation*}
{\mathbb{E}}(\partial _{\alpha }f(S_{n}(Y)))={\mathbb{E}}(\partial _{\alpha
}g(S_{n}^{(q)}(Y))).
\end{equation*}%
Now using (\ref{s20a}) with $N=0$ for $S_{n}^{(q)}(Y)$ we get%
\begin{equation}
{\mathbb{E}}(\partial _{\alpha }g(S_{n}^{(q)}(Y)))={\mathbb{E}}(\partial
_{\alpha }g(S_{n}^{(q)}(G)))+R_{n}={\mathbb{E}}\Big(\partial _{\alpha }f\Big(%
\frac{1}{\sqrt{n}}\sum_{k =1}^{q}C_{k }Y_{k }+S_{n}^{(q)}(G)\Big)\Big)+R_{n}
\label{D6}
\end{equation}%
with 
\begin{equation}
\left\vert R_{n}\right\vert \leq C\Big(\frac{1}{\sqrt{n}}%
L_{0}(g)+e^{-cn}L_{q}(g)\Big).  \label{D4}
\end{equation}%
Let us estimate $L_{q}(g).$ We set $\gamma _{k }=\sigma _{k }^{-1}.$ For $%
\alpha =(\alpha _{1},...,\alpha _{q})$ we have%
\begin{equation}
(\partial _{\alpha }f)\Big(\frac{1}{\sqrt{n}}\sum_{k =1}^{n}C_{k }y_{k }+x%
\Big) =\sum_{\beta_1,\ldots,\beta_q=1}^d n^{q/2}\Big(\prod_{k =1}^{q}(\gamma
_{k }C_k )^{\alpha _{k },\beta _{k }}\Big)\times \partial _{y_{1}^{\beta
_{1}}}....\partial _{y_{q}^{\beta _{q}}}\Big( f\Big(\frac{1}{\sqrt{n}}%
\sum_{k =1}^{n}C_{k }y_{k }+x\Big)\Big),  \label{D2'}
\end{equation}%
in which we have assumed that the $Y_k $'s take values in ${\mathbb{R}}^m$.
So%
\begin{align*}
\partial _{\alpha }g(x) &={\mathbb{E}}\Big((\partial _{\alpha }f)\Big(\frac{1%
}{\sqrt{n}}\sum_{k =1}^{q}C_{k }Y_{k }+x\Big)\Big) \\
&=n^{q/2}\sum_{\beta_1,\ldots,\beta_q=1}^m\Big(\prod_{k =1}^{q}(\gamma _{k
}C_k )^{\alpha _{k },\beta _{k }}\Big)\int_{{\mathbb{R}}^{qm}}\partial
_{y_{1}^{\beta _{1}}}....\partial _{y_{q}^{\beta _{q}}}\Big( f\Big(\frac{1}{%
\sqrt{n}}\sum_{k =1}^{n}C_{k }y_{k }+x\Big)\Big) \prod_{k =1}^{q}p_{Y_{k
}}(y_{k })dy_{1}...dy_{q} \\
&=(-1)^{q}n^{q/2}\sum_{\beta_1,\ldots,\beta_q=1}^m\Big(\prod_{k=1}^{q}(%
\gamma _{k }C_k )^{\alpha _{k },\beta _{k }}\Big)\int_{{\mathbb{R}}^{qm}}f%
\Big(\frac{1}{\sqrt{n}}\sum_{k =1}^{n}C_{k }y_{k }+x\Big)\prod_{k
=1}^{q}\partial _{y_{k }^{\beta _{k }}}p_{Y_{k }}(y_{k })dy_{1}...dy_{q}.
\end{align*}%
It follows that 
\begin{eqnarray*}
\left\vert \partial _{\alpha }g(x)\right\vert &\leq &Cn^{q/2}L_{0}(f)\int_{{%
\mathbb{R}}^{q}}(1+\left\vert x\right\vert +\sum_{k =1}^{q}\left\vert y_{k
}\right\vert )^{l_{0}(f)}\prod_{k =1}^{q}\left\vert \nabla p_{Y_{k }}(y_{k
})\right\vert dy_{1}...dy_{q} \\
&\leq &Cn^{q/2}L_{0}(f)(1+\left\vert x\right\vert )^{l_{0}(f)}\prod_{k
=1}^{q}m_{1,l_{0}(f)}(p_{Y_{k }}).
\end{eqnarray*}%
We conclude that $l_{q}(g)=l_{0}(f)$ and $L_{q}(g)\leq
Cn^{q/2}L_{0}(f)\prod_{k =1}^{q}m_{1,l_{0}(f)}(p_{Y_{k }}).$ The same is
true for $q=0$ and so (\ref{D4}) gives 
\begin{equation*}
\left\vert R_{n}\right\vert \leq CL_{0}(f)\prod_{k
=1}^{q}m_{1,l_{0}(f)}(p_{Y_{k }})\Big(\frac{1}{\sqrt{n}}+n^{q/2}e^{-cn}\Big)%
\leq CL_{0}(f)\prod_{k =1}^{q}m_{1,l_{0}(f)}(p_{Y_{k }})\times \frac{1}{%
\sqrt{n}}
\end{equation*}%
the last inequality being true if $n^{q/2}e^{-cn}\leq n^{-1/2}.$

So (\ref{D6}) says that we succeed to replace $Y_{k },q+1\leq k \leq n$ by $%
G_{k },q+1\leq k \leq n$ and the price to be paid is $CL_{0}(f)\prod_{k
=1}^{q}m_{1,l_{0}(f)}(p_{Y_{k }})\times \frac{1}{\sqrt{n}}.$ Now we can do
the same thing and replace $Y_{k },1\leq k \leq q$ by $G_{k },1\leq k \leq q$
and the price will be the same (here we use $C_{k }G_{k },k =q+1,...,2q $
instead of $C_{k }Y_{k },k =1,...,q).$ So (\ref{D2}) is proved for
polynomials $P$ of degree $i =0.$

We assume now that (\ref{D2}) holds for every polynomials of degree less or
equal to $i -1$ and we prove it for a polynomial $P$ of order $i $. We have%
\begin{equation*}
\partial _{\alpha }(P\times f)=\sum_{(\beta ,\gamma )=\alpha }\partial
_{\beta }P\times \partial _{\gamma }f
\end{equation*}%
so that%
\begin{equation*}
P\times \partial _{\alpha }f=\partial _{\alpha }(P\times f)-\sum _{\substack{
(\beta ,\gamma )=\alpha  \\ \left\vert \beta \right\vert \geq 1 }}\partial
_{\beta }P_{i }\times \partial _{\gamma }f.
\end{equation*}%
Since $\left\vert \beta \right\vert \geq 1$ the polynomial $\partial _{\beta
}P$ has degree at most $i -1.$ Then the recurrence hypothesis ensures that (%
\ref{D2}) holds for $\partial _{\beta }P\times \partial _{\gamma }f.$
Moreover, using again (\ref{D2}) for $g=P\times f$ we obtain (\ref{D2}) in
which $L_{0}(g)\leq L_{0}(P)L_{0}(f)$ and $l_{0}(g)\leq l_{0}(P)+l_{0}(f)$
appear. So \textbf{A.} is proved.

\smallskip

Let us prove \textbf{B}. We denote $f_{x}(y)=\prod_{k =1}^{d}1_{(x,\infty
)}(y)$ and, for a multiindex $\alpha =(\alpha _{1},...,\alpha _{q})$ we
denote $\overline{\alpha }=(\alpha _{1},...,\alpha _{q},1,...,d).$ Then,
using a formal computation (which may de done rigorously by means of a
regularization procedure) we obtain%
\begin{eqnarray*}
P(x)\partial _{\alpha }p_{S_{n}}(x) &=&\int \delta _{0}(y-x)P(y)\partial
_{\alpha }p_{S_{n}}(y)dy \\
&=&(-1)^{q}\sum_{(\beta ,\gamma )=\alpha }\int \partial _{\beta }\delta
_{0}(y-x)\partial _{\gamma }P(y)p_{S_{n}}(y)dy \\
&=&(-1)^{q}\sum_{(\beta ,\gamma )=\alpha }\int \partial _{\overline{\beta }%
}f_{x}(y)\partial _{\gamma }P(y)p_{S_{n}}(y)dy \\
&=&(-1)^{q}\sum_{(\beta ,\gamma )=\alpha }{\mathbb{E}}(\partial _{\overline{%
\beta }}f_{x}(S_{n}(Y))\partial _{\gamma }P(S_{n}(Y))).
\end{eqnarray*}%
A similar computation holds with $S_{n}(Y)$ replaced by $S_{n}(G).$ So we
have%
\begin{eqnarray*}
&&\left\vert P(x)(\partial _{\alpha }p_{S_{n}}(x)-\partial _{\alpha }\gamma
(x)\right\vert \\
&\leq &\sum_{(\beta ,\gamma )=\alpha }\left\vert {\mathbb{E}}(\partial _{%
\overline{\beta }}f_{x}(S_{n}(Y))\partial _{\gamma }P(S_{n}(Z)))-{\mathbb{E}}%
(\partial _{\overline{\beta }}f_{x}(S_{n}(G))\partial _{\gamma
}P(S_{n}(G)))\right\vert \\
&\leq &\frac{C_{q+d}(P)}{\sqrt{n}}\prod_{k
=1}^{q+d}m_{1,l_{0}(f)+l_{0}(P)}(p_{Y_{k }})
\end{eqnarray*}%
the last inequality being a consequence of (\ref{D2}). $\square $

\begin{remark}
We would like to obtain Edgeworth's expansions as well -- but there is a
difficulty: when we use the expansion for $S_{n}^{(q)}(Z)$ we are in the
situation when the covariance matrix of $S_{n}^{(q)}(Z)$ is not the identity
matrix. So the coefficients of the expansion are computed using a correction
(see the definition of $\overline{\Delta }_{k}$ in Proposition \ref%
{Normalization}). And this correction produces an error of order $n^{-1/2}.$
This means that we are not able to go beyond this level (at least without
supplementary technical effort).
\end{remark}

\section{Examples}

\label{examples}

\subsection{An invariance principle related to the local time}

\label{sect-local}

In this section we consider a sequence of independent identically
distributed, centered random variables $Y_{k}$, $k\in {\mathbb{N}}$, with
finite moments of any order and we denote 
\begin{equation*}
S_{n}(k,Y)=\frac{1}{\sqrt{n}}\sum_{i=1}^{k}Y_{i}.
\end{equation*}%
Our aim is to study the asymptotic behavior of the expectation of 
\begin{equation*}
L_{n}(Y)=\frac 1n\sum_{k=1}^{n}\psi _{\varepsilon _{n}}(S_{n}(k,Y))\quad %
\mbox{with}\quad \psi _{\varepsilon _{n}}(x)=\frac{1}{2\varepsilon _{n}}%
1_{\{\left\vert x\right\vert \leq \varepsilon _{n}\}}.
\end{equation*}%
So $L_{n}(Y)$ appears as the occupation time of the random walk $%
S_{n}(k,Y),k=1,...,n,$ and consequently, as $\varepsilon _{n}\to 0,$ one
expects that it has to be close to the local time in zero at time 1, denoted
by $l_{1},$ of the Brownian motion. In fact, we prove now that ${\mathbb{E}}%
(L_{n}(Y))\rightarrow {\mathbb{E}}(l_{1})$ as $n\to\infty$.

\begin{theorem}
\label{OT}Let $\varepsilon _{n}=n^{-\frac{1}{2}(1-\rho )}$ with $\rho \in
(0,1).$ We consider a centered random variable $Y\in {\mathfrak{D}}%
(r,\varepsilon )$ which has finite moments of any order and we take a
sequence $Y_{i},i\in {\mathbb{N}}$ of independent copies of $Y.$ We define%
\begin{equation*}
N(Y)=\max \{2k:{\mathbb{E}}(Y^{2k})={\mathbb{E}}(G^{2k})\}-1\geq 1
\end{equation*}%
and we denote $p_{N(Y)}=8(1+(N(Y)+1)(N(Y)+3))(4+(N(Y)+1)(N(Y)+3))$. 
For every $\eta <1$ there exists a constant ${\mathrm{C}}$ depending on $%
r,\varepsilon ,\rho ,\eta $ and on $\left\Vert Y\right\Vert _{p_{N(Y)}}$
such that 
\begin{equation}  \label{OT3-1}
\big| {\mathbb{E}}(L_{n}(Y))-{\mathbb{E}}(L_n(G))\big| \leq \frac{{\mathrm{C}%
}}{n^{\frac{1}{2}+\frac{\eta \rho N(Y)}{2}}}.
\end{equation}%
The above inequality holds for $n$ which is sufficiently large in order to
have%
\begin{equation}  \label{OT1}
n^{\frac{1}{2}}\exp \Big(-\frac{{\mathfrak{m}}^2_r }{32}\times n^{\rho \eta }%
\Big) \leq \frac{1}{n^{\frac{1}{2}(N(Y)+1)\eta \rho }}
\end{equation}
As a consequence, we have 
\begin{equation}  \label{OT3}
\lim_{n\to\infty}{\mathbb{E}}(L_n(Y))={\mathbb{E}}(l_1),
\end{equation}%
$l_1$ denoting the local time in the point $0$ at time 1 of a Brownian
motion.
\end{theorem}

\textbf{Proof.} All over this proof we denote by ${\mathrm{C}}$ a constant
which depends on $r,\varepsilon ,\rho ,\eta $ and on $\left\Vert
Y\right\Vert _{p_{N(Y)}}$ (as in the statement of the lemma) and which may
change from a line to another.

\smallskip

\textbf{Step 1}. We take $k_{n}=n^{\eta \rho }.$ Suppose first that $k\leq
k_{n}.$ We write 
\begin{equation*}
{\mathbb{E}}(\psi _{\varepsilon _{n}}(S_{n}(k,Y)))=\frac{1}{\varepsilon _{n}}%
\Big(1-{\mathbb{P}}(\left\vert S_{n}(k,Y)\right\vert \geq \varepsilon _{n})%
\Big)
\end{equation*}%
so that%
\begin{equation*}
\left\vert {\mathbb{E}}(\psi _{\varepsilon _{n}}(S_{n}(k,Y)))-{\mathbb{E}}%
(\psi _{\varepsilon _{n}}(S_{n}(k,G)))\right\vert \leq \frac{1}{\varepsilon
_{n}}\bigl({\mathbb{P}}(\left\vert S_{n}(k,Y)\right\vert \geq \varepsilon
_{n})+{\mathbb{P}}(\left\vert S_{n}(k,G)\right\vert \geq \varepsilon _{n})%
\bigr).
\end{equation*}%
Using Chebyshev's inequality and Burkholder's inequality we obtain for every 
$p\geq 2$ 
\begin{align*}
{\mathbb{P}}(\left\vert S_{n}(k,Y)\right\vert \geq \varepsilon _{n}) &={%
\mathbb{P}}\Big(\Big\vert \sum_{i=1}^{k}Y_{i}\Big\vert \geq \varepsilon _{n}%
\sqrt{n}\Big) \leq \frac{1}{(\varepsilon _{n}\sqrt{n})^{p}}{\mathbb{E}}\Big(%
\Big\vert \sum_{i=1}^{k}Y_{i}\Big\vert ^{p}\Big) \\
\leq &\frac{{\mathrm{C}}}{(\varepsilon _{n}\sqrt{n})^{p}}\Big(%
\sum_{i=1}^{k}\Vert Y_{i}\Vert _{p}^{2}\Big)^{p/2} \leq \frac{{\mathrm{C}}
k^{p/2}}{(\varepsilon _{n}\sqrt{n})^{p}} =\frac{{\mathrm{C}}}{\varepsilon
_{n}^{p}}\times \Big(\frac{k}{n}\Big)^{p/2}.
\end{align*}%
And the same estimate holds with $Y_{i}$ replaced by $G_{i}.$ We conclude
that

\begin{align*}
&\Big\vert {\mathbb{E}}\Big(\frac{1}{n}\sum_{k=1}^{k_{n}}\psi _{\varepsilon
_{n}}(S_{n}(k,Y))\Big)-{\mathbb{E}}\Big(\frac{1}{n}\sum_{k=1}^{k_{n}}\psi
_{\varepsilon _{n}}(S_{n}(k,G)) \Big)\Big\vert \leq \frac{{\mathrm{C}}}{%
\varepsilon _{n}^{p+1}}\times \frac 1n\sum_{k=1}^{k_{n}}\Big(\frac{k}{n}\Big)%
^{p/2} \\
&\leq \frac{{\mathrm{C}}}{\varepsilon _{n}^{p+1}}\times
\int_{0}^{k_{n}/n}x^{p/2}dx=\frac{{\mathrm{C}}}{\varepsilon _{n}^{p+1}}%
\times \Big(\frac{k_{n}}{n}\Big)^{\frac{p}{2}+1} \\
&=\frac{{\mathrm{C}}}{n^{\frac{p\rho }{2}(1-\eta )+\frac{1}{2}-(\eta -\frac{1%
}{2})\rho }}\leq \frac{{\mathrm{C}}}{n^{\frac{p\rho }{2}(1-\eta )}}.
\end{align*}%
We take $p=\frac{1+\rho \eta N(Y)}{\rho (1-\eta )}$ and we obtain 
\begin{equation*}
\Big\vert {\mathbb{E}}\Big(\frac{1}{n}\sum_{k=1}^{k_{n}}\psi _{\varepsilon
_{n}}(S_{n}(k,Y))\Big)-{\mathbb{E}}\Big(\frac{1}{n}\sum_{k=1}^{k_{n}}\psi
_{\varepsilon _{n}}(S_{n}(k,G)) \Big)\Big\vert \leq \frac {\mathrm{C}}{%
n^{\frac 12+\frac{N(Y)}2\eta\rho}}.
\end{equation*}

\textbf{Step 2}. We fix now $k\geq k_{n}$ and we apply our Edgeworth
development (\ref{s20}) to 
\begin{equation*}
\frac{1}{\sqrt{k}}\sum_{i=1}^{k}Y_{i}.
\end{equation*}
In particular the constants ${\mathrm{C}}_{p}(Y)$ defined in (\ref{bd00})
are given by ${\mathrm{C}}_{p}(Y)=\| Y\|_{p}^{p}.$ We denote%
\begin{equation}  \label{s-1}
h_{\alpha ,n}(x)=\int_{-\infty }^{\alpha x}\psi _{\varepsilon
_{n}}(y)dy=h_{1,n}(\alpha x).
\end{equation}%
This gives $\psi _{\varepsilon _{n}}(x)=h_{1,n}^{\prime }(x)$ and $h_{\alpha
,n}^{\prime }(x)=\alpha h_{1,n}^{\prime }(\alpha x).$ Moreover, $\|
h_{\alpha,n}\|_{\infty }\leq 1$ and $\| h_{\alpha,n}^{\prime }\| _{\infty
}\leq |\alpha|/\varepsilon _{n}$, so that%
\begin{equation*}
L_{0}(h_{\alpha,n})=1\quad\mbox{and}\quad L_{1}(h_{\alpha,n})=
|\alpha|\times \frac{1}{\varepsilon _{n}}.
\end{equation*}%
We now write%
\begin{eqnarray*}
{\mathbb{E}}(\psi _{\varepsilon _{n}}(S_{n}(k,Y))) &=&{\mathbb{E}}%
(h_{1,n}^{\prime }(S_{n}(k,Y))) ={\mathbb{E}}\Big(h_{1,n}^{\prime }\Big(%
\sqrt{\frac{k}{n}}\frac{1}{\sqrt{k}}\sum_{i=1}^{k}Y_{i}\Big)\Big) \\
&=&\sqrt{\frac{n}{k}}{\mathbb{E}}\Big(h_{\sqrt{\frac{k}{n}},n}^{\prime }\Big(%
\frac{1}{\sqrt{k}}\sum_{i=1}^{k}Y_{i}\Big)\Big).
\end{eqnarray*}%
We use now (\ref{s20}) with $f=h_{\sqrt{\frac{k}{n}},n}$ and here $\partial
_{\gamma }$ is the first order derivative. Then, by (\ref{s20}) with $N=N(Y)$%
\begin{equation*}
{\mathbb{E}}(\psi _{\varepsilon _{n}}(S_{n}(k,Y))=\sqrt{\frac{n}{k}}\Big({%
\mathbb{E}}\big(h_{\sqrt{\frac{k}{n}},n}^{\prime }(W_{1})\Phi
_{k,N(Y)}(W_{1})\big)+R_{N(Y)}(k)\Big)
\end{equation*}%
where $W$ denotes a Brownian motion and with 
\begin{eqnarray*}
\left\vert R_{N(Y)}(k)\right\vert &\leq &\frac{{\mathrm{C}}}{k^{(N(Y)+1)/2}}%
L_0(h_{\sqrt{\frac{k}{n}},n})+{\mathrm{C}} L_1(h_{\sqrt{\frac{k}{n}}%
,n}^{\prime })\exp \Big(-\frac{{\mathfrak{m}}^2_r }{32}\times k\Big) \\
&\leq &\frac{{\mathrm{C}}}{k^{(N(Y)+1)/2}}+{\mathrm{C}}\sqrt{\frac{k}{n}}%
\times \frac{1}{\varepsilon _{n}}\exp \Big(-\frac{{\mathfrak{m}}^2_r }{32}%
\times k\Big).
\end{eqnarray*}%
Here ${\mathrm{C}}$ is the constant from (\ref{s20}) defined in (\ref{s20'}%
). Notice that by (\ref{OT1}), for $k\geq k_{n}=n^{\eta\rho}$ one has 
\begin{eqnarray*}
\sqrt{\frac{k}{n}}\times \frac{1}{\varepsilon _{n}}\exp \Big(-\frac{{%
\mathfrak{m}}^2_r }{32}\times k\Big) &\leq &n^{\frac{1}{2}}\exp \Big(-\frac{{%
\mathfrak{m}}^2_r }{32}\times n^{\rho \eta }\Big) \\
&\leq &\frac{1}{n^{\frac{1}{2}(N(Y)+1)\eta \rho }}=\frac{1}{%
k_{n}^{(N(Y)+1)/2}}\leq \frac{{\mathrm{C}}}{k^{(N(Y)+1)/2}},
\end{eqnarray*}%
so that $\left\vert R_{N(Y)}(k)\right\vert \leq {\mathrm{C}}
k^{-(N(Y)+1)/2}. $ Then%
\begin{align*}
&\Big\vert \sum_{k=k_{n}}^{n}\sqrt{\frac{n}{k}}R_{N(Y)}(k)\frac{1}{n}%
\Big\vert \leq \frac{{\mathrm{C}}}{n^{(N(Y)+1)/2}}\sum_{k=k_{n}}^{n}\frac{1}{%
(k/n)^{1+\frac{N(Y)}{2}}}\times \frac{1}{n} \\
&\quad \leq \frac{{\mathrm{C}}}{n^{(N(Y)+1)/2}}\int_{k_{n}/n}^{1}\frac{ds}{%
s^{1+\frac{N(Y)}{2}}} =\frac{{\mathrm{C}}}{n^{(N(Y)+1)/2}}(n/k_{n})^{\frac{%
N(Y)}{2}}=\frac{{\mathrm{C}}}{n^{\frac{1}{2}+\frac{N(Y)\rho \eta }{2}}}.
\end{align*}%
We recall now that (see (\ref{M5})) 
\begin{equation*}
\Phi _{k,N(Y)}(x)=1+\sum_{l=1}^{N(Y)}\frac 1{n^{l/2}}H_{\Gamma _{k,l}}(x)
\end{equation*}%
with $H_{\Gamma _{k,l}}(x)$ linear combination of Hermite polynomials (see (%
\ref{M2}) and (\ref{M4})). Notice that if $l$ is odd then $\Gamma _{k,l}$ is
a linear combination of differential operators of odd order (see the
definition of $\Lambda _{m,l}$ in (\ref{M1})). So $H_{\Gamma _{k,l}}$ is an
odd function (as a linear combination of Hermite polynomials of odd order)
so that $\psi _{\varepsilon _{n}}\times H_{\Gamma _{k,l}}$ is also an odd
function. Since $W_{1}$ and $-W_{1}$ have the same law, it follows that 
\begin{eqnarray*}
{\mathbb{E}}\Big(\psi _{\varepsilon _{n}}\Big(\sqrt{\frac{k}{n}}\times W_{1}%
\Big)H_{k,\Gamma _{l}}(W_{1})\Big) &=&{\mathbb{E}}\Big(\psi _{\varepsilon
_{n}}\Big(\sqrt{\frac{k}{n}}\times (-W_{1})\Big)H_{\Gamma _{k,l}}(-W_{1})%
\Big) \\
&=&-{\mathbb{E}}\Big(\psi _{\varepsilon _{n}}\Big(\sqrt{\frac{k}{n}}\times
W_{1}\Big)H_{\Gamma _{k,l}}(W_{1})\Big)
\end{eqnarray*}%
and consequently 
\begin{equation*}
\sqrt{\frac{n}{k}}\times {\mathbb{E}}\Big(h_{\sqrt{\frac{k}{n}},n}^{\prime
}(W_{1})H_{\Gamma _{k,l}}(W_{1})\Big)={\mathbb{E}}\Big(\psi _{\varepsilon
_{n}}\Big(\sqrt{\frac{k}{n}}\times W_{1}\Big)H_{\Gamma _{k,l}}(W_{1})\Big)=0.
\end{equation*}%
Moreover, by the definition of $N(Y)$, for $2l\leq N(Y)$ we have ${\mathbb{E}%
}(Y^{2l})={\mathbb{E}}(G^{2l})$ so that $H_{\Gamma _{k,2l}}=0$. We conclude
that 
\begin{equation*}
\sqrt{\frac{n}{k}}{\mathbb{E}}\Big(h_{\sqrt{\frac{k}{n}},n}^{\prime
}(W_{1})\Phi _{k,N(Y)}(W_{1})\Big)=\sqrt{\frac{n}{k}}{\mathbb{E}}\Big(h_{%
\sqrt{\frac{k}{n}},n}^{\prime }(W_{1})\Big)={\mathbb{E}}\Big(\psi
_{\varepsilon _{n}}\Big(\sqrt{\frac{k}{n}}\times W_{1}\Big)\Big)={\mathbb{E}}%
(\psi _{\varepsilon _{n}}(S_{n}(k,G))).
\end{equation*}
We put now together the results from the first and the second step and we
obtain (\ref{OT3-1}).

\smallskip

\textbf{Step 3}. We prove (\ref{OT3}). Recall first the representation
formula%
\begin{equation*}
{\mathbb{E}}\Big(\int_{0}^{1}\psi _{\varepsilon _{n}}(W_{s})ds\Big)={\mathbb{%
E}}\Big(\int \psi _{\varepsilon _{n}}(a)l_{1}^{a}da\Big),
\end{equation*}
where $l_1^a$ denotes the local time in $a\in{\mathbb{R}}$ at time 1, so
that $l_1=l_1^0$. Since $a\mapsto l_{1}^{a}$ is H\"{o}lder continuous of
order $\frac{\rho ^{\prime }}{2}$ for every $\rho ^{\prime }<1,$ we obtain 
\begin{equation}
\Big\vert {\mathbb{E}}\Big(\int_{0}^{1}\psi _{\varepsilon _{n}}(W_{s})ds\Big)%
-{\mathbb{E}}(l_{1}^{0})\Big\vert \leq \varepsilon _{n}^{\rho ^{\prime }/2}=%
\frac{1}{n^{\frac{\rho ^{\prime }(1-\rho )}{4}}}.  \label{OT4}
\end{equation}%
We prove now that, for every $\rho ^{\prime }<1$ and $n$ large enough, 
\begin{equation}
\Big\vert {\mathbb{E}}\Big(\int_{0}^{1}\psi _{\varepsilon _{n}}(W_{s})ds\Big)%
-{\mathbb{E}}(L_{n}(G))\Big\vert \leq \frac{C}{n^{\frac{1+\eta\rho}2}}.
\label{OT5}
\end{equation}%
To begin we notice that $S_{n}(k,G)$ has the same law as $W_{k/n}$, so that
we write 
\begin{equation*}
{\mathbb{E}}\Big(\int_{0}^{1}\psi _{\varepsilon _{n}}(W_{s})ds\Big)-{\mathbb{%
E}}(L_{n}(G))={\mathbb{E}}\Big(\sum_{k=1}^{n}\delta _{k}\Big),\quad%
\mbox{with}\quad \delta _{k}=\int_{k/n}^{(k+1)/n}(\psi _{\varepsilon
_{n}}(W_{s})-\psi _{\varepsilon _{n}}(W_{k/n}))ds.
\end{equation*}%
As above, we take $k_{n}=n^{\rho \eta }$ and for $k\leq k_{n},$ we have 
\begin{equation*}
{\mathbb{E}}(\delta _{k})=-\frac{1}{2\varepsilon _{n}}\int_{k/n}^{(k+1)/n}%
\big({\mathbb{P}}(|W_{s}|\geq \varepsilon _{n})-{\mathbb{P}}(|W_{k/n}|\geq
\varepsilon _{n})\big)ds.
\end{equation*}%
Since ${\mathbb{P}}(|W_{s}|\geq \varepsilon _{n})\leq C\exp (-\frac{%
\varepsilon ^{2}}{2s})$, this immediately gives%
\begin{equation*}
\vert {\mathbb{E}}(\delta _{k})\vert \leq \frac{C}{n\varepsilon _{n}}\exp %
\Big(-\frac{1}{2}\varepsilon _{n}^{2}\times \frac{n}{k+1}\Big) \leq \frac{C}{%
n\varepsilon _{n}}\exp \Big(-\frac{1}{2}\varepsilon _{n}^{2}\times \frac{n}{%
k_{n}+1}\Big) =\frac{C}{n\varepsilon _{n}}\exp \Big(-\frac{1}{2}n^{\rho
(1-\eta )}\Big)
\end{equation*}%
so that%
\begin{equation*}
\sum_{k=1}^{k_{n}}\vert {\mathbb{E}}(\delta _{k})\vert \leq \frac{C}{%
\varepsilon _{n}}\exp (-\frac{1}{8}n^{\rho (1-\eta )})\leq \frac{C}{n^{\frac{%
1+\eta\rho}2}},
\end{equation*}
for $n$ large enough.

We consider now the case $k\geq k_{n}.$ Using a formal computation, by
applying the standard Gaussian integration by parts formula, we write 
\begin{align*}
&{\mathbb{E}}(\psi _{\varepsilon _{n}}(W_{s})-\psi _{\varepsilon
_{n}}(W_{k/n})) =\frac{1}{2}\int_{k/n}^{s}{\mathbb{E}}(\psi _{\varepsilon
_{n}}^{\prime \prime }(W_{v}))dv=\frac{1}{2}\int_{k/n}^{s}{\mathbb{E}}(\psi
_{\varepsilon _{n}}^{\prime \prime }(\sqrt{v}W_{1}))dv \\
&=\int_{k/n}^{s}{\mathbb{E}}(h^{\prime \prime \prime }_{1,n}(\sqrt{v}%
W_{1})H_{3}(W_{1}))dv =\int_{k/n}^{s}\frac{1}{2v^{3/2}}{\mathbb{E}}(h_{1,{n}%
}(\sqrt{v}W_{1})H_{3}(W_{1}))dv,
\end{align*}%
in which we have used (\ref{s-1}) and where $H_{3}$ denotes the third
Hermite polynomial. The above computation is formal because $\psi
_{\varepsilon _{n}}$ is not differentiable. But, since the first and the
last term in the chain of equalities depends on $\psi _{\varepsilon _{n}}$
only (and not on the derivatives) we may use regularization by convolution
in order to do it rigorously. Notice also that the first equality is
obtained using Ito's formula and the last one is obtained using integration
by parts. It follows that%
\begin{equation*}
\vert {\mathbb{E}}(\delta _{k})\vert \leq
\int_{k/n}^{(k+1)/n}ds\int_{k/n}^{s}\frac{1}{2v^{3/2}}{\mathbb{E}}%
(h_{1,\varepsilon _{n}}(\sqrt{v}W_{1})\left\vert H_{3}(W_{1})\right\vert
)dv\leq \frac{C}{n}\int_{k/n}^{(k+1)/n}\frac{1}{v^{3/2}}dv
\end{equation*}%
and consequently%
\begin{equation*}
\sum_{k=k_{n}}^{n}\left\vert {\mathbb{E}}(\delta _{k})\right\vert \leq \frac{%
C}{n}\int_{k_{n}/n}^{1}\frac{1}{v^{3/2}}dv\leq \frac{C}{n^{\frac{1+\eta\rho}{%
2}}}.
\end{equation*}%
So (\ref{OT5}) is proved, and this together with (\ref{OT4}) and (\ref{OT3-1}%
), give (\ref{OT3}). $\square $

\subsection{Small ball estimates}

\label{small-balls}

We look to 
\begin{equation}
S_{n}(u,Y)=\frac{1}{\sqrt{n}}\sum_{k=1}^{n}C_{n,k}(u)Y_{k},\quad u\in{%
\mathbb{R}}^\ell ,  \label{1}
\end{equation}%
where $Y_{k}\in {\mathbb{R}}^{d}$, $k\in {\mathbb{N}}$, and $C_{n,k}(u)\in 
\mathrm{Mat}(d\times d)$ (so, here $m=d$).

\begin{theorem}
\label{SB} Suppose that $\{Y_{k}\}_{k\in {\mathbb{N}}}\subset {\mathfrak{D}}%
(\varepsilon ,r)$, with ${\mathrm{M}}_{p}(Y)=\sup_{k}\Vert Y_{k}\Vert
_{p}<\infty $, and that $u\mapsto C_{n,k}^{i,j}(u)$ is twice differentiable.
We assume that for every $n\in {\mathbb{N}},k\leq n$ and $u\in {\mathbb{R}}%
^{\ell }$ 
\begin{align}
& \Vert C_{n,k}\Vert _{2,\infty }:=\sum_{i,j=1}^{d}\sum_{\left\vert \alpha
\right\vert \leq 2}\Vert \partial _{u}^{\alpha }C_{n,k}^{i,j}\Vert
_{2,\infty }\leq Q_{\ast ,2}<\infty ,  \label{2} \\
& \frac{1}{n}\sum_{k=1}^{n}C_{n,k}(u)C_{n,k}^{\ast }(u)\geq \lambda _{\ast
}>0,  \label{3}
\end{align}

\medskip

\textbf{A}. There exist ${\mathrm{C}}\geq 1$ and ${\mathrm{c}}>0$ such that
for every $\eta >0$%
\begin{equation}
\sup_{u\in {\mathbb{R}}^{\ell }}{\mathbb{P}}(\left\vert
S_{n}(u,Y)\right\vert \leq \eta )\leq {\mathrm{C}}(\eta ^{d}+e^{-\mathrm{c}%
n}).  \label{4a}
\end{equation}%
\textbf{B}. Suppose that $d>\ell .$ Let $a\geq 0$ and $\theta >\frac{a\ell }{%
d-\ell }.$ Then, for every $\varepsilon >0$ 
\begin{equation}
{\mathbb{P}}(\inf_{\left\vert u\right\vert \leq n^{a}}\left\vert
S_{n}(u,Y)\right\vert \leq \frac{1}{n^{\theta }})\leq \frac{{\mathrm{C}}}{%
n^{\theta (d-\ell )-a\ell -\varepsilon }}  \label{4}
\end{equation}%
The constant ${\mathrm{C}}$ depends on ${\mathfrak{m}}_r$ (from Doeblin's
condition), on $Q_{\ast ,2},\lambda _{\ast },\ell ,d$ and on ${\mathrm{M}}%
_p(Y)$ for sufficiently large $p.$
\end{theorem}

We first prove the following lemma.

\begin{lemma}
\label{lemma-balls} Under the hypotheses of Theorem \ref{SB}, for every $%
q>\ell $, $i\in \{1,\ldots,d\}$ and $R>0$ one has 
\begin{equation}
{\mathbb{E}}(\sup_{\left\vert u\right\vert \leq R}\left\vert \partial
_{i}S_{n}(u,Y)\right\vert ^{q})\leq {\mathrm{C}} R^{\ell }Q_{\ast ,2}^{q}{%
\mathrm{M}}_{q}^{q}(Y).  \label{5}
\end{equation}%
where ${\mathrm{C}}$ is a constant which depends on $q.$
\end{lemma}

\textbf{Proof}. As an immediate consequence of Morrey's inequality one may
find a universal constant ${\mathrm{C}}$ (independent of $R)$ such that 
\begin{equation*}
\sup_{\left\vert u\right\vert \leq R}\left\vert \partial
_{i}S_{n}(u,Y)\right\vert \leq {\mathrm{C}}\Big(\int_{\left\vert
u\right\vert \leq R+1}\left\vert \partial _{i}S_{n}(u,Y)\right\vert
^{q}+\sum_{j=1}^{\ell }\left\vert \partial _{j}\partial
_{i}S_{n}(u,Y)\right\vert ^{q}du\Big)^{1/q}
\end{equation*}%
so that%
\begin{equation*}
{\mathbb{E}}(\sup_{\left\vert u\right\vert \leq R}\left\vert \partial
_{i}S_{n}(u,Y)\right\vert ^{q})\leq {\mathrm{C}}\int_{\left\vert
u\right\vert \leq R+1}\Big({\mathbb{E}}\left\vert \partial
_{i}S_{n}(u,Y)\right\vert ^{q}+\sum_{j=1}^{\ell }{\mathbb{E}}\left\vert
\partial _{j}\partial _{i}S_{n}(u,Y)\right\vert ^{q}\Big)du.
\end{equation*}%
Since 
\begin{equation*}
\partial _{j}\partial _{i}S_{n}(u,Y)=\frac{1}{\sqrt{n}}\sum_{k=1}^{n}%
\partial _{j}\partial _{i}C_{n,k}(u)Y_{k},
\end{equation*}%
we can use the Burkholder's inequality for martingales and we obtain%
\begin{eqnarray*}
{\mathbb{E}}\left\vert \partial _{j}\partial _{i}S_{n}(u,Y)\right\vert ^{q}
&\leq &{\mathrm{C}} {\mathbb{E}}\Big(\Big[\frac{1}{n}\sum_{k=1}^{n}\left%
\Vert \partial _{j}\partial _{i}C_{n,k}(u)(\partial_{i}\partial_{\alpha
}C_{n,k}(u))^{\ast }\right\Vert ^{2}\left\vert Y_{k}\right\vert ^{2}\Big]%
^{q/2}\Big) \\
&\leq &{\mathrm{C}} Q_{\ast ,2}^{q}{\mathbb{E}}\Big(\Big[\frac{1}{n}%
\sum_{k=1}^{n}\left\vert Y_{k}\right\vert ^{2}\Big]^{q/2}\Big)\leq {\mathrm{C%
}} Q_{\ast ,2}^{q}M_{q}^{q}(Y).
\end{eqnarray*}%
A similar estimate holds for ${\mathbb{E}}\left\vert \partial
_{i}S_{n}(u,Y)\right\vert ^{q},$ so that (\ref{5}) is proved. $\square $

\bigskip

\textbf{Proof of Theorem} \ref{SB}. \textbf{A}. Let us prove (\ref{4a}). We
take $\eta >0$ and we consider the functions%
\begin{equation}
\theta _{d,\eta }(x)=\frac{1}{(c_{d}\eta )^{d}}1_{\left\vert x\right\vert
\leq \eta },\quad \Theta _{d,\eta }(x)=\int_{-\infty
}^{x_{1}}dx_{2}...\int_{-\infty }^{x_{d-1}}dx_{d}\theta _{d,\eta }(x)
\label{4b}
\end{equation}%
with $c_{d}$ such that $\int_{{\mathbb{R}}^{d}}\theta _{d,\eta }(x)dx=1.$
Then $\partial _{1}....\partial _{d}\Theta _{d,\eta }=\theta _{d,\eta }$ so
that 
\begin{equation*}
{\mathbb{P}}(\left\vert S_{n}(u,Y)\right\vert \leq \eta )=(c_{d}\eta )^{d}{%
\mathbb{E}}(\theta _{d,\eta }(S_{n}(u,Y))=(c_{d}\eta )^{d}{\mathbb{E}}%
(\partial _{1}....\partial _{d}\Theta _{d,\eta }(S_{n}(u,Y)).
\end{equation*}%
We denote $S_{n}(t,G)$ the sum from (\ref{1}) in which $Y_{k},k\in {\mathbb{N%
}} ,$ are replaced by standard normal random variables and we use Theorem %
\ref{MAIN}, specifically (\ref{s20}), in order to obtain 
\begin{equation*}
{\mathbb{E}}(\partial _{1}....\partial _{d}\Theta _{d,\eta }(S_{n}(u,Y))={%
\mathbb{E}}(\partial _{1}....\partial _{d}\Theta _{d,\eta
}(S(u,G))+\varepsilon _{n}(\eta )
\end{equation*}%
where 
\begin{equation*}
\left\vert \varepsilon _{n}(\eta )\right\vert \leq {\mathrm{C}}\Big(\frac{1}{%
n^{1/2}}+\eta ^{-d}e^{-{\mathrm{c}} n}\Big).
\end{equation*}%
Here ${\mathrm{C}}$ is a constant which depends on ${\mathfrak{m}}_r$ (from
Doeblin's condition) on $Q_{\ast ,0}$ and on ${\mathrm{M}}_{p}(Y)$ for a
sufficiently large $p.$ We conclude that%
\begin{align*}
{\mathbb{P}}(\left\vert S_{n}(u,Y)\right\vert \leq \eta ) &=(c_{d}\eta )^{d}{%
\mathbb{E}}(\partial _{1}....\partial _{d}\Theta _{d,\eta }(S_{n}(u,Y))) \\
&\leq (c_{d}\eta )^{d}\big({\mathbb{E}}(\partial _{1}....\partial _{d}\Theta
_{d,\eta }(S(u,G)))+\left\vert \varepsilon _{n}(\eta )\right\vert\big) \\
&={\mathrm{C}} {\mathbb{P}}(\left\vert S(u,G)\right\vert \leq \eta )+{%
\mathrm{C}} \eta ^{d}\left\vert \varepsilon _{n}(\eta )\right\vert .
\end{align*}%
Since $S_{n}(u,G)$ is a non degenerate Gaussian random variable we have ${%
\mathbb{P}}(\vert S_{n}(u,G)\vert \leq \eta )\leq {\mathrm{C}}\eta
^{d}\lambda _{\ast }^{-d/2}$ and finally we get 
\begin{equation*}
{\mathbb{P}}(\left\vert S_{n}(t_{\alpha },Y)\right\vert \leq \eta )\leq {%
\mathrm{C}}\eta ^{d}(\lambda _{\ast }^{-d/2}+\left\vert \varepsilon
_{n}(\eta )\right\vert )\leq {\mathrm{C}}(\eta ^{d}+e^{-cn}).
\end{equation*}

\medskip

\textbf{B}. We denote $R_{n}=n^{a},$ $\delta _{n}=n^{-\theta }$ and we take $%
h>0$ (to be chosen later on). For $\alpha \in {\mathbb{Z}}^{\ell }$ we
denote $t_{\alpha }=(t_{\alpha _{1}},...,t_{\alpha _{\ell }})=(h\alpha
_{1},\ldots,h\alpha _{\ell })$ and $I_{\alpha }=[t_{\alpha _{1}},t_{\alpha
_{1}+1})\times\cdots\times [t_{\alpha _{\ell }},t_{\alpha _{\ell }+1}),$ so,
if $\left\vert u\right\vert \leq R_{n}$ then $u\in \cup _{\left\vert
t_{\alpha }\right\vert \leq R_{n}}I_{\alpha }.$ Moreover we denote%
\begin{equation*}
\omega _{n}=\inf_{\left\vert u\right\vert \leq R_{n}}\left\vert
S_{n}(u,Y)\right\vert ,\quad \omega _{n,\alpha }=\inf_{u\in I_{\alpha
}}\left\vert S_{n}(u,Y)\right\vert
\end{equation*}%
and we have 
\begin{equation*}
\omega _{n}\geq \min_{\left\vert t_{\alpha }\right\vert \leq R_{n}}\omega
_{n,\alpha }.
\end{equation*}%
If $\omega _{n,\alpha }<\delta _{n}$ then there is some $u_{\alpha }\in
I_{\alpha }$ such that $\left\vert S_{n}(u_{\alpha })\right\vert \leq \delta
_{n}$. So, with $U_{n}=\sup_{\left\vert u\right\vert \leq R_{n}}\left\vert
\nabla S_{n}(u,Y)\right\vert ,$ we have%
\begin{equation*}
\left\vert S_{n}(t_{\alpha })\right\vert \leq \left\vert S_{n}(u_{\alpha
})\right\vert +U_{n}h\leq \delta _{n}+U_{n}h.
\end{equation*}%
Now we take $\lambda >0$ (to be chosen later on) and we write, with $q>\ell $%
,%
\begin{align*}
{\mathbb{P}}(\omega _{n} \leq \delta _{n}) &\leq {\mathbb{P}}(\omega
_{n}\leq \delta _{n},U_{n}\leq \lambda )+{\mathbb{P}}(U_{n}\geq \lambda ) \\
&\leq \sum_{\left\vert t_{\alpha }\right\vert \leq R_{n}}{\mathbb{P}}(\omega
_{n,\alpha }\leq \delta _{n},U_{n}\leq \lambda )+\lambda ^{-q}{\mathbb{E}}%
(U_{n}^{q}) \\
&\leq (R_{n}/h)^{\ell }\max_{\left\vert t_{\alpha }\right\vert \leq R_{n}}{%
\mathbb{P}}(\left\vert S_{n}(t_{\alpha },Y)\right\vert \leq \delta
_{n}+\lambda h)+{\mathrm{C}} \lambda ^{-q}R_{n}^{\ell }Q_{\ast
,2}^{q}M_{q}^{q}(Y) \\
&\leq {\mathrm{C}}(R_{n}/h)^{\ell }((\delta _{n}+\lambda h)^{d}+e^{-cn})+{%
\mathrm{C}}\lambda ^{-q}R_{n}^{\ell }Q_{\ast ,2}^{q}M_{q}^{q}(Y),
\end{align*}%
in which we have used (\ref{5}) and (\ref{4a}). We recall that $R_{n}=n^{a}$
and $\delta _{n}=n^{-\theta }.$ We take $\lambda =n^{\varepsilon }$ for a
sufficiently small $\varepsilon >0$ and $h=n^{-(\theta +\varepsilon )}.$
Then, for large enough $q,$ we get%
\begin{equation*}
{\mathbb{P}}(\omega _{n}\leq \delta _{n})\leq Cn^{(a+\theta +\varepsilon
)\ell }\times n^{-\theta d}+Cn^{-\varepsilon q}\times n^{\ell a}\leq
Cn^{-(\theta d-(a+\theta +\varepsilon )\ell )}.
\end{equation*}%
$\square $

\subsection{Expected number of roots for trigonometric polynomials: an
invariance principle}

\label{ROOTS}

In this section we look to trigonometric polynomials with random
coefficients of the form 
\begin{equation*}
Q_{n}(t,Y)=\sum_{k=1}^{n}\big(Y_{k}^{1}\cos (kt)+Y_{k}^{2}\sin (kt)\big)
\end{equation*}%
where $Y_{k}=(Y_{k}^{1},Y_{k}^{2}),k\in {\mathbb{N}},$ are independent
centered random variables such that $Y_k\in{\mathfrak{D}}(\varepsilon, r)$
for each $k$. Our aim is to estimate the asymptotic behavior, as $%
n\rightarrow \infty $, of the expected number of zeros in the interval $%
(0,\pi )$ of these polynomials. This clearly coincide with the number of
zeros in $(0,n\pi )$ of the renormalized polynomials%
\begin{equation}
P_{n}(t,Y)=\frac{1}{\sqrt{n}}\sum_{k=1}^{n}\Big(Y_{k}^{1}\cos \Big(\frac{kt}{%
n}\Big)+Y_{k}^{2}\sin \Big(\frac{kt}{n}\Big)\Big).  \label{p1}
\end{equation}%
So we denote by $N_{n}(Y)$ the number of zeros of $P_{n}(t,Y)$ in $(0,n\pi )$%
. It is known that if we replace $Y_{k}$ by $G_{k},$ independent standard
normal random variables then (see \cite{D,F}) 
\begin{equation*}
\lim_{n}\frac{1}{n}{\mathbb{E}}(N_{n}(G))=\frac 1{\sqrt 3}.
\end{equation*}%
Our aim is to prove that this remains true for any sequence $Y_{k},k\in {%
\mathbb{N}}$ of independent but non necessarily identically distributed
random variables. So we will prove:

\begin{theorem}
\label{roots} Suppose that $Y=(Y_{k})_{k\in {\mathbb{N}}}$ is a sequence of
independent random variables in $\subset{\mathfrak{D}}(\varepsilon,r)$,
having finite moments of any order. Then 
\begin{equation}
\Big\vert \frac{1}{n}{\mathbb{E}}(N_{n}(Y))-\frac{1}{n}{\mathbb{E}}(N_{n}(G))%
\Big\vert \leq \frac{{\mathrm{C}}}{\sqrt{n}}.  \label{p2}
\end{equation}
\end{theorem}

\textbf{Proof}. The first ingredient in the proof is Kac-Rice lemma that we
recall now. Let $f:[a,b]\rightarrow {\mathbb{R}}$ be a differentiable
function and set 
\begin{equation*}
\mbox{ $\omega _{a,b}(f)=\inf_{x\in \lbrack a,b]}(\left\vert
f(x)\right\vert +\left\vert f^{\prime }(x)\right\vert )$ and $\delta
_{a,b}(f)=\min \{\left\vert f(a)\right\vert ,\left\vert f(b)\right\vert
,\omega _{a,b}(f)\}.$}
\end{equation*}
We denote by $N_{a,b}(f)$ the number of solutions of $f(t)=0$ for $t\in
\lbrack a,b]$ and 
\begin{equation*}
I_{a,b}(f,\delta)=\int_{a}^{b}\left\vert f^{\prime }(t)\right\vert
1_{\{\left\vert f(t)\right\vert \leq \delta \}}\frac{dt}{2\delta }%
,\quad\delta>0.
\end{equation*}
The Kac-Rice lemma says that if $\delta_{a,b}(f)>0$ then 
\begin{equation}
N_{a,b}(f)=I_{a,b}(f,\delta) \quad \mbox{for}\quad \delta \leq\delta
_{a,b}(f).  \label{p3}
\end{equation}%
Notice that we also have, for every $\delta >0$,%
\begin{equation}
I_{a,b}(\delta ,f)\leq 1+N_{a,b}(f^{\prime }).  \label{p3'}
\end{equation}%
Indeed, we may assume that $N_{a,b}(f^{\prime })=p<\infty $ and then we take 
$a=a_{0}\leq a_{1}<....<a_{p}\leq a_{p+1}=b$ to be the roots of $f^{\prime
}. $ Since $f$ is monotonic on each $(a_{i},a_{i+1})$ one has $%
I_{a_{i},a_{i+1}}(\delta ,f)\leq 1$ so (\ref{p3'}) holds.

We will use this result for $f(t)=P_{n}(t,Y)$ so we have $N_{n}(Y)=N_{0,n\pi
}(P_{n}(t,Y)).$ We denote $\delta _{n}(Y)=\delta _{0,n\pi }(P_{n}(t,Y))),$
we take $\theta =3$ and we write%
\begin{align*}
\frac{1}{n}{\mathbb{E}}(N_{n}(Y)) = &\frac{1}{n}{\mathbb{E}}%
(N_{n}(Y)1_{\{\delta _{n}(Y)\leq n^{-\theta }\}})-\frac{1}{n}{\mathbb{E}}%
(I_{0,n\pi }(\delta ,P_{n}(\cdot,Y))1_{\{\delta _{n}(Y)\leq n^{-\theta }\}})+%
\frac{1}{n}{\mathbb{E}}(I_{0,n\pi }(\delta ,P_{n}(\cdot,Y))) \\
=:&A_{n}(Y)-A_{n}^{\prime }(Y)+B_{n}(Y).
\end{align*}%
A trigonometric polynomial of order $n$ has at most $2n$ roots on $(0,\pi ).$
So the number of roots of $P_{n}(t,Y)$ on $(0,\pi )$ is upper bonded by $2n$%
, so that consequently $N_{n}(Y)\leq 2n.$ It follows that $A_{n}(Y)\leq 2{%
\mathbb{P}}(\delta _{n}(Y)\leq n^{-\theta }).$ Since $P_{n}^{\prime }$ is
also a trigonometric polynomial of order $n$, by (\ref{p3'}) we also have $%
I_{0,n\pi }(\delta ,P_{n}(.,Y))\leq 1+N_{0,n\pi }(P_{n}^{\prime }(t,Y))\leq
2n+1.$ It follows that $\left\vert A_{n}^{\prime }(Y)\right\vert \leq 3{%
\mathbb{P}}(\delta _{n}(Y)\leq n^{-\theta }).$

We will use Theorem \ref{SB} and Theorem \ref{MAIN} for $%
S_{n}(t,Y)=(P_{n}(t,Y),P_{n}^{\prime }(t,Y))$, so we have to check the
hypotheses there. Notice that in this case we have $\ell=1,d=2$ and 
\begin{equation*}
C_{n,k}(t)=\left( 
\begin{tabular}{ll}
$\cos (\frac{kt}{n})$ & $\sin (\frac{kt}{n})$ \\ 
$-\frac{k}{n}\sin (\frac{kt}{n})$ & $\frac{k}{n}\cos (\frac{kt}{n})$%
\end{tabular}%
\right) .
\end{equation*}%
First, (\ref{2}) trivially holds. Moreover, for every $\xi \in {\mathbb{R}}%
^{2}$ one has $\vert C_{n,k}(t)\xi \vert ^{2}=\frac{k^{2}}{n^{2}}$ so that 
\begin{equation*}
\frac{1}{n}\sum_{k=1}^{n}\left\vert C_{n,k}(t)\xi \right\vert ^{2}=\frac{1}{n%
}\sum_{k=1}^{n}\frac{k^{2}}{n^{2}}\left\vert \xi \right\vert ^{2}\geq
\int_{0}^{1}x^{2}dx\times \left\vert \xi \right\vert ^{2}=\frac{1}{3}%
\left\vert \xi \right\vert ^{2}.
\end{equation*}%
This means that (\ref{3}) holds with $\lambda _{\ast }=\frac{1}{3}$ and we
are able to use (\ref{4}) in order to get $A_{n}(Y)\leq {\mathrm{C}}/n$ and $%
\left\vert A_{n}^{\prime }(Y)\right\vert \leq {\mathrm{C}}/n$. Moreover, by (%
\ref{p3}) 
\begin{equation*}
B_{n}(Y)=\frac{1}{n}{\mathbb{E}}\Big(\int_{0}^{n\pi }\left\vert
P_{n}^{\prime }(t,Y)\right\vert 1_{\{\left\vert P_{n}(t,Y)\right\vert \leq
\delta _{n}\}}\frac{dt}{\delta _{n}}\Big)\quad \mbox{with}\quad \delta _{n}=%
\frac{1}{n^{\theta }}.
\end{equation*}%
We now use the Theorem \ref{MAIN} applied to $\Psi _{\delta
_{n}}(x_{1},x_{2})=\left\vert x_{2}\right\vert \Theta _{1,\delta
_{n}}(x_{1}) $ with $\Theta _{1,\delta _{n}}$ defined in (\ref{4b}). Then%
\begin{equation*}
B_{n}(Y)=\frac{1}{n}{\mathbb{E}}\Big(\int_{0}^{n\pi }\partial _{1}\Psi
_{\delta _{n}}(S_{n}(t,Y))dt\Big)\quad \mbox{with}\quad \delta _{n}=\frac{1}{%
n^{\theta }}.
\end{equation*}%
We have $\left\Vert \Psi _{\delta _{n}}\right\Vert _{1,\infty }\leq \delta
_{n}^{-1}$ so, using (\ref{s20}) we get 
\begin{equation*}
\left\vert {\mathbb{E}}\big(\partial _{1}\Psi _{\delta _{n}}(S_{n}(t,Y))\big)%
-{\mathbb{E}}\big(\partial _{1}\Psi _{\delta _{n}}(S_{n}(t,G))\big)%
\right\vert \leq C\Big(\frac{1}{\sqrt{n}}+n^{3}e^{-cn}\Big)
\end{equation*}%
and this gives $\left\vert B_{n}(Y)-B_{n}(G)\right\vert \leq Cn^{-1/2}.$ As
above we have $A_{n}(G)\leq Cn^{-1}$ and $\left\vert A_{n}^{\prime
}(G)\right\vert \leq Cn^{-1}$ so we finally obtain (\ref{p2}). $\square $


\section{The case of smooth test functions}

\label{smooth}

We first study a variant of our main Theorem \ref{MAIN}, namely, we assume
that $q=1$ therein and we ask for a smooth function $f$. In this case,
thanks to the regularity assumption for $f$, we do not need any Doeblin's
condition. This will be used in a second step, where we will be able to
relax the smoothness assumption for $f$ by means of a regularization result
from Malliavin calculus.

We come back to the notation introduced in Section \ref{main-results}. We
just recall here the corrector polynomial $\Phi _{n,N}$ defined in (\ref{M5}):%
\begin{equation*}
\Phi _{n,N}(x)=1+\sum_{k=1}^{N}\frac{1}{n^{k/2}}H_{\Gamma _{n,k}}(x),
\end{equation*}%
where $H_{\Gamma _{n,k}}$ is the Hermite polynomial associated with the
differential operator $\Gamma _{n,k}$ defined in (\ref{M2}).

The result we prove in this section is the following:

\begin{theorem}
\label{Smooth}Let $N\in {\mathbb{N}}$ be given. 
Suppose that the normalization property (\ref{MR2'}) and the moment bounds (%
\ref{bd00}) both hold (the latter being sufficient for $p\leq N+3$). Then
for every $f\in C_{p}^{2N(\lfloor N/2\rfloor+N+5)}({\mathbb{R}}^{d})$ 
\begin{equation}  \label{M6}
\begin{array}{l}
\displaystyle \left\vert {\mathbb{E}}(f(S_{n}(Y))-{\mathbb{E}}(f(W)\Phi
_{n,N}(W))\right\vert\smallskip \\ 
\displaystyle \leq \mathcal{H}_N{\mathrm{C}}_{2(N+3)}^{2N(N+2\lfloor
N/2\rfloor)}(Y)(1+{\mathrm{C}}_{2l_{\widehat{N}}(f)}(Y))^{2N+3} 2^{(N+2)(l_{%
\widehat{N}}(f)+1)}L_{\widehat{N}}(f)\times\frac 1{n^{\frac{N+1}2}}%
\end{array}%
\end{equation}%
in which $\widehat{N}=N(2\lfloor N/2\rfloor+N+5)$, $\mathcal{H}_N$ is a
positive constant depending on $N$ and $W$ denotes a standard normal random
variable in ${\mathbb{R}}^d$. As a consequence, taking $f(x)=x^\beta$ with $%
|\beta|=k$, one gets 
\begin{equation}  \label{M6'}
|{\mathbb{E}}(S_{n}(Y)^\beta)-{\mathbb{E}}(W^\beta\Phi _{N}(W))| \leq 
\mathcal{H}_N{\mathrm{C}}_{2(N+3)}^{2N(N+2\lfloor N/2\rfloor)}(Y)(1+{\mathrm{%
C}}_{2k}(Y))^{2N+3} 2^{(N+2)(k+1)}\times\frac 1{n^{\frac{N+1}2}}.
\end{equation}
\end{theorem}

\bigskip In order to give the proof of Theorem \ref{Smooth}, we introduce a
decomposition allowing us to work with suitable semigroups. But first, in
order to simplify the forthcoming notation, we set 
\begin{equation*}
Z_{n,k}=\frac{1}{\sqrt{n}}C_{n,k}Y_{k}\quad \mbox{and}\quad G_{n,k}=\frac{1}{%
\sqrt{n}}C_{n,k}G_{k},
\end{equation*}%
so that 
\begin{equation}
S_{n}(Y)=\sum_{k=1}^{n}Z_{n,k}=:\overline{S}_{n}(Z)\quad \mbox{and}\quad
S_{n}(G)=\sum_{k=1}^{n}G_{n,k}=:\overline{S}_{n}(G).  \label{Sbar}
\end{equation}%
Notice that the covariance matrices of $Z_{n,k}$ and $G_{n,k}$ are both
given by 
\begin{equation*}
{\mathrm{Cov}}(Z_{n,k})={\mathrm{Cov}}(G_{n,k})=\overline{\sigma }_{n,k}=%
\frac{1}{n}\sigma _{n,k},
\end{equation*}%
so the normalization condition (\ref{MR2'}) reads%
\begin{equation*}
\sum_{k=1}^{n}\overline{\sigma }_{n,k}=\mathrm{Id}_{d}.
\end{equation*}

\textbf{Sketch of the proof}. The proof of the above theorem is rather long
and technical, so, in order to orient the reader, we give first a sketch of
it. The strategy is based on the classical Lindeberg method but it turns
out that it is convenient to do it in terms of semigroups (the so called
Trotter's method). We define the Markov semigroup%
\begin{equation}
P_{k,p}^{Z,n}f(x)=E(f(x+\sum_{i=k}^{p-1}Z_{n,i}))  \label{skech1}
\end{equation}%
with the convention $P_{k,k}^{Z,n}f=f.$ Then Lindberg's decomposition
gives 
\begin{equation}
P_{k,n+1}^{Z,n}-P_{k,n+1}^{G,n}=%
\sum_{r=k}^{n}P_{r+1,n+1}^{Z,n}(P_{r,r+1}^{Z,n}-P_{r,r+1}^{G,n})P_{k,r}^{G,n}.
\label{skech2}
\end{equation}%
We use now Taylor expansion of order three. The terms of order one and two
cancel (because the moments of order one and two of $Y_{r}$ and $G_{r}$
coincide) and we obtain%
\begin{eqnarray}
\delta _{n,r}f(x)
&:&=(P_{r,r+1}^{Z,n}-P_{r,r+1}^{G,n})f(x)=E(f(x+Z_{n,r})-E(f(x+G_{n,r})
\label{skech3} \\
&=&\frac{1}{6}\sum_{\left\vert \alpha \right\vert =3}\int_{0}^{1}E(\partial
^{\alpha }f(\lambda Z_{n,k}+(1-\lambda )G_{n,k})(Z_{n,k}^{\alpha
}-G_{n,k}^{\alpha })d\lambda  \label{skech3a}
\end{eqnarray}%
so one obtains $\left\Vert \delta _{n,r}f(x)\right\Vert _{\infty }\leq
C\left\Vert f\right\Vert _{3,\infty }\frac{1}{n^{3/2}}.$ We insert this in (%
\ref{skech2}) and we obtain $P_{k,n+1}^{Z,n}-P_{k,n+1}^{G,n}\sim n\times 
\frac{1}{n^{3/2}}=\frac{1}{n^{1/2}}.$ This is the proof of the classical $%
CLT.$ Now, if we want to obtain Edgeworth development of order $N$, we have to
go further. First we iterate (\ref{skech2}) and we obtain 
\begin{equation*}
P_{k,n+1}^{Z,n}f=P_{k,n+1}^{G,n}f+\sum_{l=1}^{N}%
\sum_{1<r_{1}<...<r_{l}<n}^{n}T_{r_{1},...,r_{l}}f+R_{n}^{N}f
\end{equation*}%
with%
\begin{eqnarray*}
T_{r_{1},...,r_{l}}f &=&P_{r_{l}+1,n+1}^{G,n}\prod_{\hat{\imath}%
=1}^{l-1}\left( \delta _{n,r_{i}}P_{r_{i-1},r_{i}}^{G,n}\right) f\quad and \\
R_{n}^{N}f &=&\sum_{1<r_{1}<...<r_{N+1}<n}^{n}P_{r_{N+1}+1,n+1}^{Z,n}\prod_{%
\hat{\imath}=1}^{N}\left( \delta _{n,r_{i}}P_{r_{i-1},r_{i}}^{G,n}f\right) .
\end{eqnarray*}
Since each of $\delta _{n,r}f$ is of order $n^{-3/2}$ it follows that $%
R_{n}^{N}f$ is of order $n^{N+1}\times n^{-\frac{3}{2}(N+1)}=n^{-(N+1)/2}.$

We look now to $T_{r_{1},...,r_{l}}f.$ We expect this term to be of order $%
n^{-\frac{3}{2}l},$ and indeed, it is. But we notice that $\delta _{n,r_{i}}$
contains information on the whole law of $Y,$ and not only on its moments
(see (\ref{skech3})). So, if we want to obtain the real coefficient of order 
$n^{-l/2}$ in the Edgeworth development, we have to replace $\delta _{n,r}$ by
some $\widehat{\delta }_{n,r}$ which depends on the moments only - and this
is done by using Taylor expansion as in (\ref{skech3}). But, as we want to
obtain a final error of order $n^{-(N+1)/2},$ the development of order three
is no more sufficient and we need now a development of order $N+2.$ This
will involve differential operators of order less or equal to $N+2$ with
coefficients computed by using the difference of moments given in (\ref{bd1}%
). Here is that the Hermite polynomials come on, due to the following
integration by parts formula: if $G$ is a standard normal random variable
then $E(\partial ^{\alpha }f(G))=E(f(G)H_{\alpha }(G))$ where $\partial
^{\alpha }$ is the differential operator associated to the multi-index $%
\alpha $ and $H_{\alpha }$ is the Hermite polynomial corresponding to $\alpha
.$ Collecting the terms of order $n^{-l/2}$ form all the $%
T_{r_{1},...,r_{l}}^{\prime }{}^{\prime }s$ we get the corrector of order $l$
in the Edgeworth expansion. These are the main ideas of the proof. However,
the precise description of the coefficients of the Edgeworth expansion, turns
out to be a very technical matter. We do all this through the following
lemmas in this section.

We go no on and give the complete proofs. Let $N\in \{0,1,\ldots \}.$ For $%
D_{n,r}$ given in (\ref{bd1}), we define%
\begin{equation}
T_{n,N,r}^{0}f(x)=\sum_{l=1}^{N+2}\frac{1}{n^{l/2}\,l!}D_{n,r}^{(l)}f(x).
\label{bd2}
\end{equation}%
Since $D_{n,r}^{(l)}\equiv 0$ for $l=0,1,2$, the above sum actually begins
with $l=3$ and of course this is the basic fact. Then, with the convention $%
\sum_{l=3}^{2}=0,$ we have 
\begin{equation*}
T_{n,N,r}^{0}f(x)=\sum_{l=3}^{N+2}\frac{1}{n^{l/2}\,l!}D_{n,r}^{(l)}f(x).
\end{equation*}

We also define%
\begin{equation}  \label{bd3}
\begin{array}{rcl}
T_{n,N,r}^{1,Z}f(x) & = & \displaystyle \frac{1}{(N+2)!}\sum_{\left\vert
\alpha \right\vert =N+3}\int_{0}^{1}(1-\lambda) ^{N+2}{\mathbb{E}}%
(\partial_{\alpha }f(x+\lambda Z_{n,r})Z_{n,r}^{\alpha })d\lambda \quad %
\mbox{and} \smallskip \\ 
T_{n,N,r}^{1}f(x) & = & \displaystyle %
T_{n,N,r}^{1,Z}f(x)-T_{n,N,r}^{1,G}f(x).%
\end{array}%
\end{equation}
For a matrix $\sigma \in \mathrm{Mat}({d\times d})$ we recall the Laplace
operator $L_{\sigma }$ associated to $\sigma $ (see (\ref{bd0a})) and we
define 
\begin{align}
&h_{N,\sigma }^{0}f(x)=f(x)+\sum_{l=1}^{\lfloor N/2\rfloor}\frac{(-1)^{l}}{%
n^l\,2^{l}l!}L_{\sigma }^{l}f(x),  \label{bd4} \\
&h_{N,\sigma }^{1}f(x)=\frac{(-1)^{\lfloor N/2\rfloor+1}}{n^{\lfloor
N/2\rfloor+1}\,2^{\lfloor N/2\rfloor+1}\lfloor N/2\rfloor!}%
\int_{0}^{1}s^{\lfloor N/2\rfloor}{\mathbb{E}}(L_{\sigma }^{\lfloor
N/2\rfloor+1}f(x+\sigma ^{1/2}\sqrt s\,W))ds.  \label{bd5}
\end{align}%
In (\ref{bd5}), $W$ stands for a standard Gaussian random variable. Then we
define%
\begin{equation}
U_{n,N,r}^{0}f(x)={\mathbb{E}}(h_{N,\sigma _{n,r}}^{0}f(x+G_{n,r}))\quad%
\mbox{and}\quad U_{n,N,r}^{1}f(x)=h_{N,\sigma _{n,r}}^{1}f(x).  \label{bd9}
\end{equation}

We now put our problem in a semigroup framework. For a sequence $X_k$, $%
k\geq 1$, of independent r.v.'s, for $1\leq k\leq p$ we define 
\begin{equation}  \label{bd6}
\mbox{$P_{k,k}^Xf(x)=f$ and for $p>k\geq 1$ then $P_{k,p}^{X}f(x)=
\displaystyle \E\Big(f\Big(x+\sum_{i=k}^{p-1}X_{i}\Big)\Big)$.}
\end{equation}%
By using independence, we have the semigroup and the commutative property: 
\begin{equation}
P_{k,p}^X=P_{r,p}^XP_{k,r}^X=P_{k,r}^XP_{r,p}^X\quad k\leq r\leq p.
\label{bd6bis}
\end{equation}
We use $P_{k,p}^{X}$ with $X_k=Z_{n,k}$ and $X_k=G_{n,k}$, that we call $%
P_{k,p}^{Z,n}$ and $P_{k,p}^{G,n}$ because each local random variables
depend on $n$.

Moreover, for $m=1,\ldots,N$ we denote 
\begin{equation}  \label{bd8}
\begin{array}{rcl}
Q_{n,N,r_{1},...,r_{m}}^{(m)} & = & \displaystyle \sum_{%
\mbox{\scriptsize{
$\begin{array}{c}
\sum_{i=1}^mq_i+\sum_{i=1}^mq'_i>0\\
q_{i},q_{i}^{\prime }\in\{0,1\}
\end{array}$}}} \prod_{i=1}^{m}U_{n,N,r_{i}}^{q_{i}^{\prime }}
\prod_{j=1}^{m}T_{n,N,r_{j}}^{q_{j}}\quad \mbox{and}\smallskip \\ 
R_{n,N,k}^{(m)} & = & \displaystyle \sum_{k\leq r_{1}<\cdots <r_{m}\leq
n}P_{r_{m}+1,n}^{G,n}P_{r_{m-1}+1,r_{m}}^{G,n}\cdots
P_{r_{1}+1,r_{2}}^{G,n}P_{k,r_{1}}^{G,n}Q_{n,N,r_{1},...,r_{m}}^{(m)}.%
\end{array}%
\end{equation}%
Notice that in the first sum above the conditions $q_{i},q_{i}^{\prime }\in
\{0,1\}$ and $q_{1}+\cdots +q_{m}+q_{1}^{\prime }+\cdots +q_{m}^{\prime }>0$
say that at least one of $q_{i},q_{i}^{\prime },i=1,...,m$ is equal to one.
We notice that the operators $T_{n,N,r_{i}}^{1}$ and $U_{n,N,\sigma
_{r_{i}}}^{1}$ represent ``remainders'' and they are supposed to give small
quantities of order $n^{-\frac{1}{2}(N+1)}$. So the fact that at least one $%
q_{i}$ or $q_{i}^{\prime }$ is non null means that the product $%
(\prod_{i=1}^{m}U_{n,N,r_{i}}^{q_{i}^{\prime
}})(\prod_{i=1}^{m}T_{n,N,r_{i}}^{q_{i}})$ has at least one term which is a
remainder (so is small), and consequently $R_{n,N,k}^{(m)}$ is a remainder
also.

Finally we define%
\begin{equation}  \label{bd8''}
\begin{array}{rl}
Q_{n,N,r_{1},...,r_{N+1}}^{(N+1)} & = \displaystyle %
\prod_{i=1}^{N+1}(T_{N,r_{i}}^{0}+T_{N,r_{i}}^{1}) \quad\mbox{and}\smallskip
\\ 
R_{n,N,k}^{(N+1)} & = \displaystyle \sum_{k\leq r_{1}<\cdots <r_{N+1}\leq
n}P_{r_{N+1}+1,n}^{Z,n}P_{r_{N}+1,r_{N+1}}^{G,n}\cdots
P_{r_{1}+1,r_{2}}^{G,n}P_{k,r_{1}}^{G,n}Q_{n,N,r_{1},...,r_{N+1}}^{(N+1)}%
\end{array}%
\end{equation}

As a preliminary result for Theorem \ref{Smooth}, we study the following
``backward Taylor formula'':

\begin{lemma}
\label{lemma-h4} Let $\mathcal{N}_{k}$, $k\in{\mathbb{N}}$, denote
independent centered Gaussian random variables in ${\mathbb{R}}^d$ with
covariance matrix $\sigma _{k}$, $k\in{\mathbb{N}}$, and set $%
S_{p}=\sum_{k=1}^{p}\mathcal{N}_{k}.$ For $\sigma\in\mathrm{Mat}(d\times d)$%
, we define 
\begin{align*}
H_{N,\sigma }^{0}\phi (x)=\phi (x)+\sum_{l=1}^{N}\frac{(-1)^{l}}{2^{l}l!}%
L_{\sigma }^{l}\phi (x)\mbox{ and } H_{N,\sigma }^{1}\phi (x)=\frac{%
(-1)^{N+1}}{2^{N+1}N!}\int_{0}^{1}s^{N}{\mathbb{E}}(L_{\sigma }^{N+1}\phi
(x+\sigma^{1/2}W_{s}))ds.
\end{align*}%
where $W$ is a $d-$dimensional Brownian motion independent of $S_{p}.$ Then
for every $\phi\in C^{2N+2}({\mathbb{R}}^d)$ one has 
\begin{equation}
{\mathbb{E}}(\phi (S_{p}))={\mathbb{E}}(H_{N,\sigma _{p+1}}^{0}\phi
(S_{p+1}))+{\mathbb{E}}(H_{N,\sigma _{p+1}}^{1}\phi (S_{p}))  \label{h4}
\end{equation}
\end{lemma}

\textbf{Proof.} We use the following property: for every $N\in{\mathbb{N}}$, 
$N\geq 0$, and $g\in C_b^{2N+2}({\mathbb{R}}^d)$ one has 
\begin{equation}  \label{h2}
g(0)={\mathbb{E}}(\sigma^{1/2} W_1))+\sum_{l=1}^N \frac{(-1)^l}{2^\ell l!}{%
\mathbb{E}}(L_\sigma^l g(\sigma^{1/2}W_1)) + \frac{(-1)^{N+1}}{2^{N+1} N!}%
\int_0^1 s^N {\mathbb{E}}(L^{N+1}_\sigma g(\sigma^{1/2}W_s))ds,
\end{equation}
in which $W$ denotes a standard Brownian motion in ${\mathbb{R}}^d$ and $%
L_\sigma$ is given in (\ref{bd0a}). The decomposition (\ref{h2}) is proved
in \cite{[CLT]} (see Appendix C therein) in the case $\sigma=\mathrm{Id}$
and (\ref{h2}) represents a straightforward generalization to any covariance
matrix $\sigma$.

We notice that $\mathcal{N}_{p+1}$ has the same law as $%
\sigma^{1/2}_{p+1}W_{1}$. We denote $\psi _{\omega }(x)=\phi (S_{p}(\omega
)+x).$ Then, using the independence property and (\ref{h2}) we obtain 
\begin{eqnarray*}
{\mathbb{E}}(\psi _{\omega }(0)) &=&{\mathbb{E}}(\psi _{\omega
}(\sigma^{1/2}_{p+1}W_{1}))+\sum_{l=1}^{N}\frac{(-1)^{l}}{2^{l}l!}{\mathbb{E}%
}(L_{\sigma _{p+1}}^{l}\psi _{\omega }(\sigma^{1/2}_{p+1}W_{1})) \\
&&+\frac{(-1)^{N+1}}{2^{N+1}N!}\int_{0}^{1}s^{N}{\mathbb{E}}(L_{\sigma
_{p+1}}^{N+1}\psi _{\omega }(\sigma^{1/2}_{p+1}W_{s})))ds.
\end{eqnarray*}%
Since ${\mathbb{E}}(L_{\sigma _{p+1}}^{l}\psi _{\omega
}(\sigma^{1/2}_{p+1}W_{1}))={\mathbb{E}}(L_{\sigma _{p+1}}^{l}\phi
(S_{p+1})) $ and ${\mathbb{E}}(L_{\sigma _{p+1}}^{N}\psi _{\omega }
(\sigma^{1/2}_{p+1}W_{1})) ={\mathbb{E}}(L_{\sigma_{p+1}}^{N+1}\phi
(S_{p}+\sigma^{1/2}_{p+1}W_{s}))$ the above formula is exactly (\ref{h4}). $%
\square $

\medskip

We are now able to give our first result:

\begin{proposition}
\label{prop1} Let $N\geq 1$ and let $T_{n,N,r}^{0}$, $h_{n,N,\sigma
_{r}}^{0} $ and $R^{(m)}_{n,N,k}$, $m=1,\ldots,N+1$, be given through (\ref%
{bd2}), (\ref{bd4}) and (\ref{bd8})--(\ref{bd8''}). Then for every $1\leq
k\leq n+1$ and $f\in C_{p}^{N(2\lfloor N/2\rfloor+N+5)}({\mathbb{R}}^{d})$
one has 
\begin{equation}  \label{bd7}
P_{k,n+1}^{Z,n}f =P_{k,n+1}^{G,n}f+\sum_{m=1}^{n}\sum_{k\leq r_{1}<\cdots
<r_{m}\leq n}P_{k,n+1}^{G,n}\Big(\prod_{i=1}^{m}T_{n,N,r_{i}}^{0}\Big) \Big(%
\prod_{j=1}^{m}h_{N,\sigma _{n,r_{j}}}^{0}\Big)f +
\sum_{m=1}^{N+1}R_{n,N,k}^{(m)}f.
\end{equation}
\end{proposition}

\textbf{Proof}. \textbf{Step 1 (Lindeberg method).} We use the Lindeberg
method in terms of semigroups: for $1\leq k\leq n+1$%
\begin{equation*}
P_{k,n+1}^{Z,n}-P_{k,n+1}^{G,n}=%
\sum_{r=k}^{n}P_{r+1,n+1}^{Z,n}(P_{r,r+1}^{Z,n}-P_{r,r+1}^{G,n})P_{k,r}^{G,n}.
\end{equation*}%
Then we define 
\begin{equation}  \label{Akp}
A_{k,p}= 1_{1\leq k\leq p-1\leq
n}\,(P_{p-1,p}^{Z,n}-P_{p-1,p}^{G,n})P_{k,p-1}^{G,n}
\end{equation}
(here $n$ is fixed so we do not stress the dependence of $A_{k,p}$ on $n$)
and the above relation reads%
\begin{equation}
P_{k,n+1}^{Z,n}=P_{k,n+1}^{G,n}+\sum_{r=k}^{n}P_{r+1,n+1}^{Z,n}A_{k,r+1}.
\label{BD2}
\end{equation}
We will write (\ref{BD2}) as a discrete time Volterra type equation (this is
inspired from the approach to the parametrix method given in \cite{[BH]}:
see equation (3.1) there). For a family of operators $F_{k,p}$, $k\leq p$ we
define $AF$ by 
\begin{equation*}
(AF)_{k,p}=\sum_{r=k}^{p-1}F_{r+1,p}A_{k,r+1}
\end{equation*}%
and we write (\ref{BD2}) in functional form:%
\begin{equation}
P^{Z,n}=P^{G,n}+AP^{Z,n}.  \label{BD3}
\end{equation}%
By iteration, 
\begin{equation}
P^{Z,n}=P^{G,n}+AP^{G,n}+\cdots +A^{N}P^{G,n}+A^{N+1}P^{Z,n}.  \label{BD4}
\end{equation}
By the commutative property in (\ref{bd6bis}), straightforward computations
give 
\begin{equation}  \label{AmF}
\begin{array}{l}
(A^mP^{G,n})_{k,p} =\displaystyle 1_{k\leq p-m} \!\!\!\!\!\!\!\!\!\!\!\!\!
\sum_{k\leq r_1<\cdots<r_{m}\leq p-2} \!\!\!\!\!\!\!\!\!\!\!\!\!
P^{G,n}_{r_{m}+1,p-1} P^{G,n}_{r_{m-1}+1,r_m} \cdots P^{G,n}_{r_{1}+1,r_2}
P^{G,n}_{k,r_{1}}(P^{Z,n}_{p-1,p}-P^{G,n}_{p-1,p})\times\smallskip \\ 
\displaystyle \qquad\times (P^{Z,n}_{r_{m},r_{m}+1}-P^{G,n}_{r_{m},r_{m}+1})
(P^{Z,n}_{r_{m-1},r_{m-1}+1}-P^{G,n}_{r_{m-1},r_{m-1}+1}) \cdots
(P^{Z,n}_{r_{1},r_{1}+1}-P^{G,n}_{r_{1},r_{1}+1}).%
\end{array}%
\end{equation}


\textbf{Step 2 (Taylor formula).} The drawback of (\ref{BD4}) is that $A$
depends on $P^{Z,n}$ also, see (\ref{Akp}), so we use the Taylor's formula
in order to eliminate this dependence We take into account (\ref{bd00}) and
we consider a Taylor approximation at the level of an error of order $n^{-%
\frac{N+2}{2}}$. We use the following expression for the Taylor's formula:
for $f\in C^{\infty}_p({\mathbb{R}}^d)$, 
\begin{equation*}
f(x+y)=f(x)+\sum_{p=1}^{N+2}\frac{1}{p!}\sum_{\left\vert \alpha \right\vert
=p}\partial _{\alpha }f(x)y^{\alpha }+ \frac 1{(N+2)!}\sum_{\left\vert
\alpha \right\vert =N+3}y^{\alpha}\int_0^1(1-\lambda)^{N+2}\partial _{\alpha
}f(x+\lambda y)d\lambda
\end{equation*}
Then we have, with $D_{n,r}^{(l)}$ defined in (\ref{bd1}), 
\begin{align*}
&(P_{r,r+1}^{Z,n}-P_{r,r+1}^{G,n})f(x) ={\mathbb{E}}(f(x+Z_{n,r}))-{\mathbb{E%
}}(f(x+G_{n,r})) =\sum_{l=1}^{N+2}\frac{1}{l!}D_{r}^{(l)}f(x)+ \\
&\qquad +\frac{1}{(N+2)!}\sum_{\left\vert \alpha \right\vert
=N+3}\int_{0}^{1}(1-\lambda)^{N+2}\big[{\mathbb{E}}(\partial_{\alpha
}f(x+\lambda Z_{n,r})Z_{n,r}^{\alpha })-{\mathbb{E}}(\partial_{\alpha
}f(x+\lambda G_{n,r})G_{n,r}^{\alpha })\big]d\lambda \\
&=T_{n,N,r}^{0}f(x)+ T_{n,N,r}^{1}f(x).
\end{align*}%
By using the independence property, one can apply commutativity and by using
(\ref{AmF}) we have 
\begin{equation}
(A^{m}F)_{k,r+1}=1_{k\leq r+1-m} \!\!\!\!\sum_{k\leq r_{1}<\cdots <r_{m}\leq
r} \!\!\!\!F_{r_{m}+1,r+1}P_{r_{m-1}+1,r_{m}}^{G,n}\cdots
P_{r_{1}+1,r_{2}}^{G,n}P_{k,r_{1}}^{G,n} \prod_{j=1}^m
(T_{n,N,r_{j}}^0+T_{n,N,r_{j}}^1).  \label{bd10}
\end{equation}%
Notice that the operator in (\ref{bd10}) acts on $f\in C^{m(N+3)}$. 

\medskip

\textbf{Step 3 (Backward Taylor formula).} Since 
\begin{equation*}
P^{G,n}_{r_{m}+1,n+1}P^{G,n}_{r_{m-1}+1,r_m}\cdots
P^{G,n}_{r_{1}+1,r_{2}}P^{G,n}_{k,r_{1}}f(x) ={\mathbb{E}}\Big(f\Big(%
x+\sum_{i=k}^n G_{n,k}-\sum_{j=1}^m G_{n,r_j}\Big)\Big),
\end{equation*}
the chain $P_{r_{m}+1,n}^{G,n}\cdots
P_{r_{1}+1,r_{2}}^{G,n}P_{k,r_{1}}^{G,n} $ contains all the steps, except
for the steps corresponding to $r_{i},i=1,...,m$ (remark that for each $i$, $%
P_{r_{i},r_{i}+1}^{G,n}$ is replaced with $%
T_{n,N,r_{i}}^{0}+T_{n,N,r_{i}}^{1}$). In order to ``insert'' such steps we
use the backward Taylor formula (\ref{h4}) up to order $\overline{N}=\lfloor
N/2\rfloor $. With $\overline{\sigma}_{n,k}=\frac 1n\sigma_{n,r}={\mathrm{Cov%
}}(G_{n,r})$, one has $H^0_{\overline{N},\overline{\sigma}_{n,r}}=
h^0_{N,\sigma_{n,r}}$ and $H^1_{\overline{N},\overline{\sigma}_{n,r}}=
h^1_{N,\sigma_{n,r}}$, $h_{N,\sigma_{n,r} }^{0}$ and $h_{N,\sigma_{n,r}
}^{1} $ being given in (\ref{bd4}) and (\ref{bd5}) respectively. So, we have 
\begin{eqnarray*}
P_{r_{1}+1,r_{2}}^{G,n}P_{k,r_{1}}^{G,n}f(x) &=&{\mathbb{E}}\Big(f\Big(%
x+\sum_{i=k}^{r_{2}-1}G_{n,i}-G_{n,r_{1}}\Big)\Big) \\
&=&{\mathbb{E}}\Big(h_{N,\sigma _{n,r_{1}}}^{0}f\Big(x+%
\sum_{i=k}^{r_{2}-1}G_{n,i}\Big)\Big)+{\mathbb{E}}\Big(h_{N,\sigma
_{n,r_{1}}}^{1}f\Big(x+\sum_{i=k}^{r_{2}-1}G_{n,i}-G_{n,r_{1}}\Big)\Big) \\
&=&P_{r_{1}+1,r_{2}}^{G,n}P_{k,r_{1}}^{G,n}(P_{r_{1},r_{1}+1}^{G,n}h_{N,%
\sigma _{n,r_{1}}}^{0}+h_{N,\sigma _{n,r_{1}}}^{1})f(x) \\
&=&P_{r_{1}+1,r_{2}}^{G,n}P_{k,r_{1}}^{G,n}(U_{n,N,r_{1}}^{0}+U_{n,N,r_{1}}^{1})f(x),
\end{eqnarray*}%
$U_{n,N,r_{1}}^{0}$ and $U_{n,N,r_{1}}^{1}$ being given in (\ref{bd9}). For
every $i=1,2,...,m$, we use this formula in (\ref{bd10}) evaluated in $r=n$
and we get%
\begin{equation}  \label{AmPG2}
\begin{array}{l}
(A^{m}P^{G,n})_{k,n+1}\smallskip \\ 
\displaystyle =\!\!\!\sum_{k\leq r_{1}<\cdots <r_{m}\leq n}
\!\!\!P_{r_{m}+1,n+1}^{G,n}\cdots P_{r_{1}+1,r_{2}}^{G,n}P_{k,r_{1}}^{G,n} %
\Big(\prod_{i=1}^{m}(U_{n,N,r_{i}}^{0}+U_{n,N,r_{i}}^{1})\Big) \Big(%
\prod_{j=1}^{m}(T_{n,N,r_{j}}^{0}+T_{n,N,r_{j}}^{1})\Big).%
\end{array}%
\end{equation}%
Notice that the above operator acts on $C_p^{m(2\lfloor N/2\rfloor+N+5)}({%
\mathbb{R}}^d)$. Our aim now is to isolate the principal term, that is the
sum of the terms where only $U_{n,N,r_{i}}^{0}$ and $T_{n,N,r_{i}}^{0}$
appear. So, we use $Q^{(m)}_{n,N,r_{1},...,r_{m}}$ defined in (\ref{bd8})
and we write%
\begin{eqnarray*}
(A^{m}P^{G,n})_{k,n+1} &=&\sum_{k\leq r_{1}<\cdots <r_{m}\leq
n}P_{r_{m}+1,n+1}^{G,n}\cdots P_{r_{1}+1,r_{2}}^{G,n}P_{k,r_{1}}^{G,n} \Big(%
\prod_{i=1}^{m}U_{n,N,r_{i}}^{0}\Big)\Big(\prod_{j=1}^{m}T_{n,N,r_{j}}^{0}%
\Big) \\
&&+\sum_{k\leq r_{1}<\cdots <r_{m}\leq n}P_{r_{m}+1,n+1}^{G,n}\cdots
P_{r_{1}+1,r_{2}}^{G,n} P_{k,r_{1}}^{G,n}Q^{(m)}_{n,N,r_{1},...,r_{m}}.
\end{eqnarray*}%
The second term is just $R_{n,N,k}^{(m)}$ in (\ref{bd8}). In order to
compute the first one we notice that for every $r^{\prime }<r<r^{\prime
\prime }$ we have 
\begin{equation*}
P_{r+1,r^{\prime \prime }}^{G,n}P_{r^{\prime
},r}^{G,n}P_{r,r+1}^{G,n}=P_{r^{\prime },r^{\prime \prime }}^{G,n}
\end{equation*}%
so that%
\begin{equation*}
P_{r_{m}+1,n+1}^{G,n}\cdots
P_{r_{1}+1,r_{2}}^{G,n}P_{k,r_{1}}^{G,n}(%
\prod_{i=1}^{m}U_{n,N,r_{i}}^{0})=P_{k,n+1}^{G,n}(\prod_{i=1}^{m}h_{N,\sigma
_{n,r_{i}}}^{0}).
\end{equation*}%
Then, for $m=1,...,N$%
\begin{equation*}
(A^{m}P^{G,n})_{k,n+1}=\sum_{k\leq r_{1}<\cdots <r_{m}\leq n}P_{k,n+1}^{G,n}%
\Big(\prod_{i=1}^{m}h_{N,\sigma _{n,r_{i}}}^{0}\Big)\Big(%
\prod_{i=1}^{m}T_{n,N,r_{i}}^{0}\Big)+R_{n,N,k}^{(m)}.
\end{equation*}%
We treat now $A^{N+1}P^{Z,n}$. Using (\ref{bd10}) we get%
\begin{align*}
&(A^{N+1}P^{Z,n})_{k,n+1} \\
& =\sum_{k\leq r_{1}<...<r_{N+1}\leq
n}P_{r_{N+1}+1,n+1}^{Z,n}P_{r_{N}+1,r_{N+1}}^{G,n}\cdots
P_{r_{1}+1,r_{2}}^{G,n}P_{k,r_{1}}^{G,n}%
\prod_{i=1}^{N}(T_{n,N,r_{i}}^{0}+T_{n,N,r_{i}}^{1}) =R_{n,N,k}^{(N+1)},
\end{align*}
which acts on $C_p^{N(N+3)}$. $\square $

\medskip

We give now some useful representations of the remainders.

\begin{lemma}
\label{Nm} Let $m\in \{1,...,N+1\}$ and $r_{1}<\cdots <r_{m}\leq n$ be
fixed. Set $N_{m}:=m(2\lfloor N/2\rfloor+N+5)$ for $m\leq N$ and $%
N_m=(N+1)(N+3)$ otherwise. Then, the operators $Q_{n,N,r_{1},%
\ldots,r_{m}}^{(m)}$ defined in (\ref{bd8}) for $m=1,\ldots,N$ and in (\ref%
{bd8''}) for $m=N+1$, can be written as 
\begin{equation}
Q_{n,N,r_{1},...,r_{m}}^{(m)}f(x)=\frac 1{n^{\frac {N+3m}2}} \sum_{3\leq
\left\vert \alpha \right\vert \leq N_{m}}a_{n,r_{1},...,r_{m}}(\alpha
)\theta _{n,r_{1},...,r_{m}}^{\alpha }\partial_{\alpha }f(x), \quad f\in
C^{N_m}_p({\mathbb{R}}^d),  \label{n1}
\end{equation}%
where $a_{n,r_{1},...,r_{m}}(\alpha )\in {\mathbb{R}}$ are suitable
coefficients with the property%
\begin{equation}  \label{n2}
\left\vert a_{n,r_{1},...,r_{m}}(\alpha )\right\vert \leq (C{\mathrm{C}}%
_2^{\lfloor N/2\rfloor+1}(Y))^m,
\end{equation}%
and $\theta _{n,r_{1},...,r_{m}}^{\alpha }:C_{p}^{\infty }({\mathbb{R}}%
^{d})\rightarrow C_{p}^{\infty }({\mathbb{R}}^{d})$ is an operator which
verifies 
\begin{equation}
|\theta^\alpha_{n,n,r_1,\ldots,r_m}\partial_\alpha f(x)| \leq \big(%
2^{l_{N_m}(f)+1}{\mathrm{C}}_{2(N+3)}^{1/2}(Z)(1+{\mathrm{C}}_{2{l_{N_m}(f)}%
}(Z))^2\big)^m L_{N_m}(f)(1+|x|)^{l_{N_m}(f)}  \label{n3}
\end{equation}%
$C>0$ being a suitable constant. Moreover, $\theta
_{r_{1},\ldots,r_{m}}^{\alpha }$ can be represented as 
\begin{equation}
\theta _{n,r_{1},...,r_{m}}^{\alpha }f(x)=\int_{({\mathbb{R}}^d)^{2m}}
f(x+y_{1}+\cdots+y_{2m})\mu _{n,r_{1},\ldots,r_{m}}^{\alpha
}(dy_{1},\ldots,dy_{2m})  \label{n4}
\end{equation}%
where $\mu _{n,r_{1},...,r_{m}}^{\alpha }$ is a finite signed measure such
that $|\mu _{n,r_{1},...,r_{m}}^{\alpha }|({\mathbb{R}}^{2md})\leq {\mathrm{C%
}}_*^m$, for a suitable constant ${\mathrm{C}}_*$ independent of $n$ and
depending just on $N$ and on the moment bounds ${\mathrm{C}}_p(Y)$ for $p$
large enough.
\end{lemma}

\textbf{Proof}. In a first step we construct the measures $\mu
_{n,r_{1},...,r_{m}}^{\alpha }$ and the operators $\theta
_{n,r_{1},...,r_{m}}^{\alpha }$ and in a second step we prove that the
coefficients $a_{n,r_{1},...,r_{m}}(\alpha )$ verify (\ref{n2}). We start by
representing $T_{n,N,r}^{0}$ defined in (\ref{bd2}). Set 
\begin{equation*}
\nu _{n,r}^{0,\alpha }(dy)=\Delta _{n,r}(\alpha )\delta _{0}(dy),\quad
|\alpha|=l\geq 3.
\end{equation*}
Notice that if $|\alpha|=l\geq 3$ then $\left\vert \Delta _{n,r }(\alpha
)\right\vert \leq 2{\mathrm{C}}_{l}(Y)$. So, we have%
\begin{equation}  \label{n5'}
\begin{array}{l}
\displaystyle T_{n,N,r}^{0}f(x) =\sum_{l=3}^{N+2}\frac{1}{n^{\frac{l}{2}}}%
\sum_{\left\vert \alpha \right\vert =l}\,\frac{1}{l!}\int_{{\mathbb{R}}^d}
\partial_{\alpha }f(x+y)\nu _{n,r}^{0,\alpha }(dy)\quad\mbox{with} \smallskip
\\ 
\displaystyle \int_{{\mathbb{R}}^d} (1+|y|)^\gamma|\nu _{n,r}^{0,\alpha
}|(dy)\leq 2{\mathrm{C}}_l(Y),\quad |\alpha|=l\leq N+2,\quad \gamma\geq 0.%
\end{array}%
\end{equation}%
Concerning $T_{n,N,r}^{1}$ in (\ref{bd3}), for $|\alpha|=N+3$ set $%
\nu^{1,\alpha}_{n,N,r}(dy)=n^{\frac{N+3}2}\int_{0}^{1}(1-\lambda)^{N+2} (%
\frac{y}{\lambda})^\alpha\big[\mu_{\lambda Z_{n,r}}(dy)-\mu_{\lambda
G_{n,r}}(dy)\big]d\lambda$, $\mu_{\lambda Z_{n,r}}$, resp. $\mu_{\lambda
G_{n,r}}$, denoting the law of $\lambda Z_{n,r}$, resp. $\lambda G_{n,r}$.
In other words, 
\begin{align*}
\nu^{1,\alpha}_{n,N,r}(A) &=n^{\frac{N+3}2}\int_0^1(1-\lambda)^{N+2}\big[{%
\mathbb{E}}(Z_{n,r}^\alpha 1_{\lambda Z_{n,r}\in A}) -{\mathbb{E}}%
(G_{n,r}^\alpha 1_{\lambda G_{n,r}\in A})\big]d\lambda,\quad |\alpha|=N+3,
\end{align*}
for every Borel set $A\subset {\mathbb{R}}^d$. Then we have 
\begin{equation}  \label{n5}
\begin{array}{l}
\displaystyle T_{n,N,r}^{1}f(x) =\frac{1}{n^{\frac{1}{2}(N+3)}}%
\sum_{\left\vert \alpha \right\vert =N+3} \frac{1}{(N+2)!}\int_{{\mathbb{R}}%
^d} \partial_{\alpha }f(x+y)\nu _{n,N,r}^{1,\alpha }(dy)\quad\mbox{with}
\smallskip \\ 
\displaystyle \int_{{\mathbb{R}}^d} (1+|y|)^\gamma|\nu _{n,N,r}^{1,\alpha
}|(dy)\leq \frac{2^{\gamma+1}}{N+3}\,{\mathrm{C}}_{2(N+3)}^{1/2}(Y)(1+{%
\mathrm{C}}_{2\gamma}(Y))^{1/2},\quad |\alpha|=N+3,\ \gamma\geq 0.%
\end{array}%
\end{equation}
We represent now the operator $U_{n,N,r}^{0}f(x)={\mathbb{E}}(h_{N,\sigma
_{n,r}}^{0}f(x+G_{n,r}))$ with $h_{N,\sigma _{n,r}}^{0}f$ defined in (\ref%
{bd4}). Notice that 
\begin{equation*}
h_{N,\sigma _{n,r}}^{0} =\mathrm{Id}+\sum_{l=1}^{\lfloor N/2\rfloor} \frac
1{n^l} \sum_{\left\vert \alpha \right\vert =2l}c_{\sigma _{n,r}}(\alpha
)\partial_{\alpha} \quad\mbox{with}\quad c_{\sigma }(\alpha )=\frac{(-1)^l}{%
2^ll!} \prod_{k=1}^l\sigma^{\alpha_{2k-1},\alpha_{2k}},\quad|\alpha|=2l>0.
\end{equation*}
So, by denoting $\rho _{\sigma _{n,r}}^{0}$ the law of $G_{n,r}$, we have 
\begin{equation}  \label{n7}
\begin{array}{l}
\displaystyle U_{n,N,r}^{0}f(x)={\mathbb{E}}(h_{N,\sigma
_{n,r}}^{0}f(x+G_{n,r})) =\sum_{l=0}^{\lfloor N/2\rfloor} \frac 1{n^l}
\sum_{\left\vert \alpha \right\vert =2l}c_{\sigma _{n,r}}(\alpha )\int_{{%
\mathbb{R}}^d} \partial_{\alpha }f(x+y)\rho _{\sigma _{n,r}}^{0}(dy)\quad%
\mbox{with} \smallskip \\ 
\displaystyle |c_{\sigma _{n,r}}(\alpha )|\leq \frac{{\mathrm{C}}_2(Y)^l}{%
2^ll!}\quad\mbox{and}\quad\int_{{\mathbb{R}}^d}(1+|y|)^\gamma |\rho _{\sigma
_{n,r}}^{0}|(dy) \leq 2^\gamma(1+{\mathrm{C}}_{\gamma}(Y)),\ \gamma\geq 0.%
\end{array}%
\end{equation}%
We now obtain a similar representation for $h_{N,\sigma }^{1}f(x)$ defined
in (\ref{bd5}). Set 
\begin{equation*}
\rho^1_\sigma(dy)=\Big(\int_0^1s^{\lfloor N/2\rfloor}\phi_{\sigma^{1/2}\sqrt
s\,W}(y)ds\Big)dy,
\end{equation*}
in which $\phi_{\sigma^{1/2}\sqrt s\,W}$ denotes the density of a centered
Gaussian r.v. with covariance matrix $s\sigma$. Then we write 
\begin{align*}
&h_{N,\sigma }^{1}f(x) = \frac 1{n^{\lfloor N/2\rfloor+1}}
\sum_{|\alpha|=2(\lfloor N/2\rfloor+1)} \!\!\!b_{N,\sigma}(\alpha )\int_{{%
\mathbb{R}}^d} \partial_{\alpha }f(x+y)\rho _{\sigma }^{1}(dy) \mbox{ with }
\\
&b_{N,\sigma}(\alpha )=\frac{(-1)^{\lfloor N/2\rfloor+1}}{2^{\lfloor
N/2\rfloor+1}\lfloor N/2\rfloor!} \prod_{k=1}^{\lfloor
N/2\rfloor+1}\sigma^{\alpha_{2k-1},\alpha_{2k}}.
\end{align*}
%
%
So, we have 
\begin{equation}  \label{n6}
\begin{array}{l}
\displaystyle U_{n,N,r}^{1}f(x)=h_{N,\sigma _{n,r}}^{1}f(x) =\frac
1{n^{\lfloor N/2\rfloor+1}} \sum_{\left\vert \alpha \right\vert =2(\lfloor
N/2\rfloor+1)}b_{N,\sigma _{n,r}}(\alpha )\int \partial _{\alpha }f(x+y)\rho
_{\sigma _{n,r}}^{1}(dy)\quad\mbox{with} \smallskip \\ 
\displaystyle |b_{N,\sigma_{n,r} }(\alpha )| \leq \frac{1}{2^{\lfloor
N/2\rfloor+1}\lfloor N/2\rfloor!} {\mathrm{C}}_2(Y)^{\lfloor N/2\rfloor+1}
\smallskip \\ 
\mbox{and}\quad \int_{{\mathbb{R}}^d}(1+|y|)^\gamma|\rho _{\sigma
_{n,r}}^{1}|(dy) \leq \frac{2^\gamma}{\lfloor N/2\rfloor+1}(1+{\mathrm{C}}%
_\gamma(Y)),\ \gamma\geq 0.%
\end{array}%
\end{equation}%
Using (\ref{n5'}), (\ref{n5}), (\ref{n7}) and (\ref{n6}) we obtain (\ref{n1}%
) with the measure $\mu _{n,r_{1},...,r_{m}}^{\alpha }$ from (\ref{n4})
constructed in the following way:%
\begin{align*}
&\int_{{\mathbb{R}}^{d\times 2m}}f(y_{1},\ldots,y_{m},\bar y_1,\ldots,\bar
y_m)\mu _{n,r_{1},...,r_{m}}^{\alpha }(dy_{1},\ldots,dy_{m},d\bar
y_1,\ldots,d\bar y_m) \\
&=\int_{{\mathbb{R}}^{d\times 2m}}f(y_{1},\ldots,y_{m},\bar y_1,\ldots,\bar
y_m)\eta _{1}(dy_{1})\cdots\eta _{m}(dy_{m})\bar\eta _{1}(d\bar
y_{1})\cdots\bar\eta _{m}(d\bar y_{m})
\end{align*}%
where $\eta _{i}$ is one of the measures $\nu _{n,r_{i}}^{q,\beta }$, $q=0,1$%
, and $\bar \eta _{i}$ is one of the measures $\rho _{\sigma _{n,r_{i}}}^{q}$%
, $q=0,1.$

Let us check that the coefficients $a_{n,r_{1},...,r_{m}}(\alpha )$ which
will appear in (\ref{n1}) verify the bounds in (\ref{n2}). Take first $m\in
\{1,...,N\}.$ Then $Q_{n,r_{1},...,r_{m}}^{(m)}$ is the sum of $%
(\prod_{i=1}^{m}U_{n,N,r_{i}}^{q_{i}^{\prime
}})(\prod_{j=1}^{m}T_{n,N,r_{j}}^{q_{j}})$ where $q_{i},q_{i}^{\prime }\in
\{0,1\}$ and at least one of them is equal to one. And $%
a_{n}^{r_{1},...,r_{m}}(\alpha )$ is the product of coefficients which
appear in the representation of $U_{n,N,r_{i}}^{q_{i}^{\prime }}$ and $%
T_{n,N,r_{j}}^{q_{j}}.$ Recall that the coefficients of $T_{n,N,r_{j}}^{0}$
are all bounded by $C n^{-3/2}$ and the coefficients of $T_{n,N,r_{j}}^{1}$
are bounded by $Cn^{-\frac{1}{2}(N+3)}$. Moreover the coefficients of $%
U_{n,N,r_{i}}^{0}$ are bounded by $C {\mathrm{C}}_{2}^{\lfloor
N/2\rfloor}(Y) $ and the coefficients of $U_{n,N,r_{i}}^{1}$ are bounded by $%
C{\mathrm{C}}_{2}^{\lfloor N/2\rfloor+1}(Y)n^{-(\lfloor N/2\rfloor+1)}$.
Therefore, $(\prod_{i=1}^{m}U_{n,N,r_{i}}^{q_{i}^{\prime
}})(\prod_{j=1}^{m}T_{n,N,r_{j}}^{q_{j}})$ is upper bounded by 
\begin{align*}
&\Big(\frac {C}{n^{\frac{1}{2}(N+3)}}\Big)^{\sum_{i=1}^mq_i} \times \Big(%
\frac{C }{n^{3/2}}\Big)^{\sum_{i=1}^m(1-q_i)} \times \Big( \frac{C{\mathrm{C}%
}_{2}^{\lfloor N/2\rfloor+1}(Y)}{n^{\lfloor N/2\rfloor+1}} \Big)%
^{\sum_{i=1}^mq^{\prime }_i} \times \big(C{\mathrm{C}}_{2}^{\lfloor
N/2\rfloor}(Y)\big)^{\sum_{i=1}^m(1-q^{\prime }_i)} \\
&\leq \Big(\frac {1}{n^{\frac{1}{2}N}}\Big)^{\sum_{i=1}^mq_i}\times \frac{C^m%
}{n^{\frac{3m}2}}\times \Big(\frac 1{n^{\lfloor N/2\rfloor+1}}\Big)%
^{\sum_{i=1}^mq^{\prime }_i}\times (C {\mathrm{C}}_{2}^{\lfloor
N/2\rfloor+1}(Y))^m \\
&\leq \frac {(C{\mathrm{C}}_2^{\lfloor N/2\rfloor+1}(Y))^m}{n^{\frac
N2\sum_{i=1}^mq_i+ (\lfloor N/2\rfloor+1)\sum_{i=1}^mq^{\prime }_i+\frac{3m}%
2}} \leq \frac {(C{\mathrm{C}}_2^{\lfloor N/2\rfloor+1}(Y))^m}{n^{\frac
N2(\sum_{i=1}^mq_i+ \sum_{i=1}^mq^{\prime }_i)+\frac{3m}2}}\leq \frac {(C{%
\mathrm{C}}_2^{\lfloor N/2\rfloor+1}(Y))^m}{n^{\frac {N+3m}2 }}.
\end{align*}
We finally prove (\ref{n3}). We have 
\begin{align*}
&|\theta^\alpha_{n,r_1,\ldots,r_m}\partial_\alpha f(x)| \leq\int_{{\mathbb{R}%
}^{d\times 2m}}|\partial_\alpha f|\Big(x+\sum_{i=1}^my_{i}+\sum_{j=1}^m\bar
y_j\Big)|\eta| _{1}(dy_{1})\cdots|\eta _{m}|(dy_{m})|\bar\eta| _{1}(d\bar
y_{1})\cdots|\bar\eta _{m}|(d\bar y_{m}) \\
&\leq L_{N_m}(f)(1+|x|)^{l_{N_m}(f)}\Big(\prod_{i=1}^m\int_{{\mathbb{R}}%
^d}(1+|y|)^{l_{N_m}(f)}|\eta_i|(dy)\Big) \Big(\prod_{i=1}^m\int_{{\mathbb{R}}%
^d}(1+|y|)^{l_{N_m}(f)}|\bar \eta_i|(dy)\Big) \\
&\leq L_{N_m}(f)(1+|x|)^{l_{N_m}(f)}\Big((2{\mathrm{C}}_{N+2}(Y))
\vee(2^{l_{N_m}(f)+1}{\mathrm{C}}_{2(N+3)}^{1/2}(Y)(1+{\mathrm{C}}_{2{%
l_{N_m}(f)}}(Y)\Big)^m \\
&\qquad\times \Big(2^{l_{N_m}(f)}(1+{\mathrm{C}}_{l_{N_m}(f)}(Y))\Big)^m \\
&\leq \big(2^{l_{N_m}(f)+1}{\mathrm{C}}_{2(N+3)}^{1/2}(Y)(1+{\mathrm{C}}_{2{%
l_{N_m}(f)}}(Y))^2\big)^m L_{N_m}(f)(1+|x|)^{l_{N_m}(f)}
\end{align*}
because ${\mathrm{C}}_{N+2}(Y)\leq {\mathrm{C}}_{2(N+3)}(Y)^{\frac{N+2}{%
2(N+3)}}\leq {\mathrm{C}}_{2(N+3)}(Y)^{\frac12}$. So the proof concerning $%
Q_{n,N,r_{1},...,r_{m}}^{(m)}$, $m=1,...,N$, is completed. The proof for $%
Q_{n,N, r_{1},...,r_{N+1}}^{(N+1)}$ is clearly the same. $\square $

\medskip

We give now the representation of the ``principal term'':

\begin{lemma}
Let the set-up of Proposition \ref{prop1} hold. Then, 
\begin{equation}
\sum_{m=1}^{N}\sum_{1\leq r_{1}<...<r_{m}\leq n}\Big(%
\prod_{i=1}^{m}T_{n,N,r_{i}}^{0}\Big)\Big(\prod_{j=1}^{m}h_{N,\sigma
_{n,r_{j}}}^{0}\Big) =\sum_{k=1}^{N}\frac 1{n^{k/2}}\Gamma _{n,k}+Q_{n,N}^{0}
\label{n8}
\end{equation}%
with $\Gamma _{n,k}$ defined in (\ref{M2}) and 
\begin{equation}  \label{n9}
\begin{array}{l}
\displaystyle Q_{n,N}^{0}=\frac 1{n^{(N+1)/2}}\sum_{N+1\leq \left\vert
\alpha \right\vert \leq N(N+2\lfloor N/2\rfloor)}c_{n,N}(\alpha
)\partial_{\alpha }\smallskip \\ 
\displaystyle \mbox{with}\quad |c_{n,N}(\alpha )| \leq(C{\mathrm{C}}_{N+1}(Y)%
{\mathrm{C}}_2(Y))^{N(N+2\lfloor N/2\rfloor)}%
\end{array}%
\end{equation}
\end{lemma}

\textbf{Proof}. Let $\Lambda _{m}$ and $\Lambda _{m,k}$ be the sets in (\ref%
{M1}). Notice that, for fixed $m$, the $\Lambda _{m,k}$'s are disjoint as $k$
varies. Suppose that $m\in \{1,...,N\}.$ Then $\Lambda _{m,k}=\emptyset $ if 
$k\notin \{m,\ldots,N(N+2\lfloor N/2\rfloor)\}$ so that $\Lambda _{m}=\cup
_{k=m}^{2N(N+2)}\Lambda _{m,k}$ and consequently%
\begin{equation*}
\cup _{m=1}^{N}\Lambda _{m}=\cup _{m=1}^{N}\cup _{k=m}^{N(N+2\lfloor
N/2\rfloor)}\Lambda _{m,k}=\cup _{k=1}^{N(N+2\lfloor N/2\rfloor)}\cup
_{m=1}^{k}\Lambda _{m,k}.
\end{equation*}%
It follows that 
\begin{eqnarray*}
&&\sum_{m=1}^{N}\sum_{1\leq r_{1}<...<r_{m}\leq n}\Big(%
\prod_{i=1}^{m}T_{n,N,r_{i}}^{0}\Big)\Big(\prod_{j=1}^{m}h_{N,\sigma
_{n,r_{j}}}^{0}\Big) \\
&=&\sum_{m=1}^{N}\sum_{l_{1},..,l_{m}=3}^{N+2}\sum_{l_{1}^{\prime
},..,l_{m}^{\prime }=0}^{\lfloor N/2\rfloor}\sum_{1\leq r_{1}<...<r_{m}\leq
n}\Big(\prod_{i=1}^{m}\frac{1}{n^{l_i/2}\,l_{i}!}D_{n,r_{i}}^{(l_{i})}\Big)%
\Big(\prod_{j=1}^{m}\frac{(-1)^{l_{j}^{\prime }}}{n^{l^{\prime
}_j}\,2^{l_{j}^{\prime }}l_{j}^{\prime }!}L_{\sigma
_{n,r_{j}}}^{l_{j}^{\prime }}\Big) \\
&=&\sum_{k=1}^{N(N+2\lfloor
N/2\rfloor)}\sum_{m=1}^{k}\sum_{(l_{1},l_{1}^{\prime
}),...,(l_{m},l_{m}^{\prime }))\in \Lambda _{m,k}}\sum_{1\leq
r_{1}<...<r_{m}\leq n}\frac 1{n^{(k+2m)/2}}\Big(\prod_{i=1}^{m}\frac{1}{%
l_{i}!}D_{n,r_{i}}^{(l_{i})}\Big)\Big(\prod_{j=1}^{m}\frac{%
(-1)^{l_{j}^{\prime }}}{2^{l_{j}^{\prime }}l_{j}^{\prime }!}L_{\sigma
_{n,r_{j}}}^{l_{j}^{\prime }}\Big) \\
&=&\sum_{k=1}^{N}\frac 1{n^{k/2}}\Gamma _{n,k}+ Q_{n,N}^{0}
\end{eqnarray*}%
with%
\begin{align*}
&Q_{n,N}^{0}=\frac 1{n^{{N+1}/2}} \\
&\quad \times \sum_{k=N+1}^{N(N+2\lfloor
N/2\rfloor)}\sum_{m=1}^{k}\sum_{(l_{1},l_{1}^{\prime
}),...,(l_{m},l_{m}^{\prime }))\in \Lambda _{m,k}}\sum_{1\leq
r_{1}<...<r_{m}\leq n}\frac {n^{(N+1)/2}}{n^{(k+2m)/2}}\Big(\prod_{i=1}^{m}%
\frac{1}{l_{i}!}D_{n,r_{i}}^{l_{i}}\Big)\Big(\prod_{j=1}^{m}\frac{%
(-1)^{l_{j}^{\prime }}}{2^{l_{j}^{\prime }}l_{j}^{\prime }!}L_{\sigma
_{n,r_{j}}}^{l_{j}^{\prime }}\Big),
\end{align*}%
which is a differential operator of the form (\ref{n9}). Moreover, 
\begin{align*}
|c_{n,N}(\alpha)| &\leq n^m\times \frac 1{n^m}\times \prod_{i=1}^m\Big(\frac{%
2 {\mathrm{C}}_{N+1}(Y)}{l_i!}\Big) \times \prod_{i=1}^m\Big( \frac{{\mathrm{%
C}}_2^{l^{\prime }_i}(Y)}{2^{l^{\prime }_i}l^{\prime }_i!}\Big) \leq (C{%
\mathrm{C}}_{N+1}(Y){\mathrm{C}}_2(Y))^m \\
& \leq (C{\mathrm{C}}_{N+1}(Y){\mathrm{C}}_2(Y))^{m}
\end{align*}
and the estimate in (\ref{n9}) follows. $\square $

\medskip

We are now ready for the

\medskip

\textbf{Proof of Theorem \ref{Smooth}.} We denote $%
P_{n}^{Z,n}=P_{1,n+1}^{Z,n}$ and $P_{n}^{G,n}=P_{1,n+1}^{G,n}$, so that,
since $S_n(Y)=\overline{S}_n(Z)$, 
\begin{equation*}
{\mathbb{E}}(f(S_{n}(Y))-{\mathbb{E}}(f(W)\Phi _{N}(W))
=P_n^{Z,n}f(0)-P_{n}^{G,n}\Big(\mathrm{Id}+\sum_{k=1}^{N}\frac
1{n^{k/2}}\Gamma _{n,k}\Big)f(0).
\end{equation*}
Putting together (\ref{bd7}) and (\ref{n8}), we can write 
\begin{equation}
P_{n}^{Z,n}f(x)=P_{n}^{G,n}\Big(\mathrm{Id}+\sum_{k=1}^{N}\frac
1{n^{k/2}}\Gamma _{k}\Big)f(x)+I_{1}f(x)+I_{2}f(x)+I_{3}f(x)  \label{n10}
\end{equation}%
with%
\begin{equation}  \label{n11}
\begin{array}{rcl}
I_{1}f(x) & = & \displaystyle \frac
1{n^{(N+1)/2}}P_{n}^{G,n}Q_{n,N}^{0}f(x), \smallskip \\ 
I_{2}f(x) & = & \displaystyle \sum_{1\leq r_{1}<...<r_{N+1}\leq
n}P_{r_{N+1}+1,n}^{Z,n}P_{r_{N}+1,r_{N+1}}^{G,n}\cdots
P_{r_{1}+1,r_{2}}^{G,n}P_{k,r_{1}}^{G,n}Q_{n,N,r_{1},...,r_{N+1}}^{(N+1)}f(x) \smallskip
\\ 
I_{3}f(x) & = & \displaystyle \sum_{m=1}^{N}\sum_{1\leq r_{1}<...<r_{m}\leq
n}P_{r_{m}+1,n}^{G,n}P_{r_{m-1}+1,r_{m}}^{G,n}\cdots
P_{r_{1}+1,r_{2}}^{G,n}P_{k,r_{1}}^{G,n}Q_{n,N,r_{1},...,r_{m}}^{(m)}f(x),%
\end{array}%
\end{equation}
so it is sufficient to study the remaining terms $I_1$, $I_2$ and $I_3$
above. In such study, we will use the following easy consequence of
Burkholder's inequality for discrete martingales: if $M_{n}=\sum_{k=1}^{n}%
\Delta _{k}$ with $\Delta _{k},k=1,...,n$ independent centered random
variables, then 
\begin{equation}
\left\Vert M_{n}\right\Vert _{p}\leq C{\mathbb{E}}\Big(\Big(%
\sum_{k=1}^{n}\left\vert \Delta _{k}\right\vert ^{2}\Big)^{p/2}\Big)^{1/p}=C%
\Big\|\sum_{k=1}^{n}\left\vert \Delta _{k}\right\vert ^{2}\Big\|%
_{p/2}^{1/2}\leq C\Big(\sum_{k=1}^{n}\left\Vert \Delta _{k}\right\Vert
_{p}^{2}\Big)^{1/2}.  \label{b}
\end{equation}

We first estimate $I_1f$. Let us set $N_\ast=N(N+2\lfloor N/2\rfloor)$. So, (%
\ref{n9}) gives 
\begin{align*}
&|I_1f(x)| \leq \frac 1{n^{(N+1)/2}} \!\!\!\!\!\sum_{N+1\leq \left\vert
\alpha \right\vert \leq N_\ast} \!\!\!\!\!\left\vert c_{n}(\alpha
)\right\vert |P_{n}^{G,n}\partial_{\alpha }f(x)| \\
&\leq \frac 1{n^{(N+1)/2}} \!\!\!\!\! \sum_{N+1\leq \left\vert \alpha
\right\vert \leq N_\ast} \!\!\!\!\!\left\vert c_{n}(\alpha )\right\vert
L_{N_\ast}(f)(1+|x|)^{l_{NN_\ast}(f)} {\mathbb{E}}\Big(\Big(1+\Big|%
\sum_{k=1}^nG_{n,k}\Big|\Big)^{l_{NN_\ast}(f)}\Big) \\
&\leq \frac 1{n^{\frac {N+1}2}} {(C{\mathrm{C}}_{N+1}(Y){\mathrm{C}}%
_2(Y))^{NN_\ast}}L_{NN_\ast}(f)(1+|x|)^{l_{NN_\ast}(f)} \times
2^{l_{NN_\ast}(f)}(1+{\mathrm{C}}_{2{l_{NN_\ast}(f)}}(Y)),
\end{align*}
in which we have used the Burkholder inequality (\ref{b}).

The study of $I_2$ and $I_3$ is similar, so we consider $I_3$. Take $m\in
\{1,...,N\}.$ We use Lemma \ref{Nm} (recall $N_m$ given therein) and in
particular (\ref{n1}): 
\begin{eqnarray*}
&&P_{r_{m}+1,n}^{G,n}P_{r_{m-1}+1,r_{m}}^{G,n}\cdots
P_{r_{1}+1,r_{2}}^{G,n}P_{1,r_{1}}^{G,n}Q_{n,N,r_{1},\ldots ,r_{m}}^{(m)}f \\
&=& \frac 1{n^{\frac {N+3m}2}} \sum_{3\leq \left\vert \alpha \right\vert
\leq N_{m}}a_{n,r_{1},...,r_{m}}(\alpha
)P_{r_{m}+1,n}^{G,n}P_{r_{m-1}+1,r_{m}}^{G,n}\cdots
P_{r_{1}+1,r_{2}}^{G,n}P_{1,r_{1}}^{G,n}\theta _{n,r_{1},...,r_{m}}^{\alpha
}\partial_{\alpha }f.
\end{eqnarray*}%
Notice that if $|g(x)|\leq L(1+|x|)^l$ then 
\begin{align*}
&|P_{r_{m}+1,n}^{G,n}P_{r_{m-1}+1,r_{m}}^{G,n}\cdots
P_{r_{1}+1,r_{2}}^{G,n}P_{1,r_{1}}^{G,n}g(x)| \leq {\mathbb{E}}\Big(L\Big(1+%
\Big|x+\sum_{k=1}^n G_{n,k}1_{k\notin\{r_1,\ldots,r_m\}}\Big|\Big)^l\Big) \\
&\leq L(1+|x|)^l{\mathbb{E}}\Big(\Big(1+\Big|\sum_{k=1}^n
G_{n,k}1_{k\notin\{r_1,\ldots,r_m\}}\Big|\Big)^{l}\Big) \leq L(1+|x|)^l\Big(%
1+\Big\|\sum_{k=1}^n G_{n,k}1_{k\notin\{r_1,\ldots,r_m\}}\Big\|_l\Big)^{l}.
\end{align*}
Since the $G_{n,k}1_{k\notin\{r_1,\ldots,r_m\}}$'s are centered and
independent, with $\|G_{n,k}1_{k\notin\{r_1,\ldots,r_m\}}\|_l\leq {\mathrm{C}%
}_l(Y)/n^{l/2}$, we can use the Burkholder inequality (\ref{b}), 
we get 
\begin{align}  \label{gl}
|P_{r_{m}+1,n}^{G,n}P_{r_{m-1}+1,r_{m}}^{G,n}\cdots
P_{r_{1}+1,r_{2}}^{G,n}P_{1,r_{1}}^{G,n}g(x)| \leq L(1+|x|)^l(1+{\mathrm{C}}%
_l^{1/l}(Y))^{l} \leq 2^l(1+{\mathrm{C}}_l(Y))L(1+|x|)^l.
\end{align}
We use now this inequality with $g=\theta^\alpha_{n,r_1,\ldots,r_m}\partial_%
\alpha f$: by applying (\ref{n3}) we get 
\begin{equation*}
|P_{r_{m}+1,n}^{G,n}P_{r_{m-1}+1,r_{m}}^{G,n}\cdots
P_{r_{1}+1,r_{2}}^{G,n}P_{1,r_{1}}^{G,n}Q_{N,r_{1},...,r_{m}}^{(m)}f(x)|
\leq \mathcal{K}_{N,m}(f)(1+|x|)^{l_{N_m}(f)}
\end{equation*}
with 
\begin{equation*}
\mathcal{K}_{N,m}(f) =2^{l_{N_m}(f)}(1+{\mathrm{C}}_{l_{N_m}(f)}(Y))\big(%
2^{l_{N_m}(f)+1} {\mathrm{C}}_{2(N+3)}^{1/2}(Y)(1+{\mathrm{C}}_{2{l_{N_m}(f)}%
}(Y))^2\big)^m L_{N_m}(f).
\end{equation*}
Moreover, using (\ref{n2}) 
\begin{eqnarray*}
&&| P_{r_{m}+1,n}^{G,n}P_{r_{m-1}+1,r_{m}}^{G,n}\cdots
P_{r_{1}+1,r_{2}}^{G,n}P_{1,r_{1}}^{G,n}Q_{n,N,r_{1},...,r_{m}}^{(m)}f(x)| \\
&\leq &\mathcal{K}_{N,m}(f)(1+|x|)^{l_{N_m}(f)}\frac{1}{n^{\frac{1}{2}(N+3m)}%
}\sum_{0\leq \left\vert \alpha \right\vert \leq N+1}\left\vert
a_{n,r_{1},...,r_{m}}(\alpha )\right\vert \\
&\leq &\mathcal{H}_N \mathcal{K}_{N,m}(f)(1+|x|)^{l_{N_m}(f)}\frac{1}{n^{%
\frac{1}{2}(N+3m)}}(C{\mathrm{C}}_2^{\lfloor N/2\rfloor+1}(Y))^m,
\end{eqnarray*}%
$\mathcal{H}_N$ denoting a constant depending on $N$ only. Since the set $%
\{1\leq r_{1}<...<r_{m}\leq n\}$ has less than $n^{m}$ elements, we get 
\begin{align*}
&|I_3f(x)| \leq N\times n^{m}\times \frac{1}{n^{\frac{1}{2}(N+3m)}}\times 
\mathcal{H}_N \mathcal{K}_{N,m}(f)(1+|x|)^{l_{N_m}(f)}(C{\mathrm{C}}%
_2^{\lfloor N/2\rfloor+1}(Y))^m \\
&\leq N\mathcal{H}_N \mathcal{K}_{N,m}(f)(1+|x|)^{l_{N_m}(f)}(C{\mathrm{C}}%
_2^{\lfloor N/2\rfloor+1}(Y))^m \times \frac{1}{n^{\frac{1}{2}(N+1)}}
\end{align*}
The estimate for $I_{2}(f)$ is analogous. So, we get 
\begin{equation*}
\sum_{i=1}^3|I_if(x)| \leq \mathcal{H}_N{\mathrm{C}}_{2(N+3)}^{2N(N+2\lfloor
N/2\rfloor)}(Y)(1+{\mathrm{C}}_{2l_{\widehat{N}}(f)}(Y))^{2N+3} 2^{(N+2)(l_{%
\widehat{N}}(f)+1)}L_{\widehat{N}}(f)(1+|x|)^{l_{\widehat{N}}(f)}\times
\frac 1{n^{\frac {N+1}2}}
\end{equation*}
with $\widehat{N}=N(2\lfloor N/2\rfloor+N+5)$, and statement (\ref{M6})
follows. Concerning (\ref{M6'}), it suffices to notice that for $%
f(x)=x^\beta $ with $|\beta|=k$ then $L_{\widehat{N}}(f)= 1$ and $l_{%
\widehat{N}}(f)=k$. $\square $

\section{The case of general test functions}

\label{general}

\subsection{Differential calculus based on the Nummelin's splitting}

\label{sect-mall}

In this section we use the variational calculus settled in \cite%
{bib:[B.C],bib:[BC0],bib:[BCl],bib:[BCl2]} in order to treat general test
functions. Let us give the definitions and the notation.

\smallskip

We fix $r,\varepsilon >0$ and we consider a sequence of independent random
variables $Y_{k}\in {\mathfrak{D}}(r,\varepsilon ),k\in {\mathbb{N}}.$ Then,
using the Nummelin's splitting (\ref{s4}) we write 
\begin{equation}
Y_{k}=\chi _{k}V_{k}+(1-\chi _{k})U_{k} ,  \label{s7}
\end{equation}%
the law of $\chi _{k}$, $V_{k}$ and $U_{k}$ being given in (\ref{s5}). We
assume that $\chi _{k}$, $V_{k}$, $U_{k}$, $k\in{\mathbb{N}}$, are
independent. We define $\mathcal{G}=\sigma (\chi _{k},U_{k},k\in {\mathbb{N}}%
).$ A random variable $F=f(\omega ,V_{1},...,V_{n})$ is called a simple
functional if $f$ is $\mathcal{G\times B}({\mathbb{R}}^{d\times n})$
measurable and for each $\omega ,f(\omega ,\cdot )\in C_{b}^{\infty }({%
\mathbb{R}}^{d\times n}).$ We denote $\mathcal{S}$ the space of the simple
functionals. Moreover we define the differential operator $D:\mathcal{S}%
\rightarrow l_{2}:=l_{2}({\mathbb{R}}^{d})$ by $D_{(k,i)}F=\chi _{k}\partial
_{v_{k}^{i}}f(\omega ,V_{1},...,V_{n}).$ Then the Malliavin covariance
matrix of $F\in (F^{1},...,F^{m})\in \mathcal{S}^{m}$ is defined as%
\begin{equation}
\sigma _{F}^{i,j}=\left\langle DF^{i},DF^{j}\right\rangle
_{l_{2}}=\sum_{k=1}^{\infty }\sum_{p=1}^{d}D_{(k,p)}F^{i}\times
D_{(k,p)}F^{j},\quad i,j=1,...,m.  \label{s8}
\end{equation}%
If $\sigma _{F}$ is invertible we denote $\gamma _{F}=\sigma _{F}^{-1}.$

Moreover, we define the iterated derivatives $D^{m}:\mathcal{S}\rightarrow
l_{2}^{\otimes m}$ by $D_{(k_{1},i_{1}),%
\ldots,(k_{m},i_{m})}^{(m)}=D_{(k_{1},i_{1})}\cdots D_{(k_{m},i_{m})} $ and
on $\mathcal{S}$ we consider the norms%
\begin{equation*}
\left\vert F\right\vert _{q}^{2}=\left\vert F\right\vert
^{2}+\sum_{m=1}^{q}\left\vert D^{m}F\right\vert _{l_{2}^{\otimes
m}}^{2}=\left\vert F\right\vert
^{2}+\sum_{m=1}^{q}\sum_{k_{1},...,k_{m}=1}^{\infty
}\sum_{i_{1},...,i_{m}=1}^{d}\left\vert
D_{(k_{1},i_{1})}....D_{(k_{m},i_{m})}F\right\vert ^{2}
\end{equation*}%
and 
\begin{equation}
\left\Vert F\right\Vert _{q,p}=({\mathbb{E}}(\left\vert F\right\vert
_{q}^{p}))^{1/p}.  \label{s9}
\end{equation}

We introduce now the Ornstein-Uhlenbeck operator $L$. We denote $%
\theta_{k,i}=\partial_i \ln
p_{V_k}(V_k)=2(V_k-y_Y)^i1_{r<|V_k-y_Y|^2<2r}a^{\prime }_r(|V_k-y_{Y_k}|^2)$%
, $p_{V_k}$ being the density of $V_k$ (see (\ref{s5})) and $a_r$ is given
in (\ref{s2}). So, we define 
\begin{equation}  \label{s10}
LF=-\sum_{k=1}^{\infty }\sum_{i=1}^{d}\big(D_{(k,i)}D_{(k,i)}F+D_{(k,i)}F%
\times \theta_{k,i}\big).
\end{equation}%
Using elementary integration by parts on ${\mathbb{R}}^{d}$ one easily
proves the following duality formula: for $F,G\in \mathcal{S}$%
\begin{equation}
{\mathbb{E}}(\left\langle DF,DG\right\rangle _{l_{2}})={\mathbb{E}}(FLG)={%
\mathbb{E}}(GLF).  \label{s10'}
\end{equation}

Finally, for $q\geq 2,$ we define%
\begin{equation}
\left\Vert \left\vert F\right\vert \right\Vert _{q,p}=\left\Vert
F\right\Vert _{q,p}+\left\Vert LF\right\Vert _{q-2,p}.  \label{s11}
\end{equation}

We recall now the basic computational rules and the integration by parts
formulas. For $\phi \in C^{1}({\mathbb{R}}^{d})$ and $F=(F^{1},...,F^{d})\in 
\mathcal{S}^d$ we have%
\begin{equation}
D\phi (F)=\sum_{j=1}^{d}\partial _{j}\phi (F)DF^{j},  \label{s11'}
\end{equation}%
and for $F,G\in \mathcal{S}$ 
\begin{equation}
L(FG)=FLG+GLF-2\left\langle DF,DG\right\rangle .  \label{s11''}
\end{equation}%
The formula (\ref{s11'}) is just the chain rule in the standard differential
calculus and (\ref{s11''}) is obtained using duality. Let $H\in \mathcal{S}.$
We use the duality relation and (\ref{s10'}) we obtain%
\begin{equation}
{\mathbb{E}}(HFLG) ={\mathbb{E}}(\left\langle D(HF),DG\right\rangle _{l_{2}})
\notag \\
={\mathbb{E}}(H\left\langle DF,DG\right\rangle _{l_{2}})+{\mathbb{E}}%
(F\left\langle DH,DG\right\rangle _{l_{2}}).  \label{s11bis}
\end{equation}
A similar formula holds with $GLF$ instead of $FLG.$ We sum them and we
obtain%
\begin{eqnarray*}
{\mathbb{E}}(H(FLG+GLF)) &=&2{\mathbb{E}}(H\left\langle DF,DG\right\rangle
_{l_{2}})+{\mathbb{E}}(\left\langle DH,D(FG)\right\rangle _{l_{2}}) \\
&=&2{\mathbb{E}}(H\left\langle DF,DG\right\rangle _{l_{2}})+{\mathbb{E}}%
(HL(FG)).
\end{eqnarray*}

We give now the integration by parts formula (this is a localized version of
the standard integration by parts formula from Malliavin calculus).

\begin{theorem}
\label{TH1} Let $\eta >0$ be fixed and let $\Psi _{\eta }\in C^{\infty }({%
\mathbb{R}})$ be such that $1_{[\eta /2,\infty )}\leq \Psi _{\eta }\leq
1_{[\eta ,\infty )}$ and for every $k\in {\mathbb{N}}$ one has $\|\Psi
_{\eta }^{(k)}\|_{\infty }\leq C\eta ^{-k}.$ Let $F\in \mathcal{S}^{d}$ and $%
G\in \mathcal{S}$. For every $\phi \in C_{p}^{\infty }({\mathbb{R}}^{d})$, $%
\eta >0$ and $i=1,...,d$%
\begin{equation}
{\mathbb{E}}(\partial _{i}\phi (F)G\Psi _{\eta }(\det \sigma _{F}))={\mathbb{%
E}}(\phi (F)H_{i}(F,G\Psi _{\eta }(\det \sigma _{F})))  \label{FD13}
\end{equation}%
with 
\begin{equation}
H_{i}(F,G\Psi _{\eta }(\det \sigma _{F}))=\sum_{j=1}^{d}\big(G\Psi _{\eta
}(\det \sigma _{F})\big)\gamma _{F}^{i,j}LF^{j}+\big\langle D(G\Psi _{\eta
}(\det \sigma _{F})\gamma _{F}^{i,j}),DF^{j}\big\rangle _{l_{2}}.
\label{FD14}
\end{equation}%
Let $m\in {\mathbb{N}},m\geq 2$ and $\alpha =(\alpha _{1},...,\alpha
_{m})\in \{1,...,d\}^{m}.$ Then%
\begin{equation}
{\mathbb{E}}(\partial _{\alpha }\phi (F)G\Psi _{\eta }(\det \sigma _{F}))={%
\mathbb{E}}(\phi (F)H_{\alpha }(F,G\Psi _{\eta }(\det \sigma _{F})))
\label{FD15}
\end{equation}%
with $H_{\alpha }(F,G\Psi _{\eta }(\det \sigma _{F}))$ defined by recurrence 
\begin{equation*}
H_{(\alpha _{1},...,\alpha _{m})}(F,G\Psi _{\eta }(\det \sigma
_{F})):=H_{\alpha _{m}}(F,H_{(\alpha _{1},...,\alpha _{m-1})}(F,G\Psi _{\eta
}(\det \sigma _{F}))).
\end{equation*}
\end{theorem}

The proof is standard, for details see e.g. \cite{[BCDi],bib:[BCl]}.

We give now useful estimates for the weights which appear in (\ref{FD15}):

\begin{lemma}
\label{L1} Let $q ,m \in {\mathbb{N}}$ and $F\in \mathcal{S}^{d}$ and $G\in 
\mathcal{S}.$ There exists a universal constant ${\mathrm{C}}\geq 1$
(depending on $d,q ,m $ only) such that for every multiindex $\alpha $ with $%
\left\vert \alpha \right\vert =q $ and every $\eta >0$ one has%
\begin{equation}
\big\vert H_{\alpha}(F,\Psi _{\eta }(\det \sigma _{F})G)\big\vert _{m }\leq 
\frac{{\mathrm{C}}}{\eta ^{2q +m }}\times {\mathcal{K}}_{q ,m }(F)\times
\left\vert G\right\vert _{m +q },  \label{FD16}
\end{equation}%
with 
\begin{equation}
{\mathcal{K}}_{q ,m }(F)=(\left\vert F\right\vert _{1,m +q +1}+\left\vert
LF\right\vert _{m +q })^{q }(1+\left\vert F\right\vert _{1,m +q +1})^{2d(2q
+m )}.  \label{FD15'}
\end{equation}

In particular, taking $m =0$ and $G=1$ we have%
\begin{equation}
\left\Vert H_{\alpha }(F,\Psi _{\eta }(\det \sigma _{F}))\right\Vert
_{p}\leq \frac{C}{\eta ^{2q }}\times \left\Vert {\mathcal{K}}_{q
,0}(F)\right\Vert _{p}  \label{FD4}
\end{equation}
\end{lemma}

The proof is straightforward but technical so we leave it for Appendix \ref%
{app:norms}.

\medskip

We go now on and we give the regularization lemma. We recall that a super
kernel $\phi :{\mathbb{R}}^{d}\rightarrow {\mathbb{R}}$ is a function which
belongs to the Schwartz space $\mathbf{S}$\ (infinitely differentiable
functions which decrease in a polynomial way to infinity), $\int \phi
(x)dx=1,$ and such that for every multiindexes $\alpha$ and $\beta$, one has 
\begin{align}
\int y^{\alpha }\phi (y)dy &=0,\quad |\alpha|\geq 1,  \label{kk1} \\
\int \left\vert y\right\vert ^{m}\left\vert \partial_\beta \phi
(y)\right\vert dy &<\infty.  \label{kk2}
\end{align}%
As usual, for $|\alpha|=m$ then $y^{\alpha }=\prod_{i=1}^{m}y_{\alpha _{i}}$%
. Since super kernels play a crucial role in our approach we give here the
construction of such an object (we follow \cite{[KK]} Section 3, Remark 1).
We do it in dimension $d=1$ and then we take tensor products. So, if $d=1$
we take $\psi \in \mathbf{S}$ which is symmetric and equal to one in a
neighborhood of zero and we define $\phi =\mathcal{F}^{-1}\psi ,$ the
inverse of the Fourier transform of $\psi .$ Since $\mathcal{F}^{-1}$ sends $%
\mathbf{S}$\ into $\mathbf{S}$ the property (\ref{kk2}) is verified. And we
also have $0=\psi ^{(m)}(0)=i^{-m}\int x^{m}\phi (x)dx$ so (\ref{kk1}) holds
as well. We finally normalize in order to obtain $\int \phi =1.$

We fix a super kernel $\phi $. For $\delta \in (0,1)$ and for a function $f$
we define 
\begin{equation*}
\phi _{\delta }(y)=\frac{1}{\delta ^{d}}\phi \Big(\frac{y}{\delta }\Big)\quad%
\mbox{and}\quad f_{\delta }=f\ast \phi _{\delta },
\end{equation*}%
the symbol $\ast$ denoting convolution. For $f\in C_{p}^{k}({\mathbb{R}}%
^{d}) $, we recall the constants $L_{k}(f)$ and $l_{k}(f)$ in (\ref{Lklk}).

\begin{lemma}
\label{R-1} Let $F\in \mathcal{S}^{d}$ and $q,m\in {\mathbb{N}}.$ There
exists some constant $C\geq 1,$ depending on $d,m$ and $q $ only, such that
for every $f\in C_{{\mathrm{{\scriptsize {pol}}}}}^{q+m}({\mathbb{R}}^{d}),$
every multiindex $\gamma $ with $\left\vert \gamma \right\vert =m$ and every 
$\eta ,\delta >0 $%
\begin{equation}  \label{1a}
\begin{array}{l}
\displaystyle \left\vert {\mathbb{E}}(\Psi _{\eta }(\det \sigma
_{F})\partial_{\gamma }f(x+F))-{\mathbb{E}}(\Psi _{\eta }(\det \sigma
_{F})\partial_{\gamma }f_{\delta }(x+F))\right\vert \smallskip \\ 
\displaystyle \quad \leq C\,2^{l_0(f)}c_{l_0(f)+q} L_{0}(f)\left\Vert
F\right\Vert _{2l_{0}(f)}^{l_{0}(f)}\mathcal{C}_{q+m}(F)\frac{\delta ^{q}}{%
\eta ^{2(q+m)}}(1+|x|)^{l_0(f)}%
\end{array}%
\end{equation}
with 
\begin{equation}
c_{p}=\int \left\vert \phi (z)\right\vert (1+\left\vert z\right\vert
)^{p}dz\quad\mbox{and}\quad\mathcal{C}_{p}(F)= \left\Vert {\mathcal{K}}%
_{p,0}(F)\right\Vert _{2},  \label{1aa}
\end{equation}%
${\mathcal{K}}_{p,0}(F)$ being defined in (\ref{FD15'}). Moreover, for every 
$p>1$%
\begin{equation}  \label{1b}
\begin{array}{l}
\displaystyle \left\vert {\mathbb{E}}(\partial_{\gamma }f(x+F))-{\mathbb{E}}%
(\partial_{\gamma }f_{\delta }(x+F))\right\vert \leq C(1+\left\Vert
F\right\Vert _{pl_{0}(f)})^{l_{0}(f)}(1+|x|)^{l_m(f)} \smallskip \\ 
\displaystyle \quad \times \Big(L_{m}(f)c_{l_m(f)}2^{l_{m}(f)}{\mathbb{P}}%
^{(p-1)/p}(\det \sigma _{F}\leq \eta )+2^{l_0(f)}c_{l_0(f)+q}\,L_{0}(f)\frac{%
\delta ^{q}}{\eta ^{2(q+m)}}\mathcal{C}_{q+m}(F)\Big).%
\end{array}%
\end{equation}
\end{lemma}

\textbf{Proof A}. Using Taylor expansion of order $q$ 
\begin{align*}
\partial_{\gamma }f(x)-\partial_{\gamma }f_{\delta }(x)& =\int
(\partial_{\gamma }f(x)-\partial_{\gamma }f(y))\phi _{\delta }(x-y)dy \\
& =\int I_{\gamma,q}(x,y)\phi _{\delta }(x-y)dy+\int R_{\gamma,q}(x,y)\phi
_{\delta }(x-y)dy
\end{align*}%
with 
\begin{align*}
I_{\gamma,q}(x,y)& =\sum_{i=1}^{q-1}\frac{1}{i!}\sum_{\left\vert \alpha
\right\vert =i}\partial_{\alpha }\partial_{\gamma }f(x)(x-y)^{\alpha }, \\
R_{\gamma,q}(x,y)& =\frac{1}{q!}\sum_{\left\vert \alpha \right\vert
=q}\int_{0}^{1}\partial_{\alpha }\partial_{\gamma }f(x+\lambda
(y-x))(x-y)^{\alpha }(1-\lambda) ^{q}d\lambda .
\end{align*}%
Using (\ref{kk1}) we obtain $\int I(x,y)\phi _{\delta }(x-y)dy=0$ and by a
change of variable we get 
\begin{equation*}
\int R_{\gamma,q}(x,y)\phi _{\delta }(x-y)dy=\frac{1}{q!}\sum_{\left\vert
\alpha \right\vert =q}\int_{0}^{1}\int dz\phi _{\delta }(z)\partial_{\alpha
}\partial_{\gamma }f(x+\lambda z)z^{\alpha }(1-\lambda) ^{q}d\lambda .
\end{equation*}%
So, we have 
\begin{align*}
& {\mathbb{E}}(\Psi _{\eta }(\det \sigma _{F})\partial_{\gamma }f(x+F))-{%
\mathbb{E}}(\Psi _{\eta }(\det \sigma _{F})\partial_{\gamma }f_{\delta
}(x+F)) \\
& ={\mathbb{E}}\Big(\int \Psi _{\eta }(\det \sigma
_{F})R_{\gamma,q}(x+F,y)\phi _{\delta }(x+F-y)dy\Big) \\
& =\frac{1}{q!}\sum_{\left\vert \alpha \right\vert =q}\int_{0}^{1}\int
dz\phi _{\delta }(z){\mathbb{E}}\big(\Psi _{\eta }(\det \sigma
_{F})\partial_{\alpha }\partial_{\gamma }f(x+F+\lambda z)\big)z^{\alpha
}(1-\lambda) ^{q}d\lambda .
\end{align*}%
Using integration by parts formula (\ref{FD15}) (with $G=1)$\ 
\begin{align*}
&\left\vert {\mathbb{E}}(\Psi _{\eta }(\det \sigma _{F})\partial_{\alpha
}\partial_{\gamma }f(x+F+\lambda z))\right\vert \\
& =\left\vert {\mathbb{E}}(f(F+\lambda z)H_{(\gamma,\alpha)}(F,\Psi _{\eta
}(\det \sigma _{F}))\right\vert \\
& \leq L_{0}(f){\mathbb{E}}((1+|x|+\left\vert z\right\vert +\left\vert
F\right\vert )^{l_{0}(f)}\left\vert H_{(\gamma,\alpha)}(F,\Psi _{\eta }(\det
\sigma _{F}))\right\vert ) \\
& \leq C(1+\left\vert x\right\vert )^{l_{0}(f)}(1+\left\vert z\right\vert
)^{l_{0}(f)}L_{0}(f)\left\Vert F\right\Vert _{2l_{0}(f)}^{l_{0}(f)}\|
H_{(\gamma,\alpha)}(F,\Psi _{\eta }(\det \sigma _{F}))\|_2.
\end{align*}%
The upper bound from (\ref{FD4}) (with $p=2$) gives 
\begin{equation*}
\left\Vert H_{\alpha }(F,\Psi _{\eta }(\det \sigma _{F}))\right\Vert
_{2}\leq \frac{C}{\eta ^{2(q+m) }}\times \left\Vert {\mathcal{K}}_{q+m
,0}(F)\right\Vert _{2}
\end{equation*}
And since $\int dz\left\vert \phi _{\delta }(z)z^{\alpha }\right\vert
(1+\left\vert z\right\vert )^{l_{0}(f)} \leq \delta ^{q} \int \left\vert
\phi (z)\right\vert (1+\left\vert z\right\vert )^{|\alpha|+l_{0}(f)}dz$ we
conclude that 
\begin{align*}
&\left\vert {\mathbb{E}}(\Psi _{\eta }(\det \sigma _{F}))\partial_{\gamma
}f(x+F))-{\mathbb{E}}(\Psi _{\eta }(\det \sigma _{F}))\partial_{\gamma
}f_{\delta }(x+F))\right\vert \\
&\leq C(1+\left\vert x\right\vert )^{l_{0}(f)}c_{l_0(f)+q}L_{0}(f)\left\Vert
F\right\Vert _{2l_{0}(f)}^{l_{0}(f)}\left\Vert {\mathcal{K}}_{q+m
,0}(F)\right\Vert _{2}\frac{C\delta ^{q}}{\eta ^{2(q+m)}}
\end{align*}
and (\ref{1a}) holds. Concerning (\ref{1b}), we set $L_{\gamma,%
\delta}=L_{0}(\partial_{\gamma }f_{\delta })\vee L_{0}(\partial_{\gamma }f)$
and $l_{\gamma,\delta}= l_{0}(\partial_{\gamma }f_{\delta })\vee
l_{0}(\partial_{\gamma }f)$. So, for every $p>1$, we have 
\begin{eqnarray*}
&&\left\vert {\mathbb{E}}\big((1-\Psi _{\eta }(\det \sigma
_{F}))\partial_{\gamma }f(x+F)\big)-{\mathbb{E}}\big((1-\Psi _{\eta }(\det
\sigma _{F}))\partial_{\gamma }f_{\delta }(x+F)\big)\right\vert \\
&\leq &2L_{\gamma,\delta}{\mathbb{E}}\big((1-\Psi _{\eta }(\det \sigma
_{F}))(1+|x|+\left\vert F\right\vert )^{l_{\gamma,\delta}}\big) \\
&\leq
&2L_{\gamma,\delta}2^{l_{\gamma,\delta}}(1+|x|)^{l_{\gamma,\delta}}(1+\left%
\Vert F\right\Vert _{pl_{0}(f_{\delta })\vee l_{0}(f)})^{l_{0}(f_{\delta
})\vee l_{0}(f)}{\mathbb{P}}^{(p-1)/p}(\det \sigma _{F}\leq \eta ).
\end{eqnarray*}%
So the proof of (\ref{1b}) will be completed as soon as we check that $%
l_0(\partial_\gamma f_\delta)\leq l_m(f)$ and $L_0(\partial_\gamma f_{\delta
}) \leq L_{m}(f)\int(1+\left\vert y\right\vert )^{l_{m}(f)}\left\vert
\phi(y)\right\vert dy=L_m(f)c_{l_m(f)}$: 
\begin{align*}
\left\vert \partial_{\gamma }f_{\delta }(x)\right\vert &=\left\vert \int
\partial_{\gamma }f(x-y)\phi _{\delta }(y)dy\right\vert \leq L_{m}(f)\int
(1+\left\vert x-y\right\vert )^{l_{m}(f)}\left\vert \phi _{\delta
}(y)\right\vert dy \\
& \leq L_{m}(f)(1+\left\vert x\right\vert )^{l_{m}(f)}\int (1+\left\vert
y\right\vert )^{l_{m}(f)}\left\vert \phi (y)\right\vert dy.
\end{align*}%
$\square $

\subsection{CLT and Edgeworth's development}

\label{Mall-CLT}

In this section we take $F=S_{n}(Y)=\frac 1{\sqrt n}\sum_{k=1}^{n}C_{n,k}Y_k$
defined in (\ref{BD1}), and we recall that $\sigma_{n,k}=C_{n,k}C^\ast_{n,k}=%
{\mathrm{Cov}}(C_{n,k}Y_k)$. From now on, we assume that $Y_{k}\in {%
\mathfrak{D}}(r,\varepsilon )$ so we have the decomposition (\ref{s7}).
Consequently%
\begin{equation*}
F=S_{n}(Y)=\frac 1{\sqrt n}\sum_{k=1}^{n}C_{n,k}Y_{k}=\frac 1{\sqrt
n}\sum_{k=1}^{n}C_{n,k}(\chi _{k}V_{k}+(1-\chi _{k})U_{k}).
\end{equation*}%
We will use Lemma \ref{R-1}, so we estimate the quantities which appear in
the right hand side of (\ref{1a}).

\begin{lemma}
Let $Y_k\in{\mathfrak{D}}(2\varepsilon, r)$ and let the moment bounds
condition (\ref{bd00}) hold. For every $k\in {\mathbb{N}}$ and $p\geq 2$
there exists a constant $C$ depending on $k,p$ only, such that 
\begin{equation}
\sup_{n}\left\Vert S_{n}(Y)\right\Vert _{k,p}\leq 2\times {\mathrm{C}}%
_{p}(Y) \quad\mbox{and}\quad \sup_{n}\left\Vert \left\vert
S_{n}(Y)\right\vert \right\Vert _{k,p}\leq C\times \frac{{\mathrm{C}}_{p}(Y)%
}{r^{k}}.  \label{s15}
\end{equation}
\end{lemma}

\textbf{Proof}. Using the Burkholder inequality (\ref{b}) and (\ref{bd00})
we obtain $\left\Vert S_{n}(Y)\right\Vert _{p}\leq C\times {\mathrm{C}}%
_{p}(Y).$ We look now to the Sobolev norms. It is easy to see that, $%
S_n(Y)^i $ denoting the $i$th component of $S_n(Y)$,%
\begin{equation*}
D_{(k,j)}S_n(Y)^{i}=\frac 1{\sqrt n}\chi_k C_{n,k}^{i,j}\quad \mbox{and}%
\quad D^{(l)}S_n(Y)=0\text{ for }l\geq 2.
\end{equation*}%
Since $\frac 1n\sum_{k=1}^{n}\left\vert \sigma _{n,k}\right\vert \leq {%
\mathrm{C}}_{2}(Y)$ it follows that%
\begin{equation*}
\left\Vert S_{n}(Y)\right\Vert _{k,p}\leq 2 {\mathrm{C}}_{p}(Y)\quad \forall
k\in {\mathbb{N}},p\geq 2.
\end{equation*}%
Moreover%
\begin{align*}
LS_n(Y)=-\frac 1{\sqrt n}\sum_{k=1}^n C_{n,k}LY_k=- \frac 1{\sqrt
n}\sum_{k=1}^{n}\chi_kC_{n,k} A_r(V_k), \\
\mbox{with}\quad A_r(V_k)=1_{r<|V_k-y_Y|^2<2r}\times 2 a^{\prime
}_r(|V_k-y_{Y_k}|^2)(V_k-y_{Y_k}).
\end{align*}%
We prove that 
\begin{equation}
\left\Vert LS_{n}(Y)\right\Vert _{k,p}\leq \frac{C}{r^{k}}\times {\mathrm{C}}%
_{p}(Y),  \label{s14}
\end{equation}
$C$ depending on $k,p$ but being independent of $n$.

Let $k=0.$ The duality relation gives ${\mathbb{E}}(LY_{k})={\mathbb{E}}%
(\left\langle D1,DY_{k}\right\rangle _{l_{2}})=0$. Since the $LY_k$'s are
independent, we can apply (\ref{b}) first and then (\ref{bd00}), so that 
\begin{equation*}
\left\Vert LS_{n}(Y)\right\Vert _{p} \leq C\Big(\frac
1n\sum_{k=1}^{n}\left\Vert C_{n,k}A_r(V_k)\right\Vert _{p}^{2}\Big)^{1/2} 
\newline
\leq C\Big({\mathrm{C}}_{2}(Y)\frac 1n\sum_{k=1}^n\left\Vert
A_r(V_k)\right\Vert _{p}^{2}\Big)^{1/2}
\end{equation*}
By (\ref{s3}) ${\mathbb{E}}(|A_r(V_k)|^{p})\leq Cr^{-p}$ so $\left\Vert
LS_{n}(Y)\right\Vert _{p}\leq Cr^{-1}\times {\mathrm{C}}_{2}(Y)$ and (\ref%
{s14}) follows for $k=0$.

We take now $k=1.$ We have%
\begin{equation*}
D_{(q,j)}LS_n(Y)^{i}=\frac 1{\sqrt n}D_{(q,j)}\big(\chi_kC_{n,q}A_r(V_q)%
\big) =\frac 1{\sqrt n}\chi_kC_{n,q}D_{(q,j)}A_r(V_q)
\end{equation*}%
so that, using again (\ref{bd00}),%
\begin{equation*}
\left\vert DLS_n(Y)\right\vert _{l_{2}}^{2} =\frac
1n\sum_{q=1}^{n}\sum_{j=1}^{d}\left\vert \chi_kC_{n,q}
D_{(q,j)}A_r(V_q)\right\vert ^{2} \leq C\times \frac{{\mathrm{C}}_2(Y)}n
\sum_{q=1}^{n}\sum_{j=1}^{d}\left\vert D_{(q,j)}A_r(V_q)\right\vert ^{2}.
\end{equation*}%
We notice that $D_{(q,j)}A_r(V_q)$ is not null for $r<|V_q-y_{Y_q}|^2<2r$
and contains the derivatives of $a_r$ up to order $2$, possibly multiplied
by polynomials in the components of $V_q-y_{Y_q}$ of order up to 2. Since $%
|V_q-y_{Y_q}|^2\leq 2r$, by using (\ref{s3}) one obtains ${\mathbb{E}}%
(\left\vert DLS_n(Y)\right\vert _{l_{2}}^{p})\leq Cr^{-2p}\times {\mathrm{C}}%
_{2}^{p/2}(Y),$ so (\ref{s14}) holds for $k=1$ also. And for higher order
derivatives the proof is similar. $\square $

\begin{remark}
For further use, we give here an upper estimate of the quantity $\|\mathcal{K%
}_{q,0}(F)\|_p$, with $\mathcal{K}_{q,0}(F)$ defined in (\ref{FD15'}), in
the case $F=S_n(Y)$ (recall that $S_n(Y)$ takes values in ${\mathbb{R}}^d$).
From (\ref{FD15'}), it follows that 
\begin{equation*}
\|\mathcal{K}_{q,0}(F)\|_p \leq \||F|\|_{q,2qp}^{q}
(1+\|F\|_{q+1,8dqp})^{4dq}.
\end{equation*}
So, for $F=S_n(Y)$ we use (\ref{s15}) and for a suitable constant $C$
depending on $d$, $q$ and $p$ only, we obtain 
\begin{equation}  \label{CqSn}
\|\mathcal{K}_{q,0}(S_n(Y))\|_p \leq C\times \frac{{\mathrm{C}}%
_{8dqp}^{(4d+1)q}(Y)}{r^{q(q+1)}}.
\end{equation}
\end{remark}

We give now estimates of the Malliavin covariance matrix. We have 
\begin{equation*}
\sigma _{S_{n}(Y)}=\frac 1n\sum_{k=1}^{n}\chi _{k}\sigma _{n,k}.
\end{equation*}%
We denote 
\begin{equation}  \label{s17-0}
\Sigma_{n}=\frac 1n\sum_{i=1}^{n}\sigma _{n,i},\quad \underline{\lambda }%
_{n}=\inf_{\left\vert \xi \right\vert =1}\left\langle \Sigma_{n}\xi ,\xi
\right\rangle ,\quad \overline{\lambda }_{n}=\sup_{\left\vert \xi
\right\vert =1}\left\langle \Sigma_{n}\xi ,\xi \right\rangle .
\end{equation}
For reasons which will be clear later on, we \textit{do not} consider here
the normalization condition $\Sigma_{n}=\mathrm{Id}_d$. We have the
following result.

\begin{lemma}
Let $\eta=(\frac{\underline{\lambda}_n{\mathfrak{m}}_r }{2(1+2\overline{%
\lambda}_n)})^d$, $\overline{\lambda}_n$ and $\underline{\lambda}_n$ being
given in (\ref{s17-0}). Then%
\begin{equation}
{\mathbb{P}}(\det \sigma _{S_{n}(Y)}\leq \eta ) \leq \frac{e^{3}\overline{c}%
_d}{9}\Big(\frac{2(1+2 \overline{\lambda}_n)}{\underline{\lambda}_n{%
\mathfrak{m}}_r }\Big)^d \exp \Big(-\frac{\underline{\lambda }_{n}^{2}{%
\mathfrak{m}}^2_r }{16\overline{\lambda }_{n}}\times n\Big),  \label{s17}
\end{equation}
$\overline{c}_d$ denoting a positive constant depending on the dimension $d$
only.
\end{lemma}

\textbf{Proof}. Since $\sigma _{n,k}=C_{n,k}C_{n,k}^{\ast }$ we have 
\begin{equation*}
\langle \sigma _{S_{n}(Y)}\xi ,\xi \rangle =\frac 1n\sum_{k=1}^{n}\chi
_{k}\left\langle \sigma _{n,k}\xi ,\xi \right\rangle =\frac
1n\sum_{k=1}^{n}\chi _{k}\left\vert C_{n,k}\xi \right\vert ^{2}.
\end{equation*}%
Take $\xi _{1},...,\xi _{N}\in S_{d-1}=:\{\xi \in {\mathbb{R}}^{d}:\left|
\xi \right| =1\}$ such that the balls of centers $\xi _{i}$ and radius $\eta
^{1/d}$ cover $S_{n-1}.$ One needs $N\leq \bar c_{d} \eta ^{-1}$ points,
where $\bar c_{d}$ is a constant depending on the dimension. It is easy to
check that $\xi \mapsto \langle \sigma _{S_{n}(Z)}\xi ,\xi \rangle $ is
Lipschitz continuous with Lipschitz constant $2\overline{\lambda }_{n}$ so
that $\inf_{\left| \xi \right| =1}\<\sigma _{S_{n}(Z)} \xi ,\xi \>\geq
\inf_{i=1,\ldots,N}\<\sigma _{S_{n}(Z)} \xi _{i},\xi _{i}\>-2\overline{%
\lambda}_n\eta^{1/d}$. Consequently, 
\begin{align*}
{\mathbb{P}}(\det \sigma _{S_{n}(Z)} \leq \eta ) &\leq {\mathbb{P}}%
(\inf_{\left\vert \xi \right\vert =1}\left\langle \sigma _{S_{n}(Z)}\xi ,\xi
\right\rangle \leq \eta ^{1/d})\leq \sum_{i=1}^N {\mathbb{P}}(\left\langle
\sigma _{S_{n}(Z)}\xi_i ,\xi_i \right\rangle \leq \eta ^{1/d}+2\overline{%
\lambda}_n\eta^{1/d}) \\
&\leq \frac{\bar c_{d}}{\eta} \max_{i=1,\ldots,N}{\mathbb{P}}\big(\langle
\sigma _{S_{n}(Z)}\xi_i ,\xi_i \rangle \leq \eta ^{1/d}(1+2\overline{\lambda}%
_n)\big).
\end{align*}
So, it remains to prove that for every $\xi\in S_{d-1}$ and for the choice $%
\eta=(\frac{\underline{\lambda}_n{\mathfrak{m}}_r }{2(1+2\overline{\lambda}%
_n)})^d$, 
\begin{equation*}
{\mathbb{P}}\big(\langle \sigma _{S_{n}(Z)}\xi ,\xi \rangle \leq (1+2%
\overline{\lambda}_n)\eta ^{1/d}) \leq\frac{2e^{3}}{9}\exp \Big(-\frac{%
\underline{\lambda }_{n}^{2}{\mathfrak{m}}^2_r }{16\overline{\lambda }_{n}}%
\times n\Big).
\end{equation*}

We recall that ${\mathbb{E}}(\chi _{k})={\mathfrak{m}}_r $ and we write%
\begin{align*}
{\mathbb{P}}\big(\langle \sigma _{S_{n}(Z)}\xi ,\xi \rangle \leq (1+2%
\overline{\lambda}_n)\eta ^{1/d}) &={\mathbb{P}}\Big(\frac
1n\sum_{k=1}^{n}(\chi _{k}-{\mathfrak{m}}_r )\left\vert C_{n,k}\xi
\right\vert ^{2}\leq (1+2\overline{\lambda}_n)\eta ^{1/d}-{\mathfrak{m}}_r
\frac 1n\sum_{k=1}^{n}\left\vert C_{n,k}\xi \right\vert ^{2}\Big) \\
&\leq {\mathbb{P}}\Big(-\frac 1n\sum_{k=1}^{n}(\chi _{k}-{\mathfrak{m}}_r
)\left\vert C_{n,k}\xi \right\vert ^{2}\geq \underline{\lambda }_{n}{%
\mathfrak{m}}_r -(1+2\overline{\lambda}_n)\eta ^{1/d}\Big)
\end{align*}%
the last equality being true because, from (\ref{s17-0}), 
\begin{equation*}
\frac 1n\sum_{k=1}^{n}\left\vert C_{k}\xi \right\vert ^{2}=\frac
1n\sum_{k=1}^{n}\left\langle \sigma _{n,k}\xi ,\xi \right\rangle
=\<\Sigma_n\xi,\xi\>\geq \underline{\lambda }_{n}\left\vert \xi \right\vert
^{2}=\underline{\lambda }_{n}.
\end{equation*}%
So, we take $\eta=(\frac{\underline{\lambda}_n{\mathfrak{m}}_r }{2(1+2%
\overline{\lambda}_n)})^d$ and we get 
\begin{align*}
{\mathbb{P}}\big(\langle \sigma _{S_{n}(Z)}\xi ,\xi \rangle \leq (1+2%
\overline{\lambda}_n)\eta ^{1/d}) &\leq {\mathbb{P}}\Big(-\sum_{k=1}^{n}(%
\chi _{k}-{\mathfrak{m}}_r )\left\vert C_{n,k}\xi \right\vert ^{2}\geq \frac{%
\underline{\lambda }_{n}{\mathfrak{m}}_r }2\times n\Big)
\end{align*}%
We now use the following Hoeffding's inequality (in the slightly more
general form given in \cite{[Ben]} Corollary 1.4): if the differences $X_k$
of a martingale $M_n$ are such that ${\mathbb{P}}(|X_k|\leq b_k)=1$ then ${%
\mathbb{P}}(M_n\geq x)\leq (2e^3/9)\exp(-|x|^2\times
n/(2(b_1^2+\cdots+b_n^2)))$. Here, we choose $X_{k}=-(\chi _{k}-{\mathfrak{m}%
}_r )\left\vert C_{n,k}\xi \right\vert ^{2}$. These are independent random
variables and $\left\vert X_{k}\right\vert \leq 2\left\vert C_{n,k}\xi
\right\vert ^{2}.$ Then 
\begin{align*}
{\mathbb{P}}\Big(-\sum_{k=1}^{n}(\chi _{k}-{\mathfrak{m}}_r )\left\vert
C_{n,k}\xi \right\vert ^{2} \geq \frac{\underline{\lambda}_{n}{\mathfrak{m}}%
_r }2\times n\Big) &\leq \frac{2e^{3}}{9}\exp \Big(-\frac{\underline{\lambda 
}_{n}^{2}{\mathfrak{m}}^2_r }{4}\times \frac{n}{4\sum_{k=1}^{n}\left\vert
C_{n,k}\xi \right\vert ^{2}}\Big) \\
&\leq\frac{2e^{3}}{9}\exp \Big(-\frac{\underline{\lambda }_{n}^{2}{\mathfrak{%
m}}^2_r }{16\overline{\lambda }_{n}}\times n\Big).
\end{align*}%
$\square $

\bigskip

We are now able to give the regularization lemma in our specific framework.

\begin{lemma}
\label{RR} Let $h ,q \in {\mathbb{N}}.$ There exists a constant $C\geq 1,$
depending just on $h ,q $, such that for every $\delta >0,$\ every
multiindex $\gamma $ with $\left\vert \gamma \right\vert =q $\ and every $%
f\in C_{p}^{q }({\mathbb{R}}^{d})$ one has%
\begin{equation}  \label{s16}
\begin{array}{ll}
& \displaystyle \left\vert {\mathbb{E}}(\partial_{\gamma }f(x+S_{n}(Y)))-{%
\mathbb{E}}(\partial_{\gamma }f_{\delta }(x+S_{n}(Y)))\right\vert\smallskip
\\ 
& \leq \displaystyle C{\mathrm{C}}_{2l_{0}(f)}^{1/2}(Y) {\mathrm{Q}}_{h ,q
}(Y) \Big(L_{q }(f) \exp \Big(-\frac{\underline{\lambda }_{n}^{2}{\mathfrak{m%
}}^2_r }{32\overline{\lambda }_{n}}\times n\Big) +L_{0}(f)\delta ^{h }\Big)%
(1+|x|)^{l_q (f)}%
\end{array}%
\end{equation}
with%
\begin{equation}
{\mathrm{Q}}_{h ,q }(Y)=\,2^{l_q (f)}c_{l_q (f)\vee(l_0(f)+h )}\frac{{%
\mathrm{C}}_{16d(h +q )}^{(4d+1)(h +q )}(Y)}{r^{(h +q )(h +q +1)}} \Big(%
1\vee \frac{2(1+2 \overline{\lambda}_n)}{\underline{\lambda}_n{\mathfrak{m}}%
_r }\Big)^{2d(h +q )},  \label{s16'}
\end{equation}
$c_p$ being given in (\ref{1aa}).
\end{lemma}

\textbf{Proof}. We will use Lemma \ref{R-1}. Since $\mathcal{C}_{h +q
}(S_n(Y))= \|\mathcal{K}_{h ,0}(S_n(Y))\|_2$, (\ref{CqSn}) gives 
\begin{equation*}
\mathcal{C}_{h +q }(S_n(Y)) \leq C\times \frac{{\mathrm{C}}_{16d(h +q
)}^{(4d+1)(h +q )}(Y)}{r^{(h +q )(h +q +1)}},
\end{equation*}
$C$ depending on $d$ and $h +q $. And by using the Burkholder inequality (%
\ref{b}), one has $\left\Vert S_{n}(Y)\right\Vert
_{2l_{0}(f)}^{l_{0}(f)}\leq {\mathrm{C}}_{2l_{0}(f)}^{1/2}(Y).$ So (\ref{1b}%
) (with $p=2$) gives%
\begin{equation*}
\begin{array}{l}
\displaystyle \left\vert {\mathbb{E}}(\partial _{\gamma }f(x+S_n(Y)))-{%
\mathbb{E}}(\partial _{\gamma }f_{\delta }(x+S_n(Y)))\right\vert \leq C{%
\mathrm{C}}_{2l_{0}(f)}^{1/2}(Y) \,2^{l_q (f)}c_{l_q (f)\vee(l_0(f)+h )}%
\frac{{\mathrm{C}}_{16d(h +q )}^{(4d+1)(h +q )}(Y)}{r^{(h +q )(h +q +1)}}%
\smallskip \\ 
\displaystyle \quad \times \Big(L_{q }(f) {\mathbb{P}}^{1/2}(\det%
\sigma_{S_n(Y)}\leq \eta)+L_{0}(f)\frac{\delta ^{h }}{\eta ^{2(h +q )}}\Big)%
(1+|x|)^{l_q (f)}.%
\end{array}%
\end{equation*}
We take now $\eta =(\frac{\underline{\lambda }_{n}{\mathfrak{m}}_r }{2(1+%
\overline{\lambda}_n)})^{d}$ and we use (\ref{s17}) in order to obtain 
\begin{equation*}
\begin{array}{l}
\displaystyle \left\vert {\mathbb{E}}(\partial _{\gamma }f(x+S_n(Y)))-{%
\mathbb{E}}(\partial _{\gamma }f_{\delta }(x+S_n(Y)))\right\vert \leq C{%
\mathrm{C}}_{2l_{0}(f)}^{1/2}(Y) \,2^{l_q (f)}c_{l_q (f)\vee(l_0(f)+h )}%
\frac{{\mathrm{C}}_{16d(h +q )}^{(4d+1)(h +q )}(Y)}{r^{(h +q )(h +q +1)}}%
\smallskip \\ 
\displaystyle \quad \times \Big(1\vee \frac{2(1+2 \overline{\lambda}_n)}{%
\underline{\lambda}_n{\mathfrak{m}}_r }\Big)^{2d(h +q )}\Big(L_{q }(f) \exp %
\Big(-\frac{\underline{\lambda }_{n}^{2}{\mathfrak{m}}^2_r }{32\overline{%
\lambda }_{n}}\times n\Big) +L_{0}(f)\delta ^{h }\Big)(1+|x|)^{l_q (f)}.%
\end{array}%
\end{equation*}
$\square $

\medskip

We are now able to characterize the regularity of the semigroup $%
P_{n}^{Z,n}: $

\begin{proposition}
Let $f\in C_p^q ({\mathbb{R}}^d)$. If $\left\vert \gamma \right\vert =q $
then 
\begin{equation}  \label{s16a}
\begin{array}{rl}
\displaystyle \left\vert {\mathbb{E}}\big(\partial_{\gamma }f(x+S_{n}(Y))%
\big)\right\vert \leq & C\times 2^{l_q (f)}{\mathrm{B}}_q (Y)(1+{\mathrm{C}}%
_{2l_q (f)}^{l_q (f)}(Y)) (1+|x|)^{l_q (f)}\times \smallskip \\ 
& \times \displaystyle \Big[L_q (f)\exp \Big(-\frac{\underline{\lambda }%
_{n}^{2}{\mathfrak{m}}^2_r }{32\overline{\lambda }_{n}}\times n\Big) +L_0(f) %
\Big]%
\end{array}%
\end{equation}
where 
\begin{equation}  \label{Bm}
{\mathrm{B}}_q (Y)=\Big(1\vee\frac{2(1+\overline{\lambda}_n)}{\underline{%
\lambda }_{n}{\mathfrak{m}}_r }\Big)^{2dq }\,\frac{{\mathrm{C}}_{16dq
}^{(4d+1)q }(Y)}{r^{q (q +1)}}
\end{equation}
and $C$ is a constant depending on $q $ only.
\end{proposition}

\textbf{Proof}. We take $\eta =(\frac{\underline{\lambda }_{n}{\mathfrak{m}}%
_r }{2(1+\overline{\lambda}_n)})^{d}$ and the truncation function $\Psi
_{\eta }$ and we write%
\begin{equation*}
{\mathbb{E}}(\partial_{\gamma }f(x+S_{n}(Y)))=I+J
\end{equation*}%
with%
\begin{equation*}
I={\mathbb{E}}(\partial_{\gamma }f(x+S_{n}(Y))(1-\Psi _{\eta }(\det \sigma
_{S_{n}}))),\quad J={\mathbb{E}}(\partial_{\gamma }f(x+S_{n}(Y))\Psi _{\eta
}(\det \sigma _{S_{n}})).
\end{equation*}%
We estimate first%
\begin{eqnarray*}
\left\vert I\right\vert &\leq &L_{q }(f){\mathbb{E}}((1+|x|+\left\vert
S_{n}(Y)\right\vert )^{l_{q }(f)}(1-\Psi _{\eta }(\det \sigma _{S_{n}}))) \\
&\leq &L_{q }(f)\big({\mathbb{E}}\big((1+|x|+\left\vert S_{n}(Y)\right\vert
^{2l_{q }(f)})\big)\big)^{1/2}{\mathbb{P}}^{1/2}(\det \sigma _{S_{n}}\leq
\eta ) \\
&\leq &CL_{q }(f)2^{l_q (f)}(1+|x|) ^{l_{q }(f)}(1+{\mathrm{C}}_{2l_{q
}(f)}^{l_{q }(f)}(Y)) \Big(\frac{2(1+2\overline{\lambda}_n)}{\underline{%
\lambda}_n{\mathfrak{m}}_r }\Big)^{d/2} \exp \Big(-\frac{\underline{\lambda }%
_{n}^{2}{\mathfrak{m}}^2_r }{32\overline{\lambda }_{n}}\times n\Big),
\end{eqnarray*}%
in which we have used the Burkholder's inequality (\ref{b}). In order to
estimate $J$ we use integration by parts and we obtain%
\begin{eqnarray*}
|J| &=&\big|{\mathbb{E}}\big(f(x+S_{n}(Y))H_{\gamma }(S_{n}(Y),\Psi _{\eta
}(\det \sigma _{S_{n}}))\big)\big| \\
&\leq &L_{0}(f){\mathbb{E}}\big((1+|x|+\left\vert S_{n}\right\vert)
^{l_{0}(f)}\left\vert H_{\gamma }(S_{n}(Y),\Psi _{\eta }(\det \sigma
_{S_{n}}))\right\vert \big) \\
&\leq &CL_{0}(f)2^{l_0(f)}(1+|x|)^{l_0(f)}(1+{\mathrm{C}}%
_{2l_{0}(f)}^{l_0(f)}(Y))\big({\mathbb{E}}(\left\vert H_{\gamma
}(S_{n}(Y),\Psi _{\eta }(\det \sigma _{S_{n}}))\right\vert ^{2})\big)^{1/2}.
\end{eqnarray*}%
Then using (\ref{FD4}) and (\ref{CqSn})%
\begin{align*}
\left\Vert H_{\alpha }(S_n(Y),\Psi _{\eta }(\det \sigma
_{S_n(Y)}))\right\Vert _{2} &\leq C\times \Big(\frac{2(1+\overline{\lambda}%
_n)}{\underline{\lambda }_{n}{\mathfrak{m}}_r }\Big)^{2dq }\times \frac{{%
\mathrm{C}}_{16dq }^{(4d+1)q }(Y)}{r^{q (q +1)}}(1+|x|)^{l_0(f)},
\end{align*}
so that 
\begin{equation*}
|J| \leq C\times L_{0}(f)2^{l_0(f)}(1+|x|)^{l_0(f)}(1+{\mathrm{C}}%
_{2l_{0}(f)}^{l_0(f)}(Y)) \times \Big(\frac{2(1+\overline{\lambda}_n)}{%
\underline{\lambda }_{n}{\mathfrak{m}}_r }\Big)^{2dq }\times \frac{{\mathrm{C%
}}_{16dq }^{(4d+1)q }(Y)}{r^{q (q +1)}}(1+|x|)^{l_0(f)}.
\end{equation*}
(\ref{s16a}) now follows from the above estimates for $I$ and $J$. $\square $

\medskip

We are now able to give the proof of the main result Theorem \ref{MAIN}.

\subsection{Proofs of the results in Section \ref{main-results}}

\label{proofs}

\subsubsection{Proof of Theorem \ref{MAIN}}

\label{sect-proof-MAIN}

\textbf{Step 1}. We assume first that $f\in C_{p}^{q+(N+1)(N+3)}({\mathbb{R}}%
^{d})$ and we prove that

\begin{equation}  \label{s20bis}
\begin{array}{l}
\displaystyle \left\vert {\mathbb{E}}(\partial _{\gamma }f(S_{n}(Y)))-{%
\mathbb{E}}(\partial _{\gamma }f(W)\Phi _{N}(W))\right\vert\smallskip \\ 
\displaystyle \quad \leq \frac C{n^{\frac {N+1}2}}\times \widehat {\mathrm{C}%
}_{q+(N+1)+(N+3),N}(f,Y) \Big[L_{q+(N+1)(N+3)}(f)\,e^{-\frac{n}{128}}
+L_0(f) \Big],%
\end{array}%
\end{equation}
where $C$ is a constant depending only on $q$ and $N$ and 
\begin{equation}  \label{s20bis'}
\begin{array}{rl}
\displaystyle \widehat{\mathrm{C}}_{p,N}(f,Y) = & {\mathrm{C}}_2^{(\lfloor
N/2\rfloor +1)(N+1)}(Y) 2^{(N+3)l_{p}(f)}{\mathrm{C}}_{2(N+3)}^{(N+1)/2}(Y)%
\big(1+{\mathrm{C}}_{2l_{p}(f)}^{l_{p}(f)\vee (N+1)}(Y)\big)\times\smallskip
\\ 
\displaystyle & \displaystyle \times\Big(1\vee\frac{8}{{\mathfrak{m}}_r }%
\Big)^{2dp}\,\frac{{\mathrm{C}}_{16dp}^{(4d+1)p}(Y)}{r^{p(p+1)}}.%
\end{array}%
\end{equation}
Notice that (\ref{s20bis}) is analogous to (\ref{s20}) but here $L_{q}(f)$
and $l_q(f)$ are replaced by $L_{q+(N+1)(N+3)}(f)$ and $l_{q+(N+1)(N+3)}(f)$%
. We will see in Step 2 how to drop the dependence on $q+(N+1)(N+3)$.

We recall (\ref{n10}) and (\ref{n11}): we have 
\begin{align*}
&{\mathbb{E}}\big(\partial_\gamma f(S_n(Y))\big) -{\mathbb{E}}\big(%
\partial_\gamma f(W)\Phi_{n,N}(W)\big) =P_{n}^{Z,n}(\partial _{\gamma
}f)(0)-P_{n}^{G,n}\Big(\mathrm{Id}+\sum_{k=1}^{N}\frac 1{n^{k/2}}\Gamma
_{n,k}\Big)(\partial _{\gamma }f)(0) \\
&=I_{1}(\partial _{\gamma }f)(0)+I_{2}(\partial _{\gamma
}f)(0)+I_{3}(\partial _{\gamma }f)(0)
\end{align*}%
with%
\begin{equation}  \label{s21}
\begin{array}{l}
\displaystyle I_{1} =P_{n}^{G,n}Q_{n,N}^{0}, \smallskip \\ 
\displaystyle I_{2} =\sum_{1\leq r_{1}<\cdots <r_{N+1}\leq
n}P_{r_{N+1}+1,n+1}^{Z,n}P_{r_{N}+1,r_{N+1}}^{G,n}\cdots
P_{r_{1}+1,r_{2}}^{G,n}P_{1,r_{1}}^{G,n}Q_{n,N,r_{1},...,r_{N+1}}^{(N+1)}
\smallskip \\ 
\displaystyle I_{3} =\sum_{m=1}^{N}\sum_{1\leq r_{1}<\cdots <r_{m}\leq
n}P_{r_{m}+1,n+1}^{G,n}P_{r_{m-1}+1,r_{m}}^{G,n}\cdots
P_{r_{1}+1,r_{2}}^{G,n}P_{1,r_{1}}^{G,n}Q_{n,N,r_{1},...,r_{m}}^{(m)}%
\end{array}%
\end{equation}
and (see (\ref{n1}) and (\ref{n9}))

\begin{equation}  \label{s22}
\begin{array}{l}
\displaystyle Q_{N,r_{1},...,r_{m}}^{(m)} =\frac 1{n^{\frac {N+3m}2}}
\sum_{3\leq \left\vert \alpha \right\vert \leq
N_{m}}a_{n,r_{1},...,r_{m}}(\alpha )\theta _{n,r_{1},...,r_{m}}^{\alpha
}\partial_{\alpha }, \smallskip \\ 
\displaystyle Q_{N,n}^{0}f(x) =\frac 1{n^{(N+1)/2}}\sum_{N+1\leq \left\vert
\alpha \right\vert \leq N(N+2\lfloor N/2\rfloor)}c_{n,N}(\alpha
)\partial_{\alpha },%
\end{array}%
\end{equation}
$N_m$ being given in Lemma \ref{Nm}: $N_{m}=m(2\lfloor N/2\rfloor+N+5)$ for $%
m\leq N$ and if $m=N+1$ then $N_{m}=(N+1)(N+3)$. By (\ref{n2}) and (\ref{n9}%
), the coefficient which appear above satisfy 
\begin{equation}  \label{s23}
\left\vert a_{n,r_{1},...,r_{m}}(\alpha )\right\vert \leq (C{\mathrm{C}}%
_2^{\lfloor N/2\rfloor+1}(Y))^m \quad\mbox{and}\quad |c_{n,N}(\alpha )| \leq
(C{\mathrm{C}}_{N+1}(Y){\mathrm{C}}_2(Y))^{N(N+2\lfloor N/2\rfloor)}.
\end{equation}%
We first estimate $I_{2}(\partial_\gamma f).$ Let us prove that for every $%
r_{1}<\cdots<r_{N+1}$ 
\begin{equation}  \label{s26}
\begin{array}{l}
\displaystyle \left\vert
P_{r_{N+1}+1,n}^{Z,n}P_{r_{N}+1,r_{N+1}}^{G,n}...P_{r_{1}+1,r_{2}}^{G,n}P_{k,r_{1}}^{G,n}Q_{r_{1},...,r_{N+1}}^{(N+1)}\partial _{\gamma }f(x)\right\vert\leq \frac C{n^{%
\frac {4N+3}2}}\times \smallskip \\ 
\displaystyle \times \widehat {\mathrm{C}}_{q+(N+1)+(N+3),N}(f,Y) \Big[%
L_{q+(N+1)(N+3)}(f)\,e^{-\frac{{\mathfrak{m}}^2_r}{128}\, n} +L_0(f) \Big] %
(1+|x|)^{l_{q+(N+1)(N+3)}(f)}%
\end{array}%
\end{equation}
where $C$ depends only on $q$ and $N$ and $\widehat{\mathrm{C}}_{p,N}(f,Y)$
is given by (\ref{s20bis'}).

Recall that $\sigma _{n,r_{i}}\leq \frac{1}{n}{\mathrm{C}}_{2}(Y).$ We take $%
n\geq 4(N+1){\mathrm{C}}_{2}(Y)$ so that%
\begin{equation}
\frac 1n\sum_{i=1}^{N+1}\sigma _{n,r_{i}}\leq \frac{1}{4}\mathrm{Id}_d.
\label{s24'}
\end{equation}%
Recall that $\frac 1n\sum_{r=1}^{n}\sigma _{n,r}=\mathrm{Id}_d.$ So we
distinguish now two cases: 
\begin{align}
\mbox{\textbf{Case 1:}}\quad\quad & \displaystyle \frac
1n\sum_{r=r_{N+1}+1}^{n}\sigma _{n,r}\geq \frac{1}{2}\mathrm{Id}_d,
\smallskip  \label{s24} \\
\mbox{\textbf{Case 2:}}\quad\quad & \displaystyle \frac
1n\sum_{r=1}^{r_{N+1}}\sigma _{n,r} \geq \frac{1}{2}\mathrm{Id}_d.
\label{s25}
\end{align}
We treat Case 1. Notice that all the operators coming on in (\ref{s21})
commute so, using also (\ref{s22}) we obtain 
\begin{eqnarray*}
&&P_{r_{N+1}+1,n+1}^{Z,n}P_{r_{N}+1,r_{N+1}}^{G,n}...P_{r_{1}+1,r_{2}}^{G,n}P_{k,r_{1}}^{G,n}Q_{N,r_{1},...,r_{N+1}}^{(N+1)}\partial _{\gamma }f(x)
\\
&=&\frac 1{n^{(4N+3)/2}}\!\!\!\!\sum_{3\leq \left\vert \alpha \right\vert
\leq (N+1)(N+3)}\!\!\!\!a_{n,r_{1},\ldots ,r_{N+1}}(\alpha )\theta
_{n,r_{1},\ldots,r_{N+1}}^{\alpha }P_{r_{N-1}+1,r_{N}}^{G,n}\cdots
P_{r_{1}+1,r_{2}}^{G,n}P_{1,r_{1}}^{G,n}P_{r_{N}+1,n}^{Z,n}\partial _{\gamma
}\partial _{\alpha }f(x).
\end{eqnarray*}%
We use now (\ref{s16a}) with $m=\left\vert \gamma \right\vert +\left\vert
\alpha \right\vert \leq q+(N+1)(N+3)$ and $S_n(Y)=\frac 1{\sqrt
n}\sum_{k=1}^nC_{n,k}Y_k\equiv \sum_{k=1}^nZ_{n,k}$ replaced by $%
\sum_{k=r_{N+1}+1}^nZ_{n,k}$, whose covariance matrix is $\frac
1n\sum_{k=r_{N+1}+1}^n\sigma_{n,r}$. Under (\ref{s24}) we have $\frac 12\leq 
\underline{\lambda}_n\leq \overline{\lambda}_n\leq 1$, so (\ref{s16a}) gives 
\begin{equation*}
\begin{array}{rl}
\displaystyle |P_{r_{N+1}+1,n+1}^{Z,n}\partial _{\gamma }\partial _{\alpha
}f(x)| \leq & C\times 2^{l_{q+(N+1)(N+3)}(f)}\widehat{\mathrm{B}}%
_{q+(N+1)(N+3)}(Y) (1+{\mathrm{C}}%
_{2l_{q+(N+1)(N+3)}(f)}^{l_{q+(N+1)(N+3)}(f)}(Y))\times \smallskip \\ 
& \times \displaystyle \Big[L_{q+(N+1)(N+3)}(f)\,e^{-\frac{{\mathfrak{m}}_r^2%
}{128}\,n} +L_0(f) \Big](1+|x|)^{l_{q+(N+1)(N+3)}(f)}%
\end{array}%
\end{equation*}
where $C$ is a constant depending only on $q$ and $N$ and $\widehat{\mathrm{B%
}}_p(Y)$ is the constant ${\mathrm{B}}_p(Y)$ given in (\ref{Bm}) with $%
\underline{\lambda}_n=1/2$ and $\overline{\lambda}_n=1$, that is, 
\begin{equation*}
\widehat{\mathrm{B}}_p(Y)=\Big(1\vee\frac{8}{{\mathfrak{m}}_r }\Big)^{2dp}\,%
\frac{{\mathrm{C}}_{16dp}^{(4d+1)p}(Y)}{r^{p(p+1)}}
\end{equation*}
Therefore, we can write 
\begin{align*}
l_0\big(P_{r_{N+1}+1,n}^{Z,n}\partial _{\gamma }\partial _{\alpha }f\big)%
=&l_{q+(N+1)(N+3)}(f), \\
L_0\big(P_{r_{N+1}+1,n}^{Z,n}\partial _{\gamma }\partial _{\alpha }f\big)%
=&C\times 2^{l_{q+(N+1)(N+3)}(f)}\widehat{\mathrm{B}}_{q+(N+1)(N+3)}(Y) (1+{%
\mathrm{C}}_{2l_{q+(N+1)(N+3)}(f)}^{l_{q+(N+1)(N+3)}(f)}(Y))\times \\
& \times \displaystyle \Big[L_{q+(N+1)(N+3)}(f)\,e^{-\frac{{\mathfrak{m}}_r^2%
}{128}\,n} +L_0(f) \Big].
\end{align*}

Now, in the proof of Theorem \ref{Smooth} we have proven that (see (\ref{gl}%
) 
\begin{equation*}
|P_{r_{N-1}+1,r_{N}}^{G,n}\cdots
P_{r_{1}+1,r_{2}}^{G,n}P_{1,r_{1}}^{G,n}g(x)| \leq 2^{l_0(g)}(1+{\mathrm{C}}%
_{l_0(g)}(Y))L_0(g)(1+|x|)^{l_0(g)}
\end{equation*}
and following the proof of Lemma \ref{Nm} we have 
\begin{equation*}
|\theta^\alpha_{n,r_1,\ldots,r_m}g(x)| \leq \big(2^{l_{0}(g)}{\mathrm{C}}%
_{2(N+3)}^{1/2}(Y)(1+{\mathrm{C}}_{2{l_{0}(g)}}(Y))^2\big)^m
L_{0}(g)(1+|x|)^{l_{0}(g)}.
\end{equation*}
So, taking all estimates, we obtain 
\begin{align*}
&|\theta _{n,r_{1},\ldots,r_{N+1}}^{\alpha }P_{r_{N-1}+1,r_{N}}^{G,n}\cdots
P_{r_{1}+1,r_{2}}^{G,n}P_{1,r_{1}}^{G,n}P_{r_{N}+1,n+1}^{Z,n}\partial
_{\gamma }\partial _{\alpha }f(x)| \\
&\leq C\times {\mathrm{D}}_{q+(N+1)(N+3)}(f,Y) \Big[L_{q+(N+1)(N+3)}(f)\,e^{-%
\frac{{\mathfrak{m}}_r^2}{128}\,n} +L_0(f) \Big] %
(1+|x|)^{l_{q+(N+1)(N+3)}(f)}
\end{align*}
where 
\begin{equation*}
{\mathrm{D}}_p(f,Y)=2^{(N+3)l_{p}(f)}{\mathrm{C}}_{2(N+3)}^{(N+1)/2}(Y)\big(%
1+{\mathrm{C}}_{2l_{p}(f)}^{l_{p}(f)\vee (N+1)}(Y)\big)\Big(1\vee\frac{8}{{%
\mathfrak{m}}_r }\Big)^{2dp}\,\frac{{\mathrm{C}}_{16dp}^{(4d+1)p}(Y)}{%
r^{p(p+1)}}
\end{equation*}
We use now formula (\ref{s22}) for $Q^{(N+1)}_{N,r_1,\ldots, r_{N+1}}$ and
the estimate (\ref{s23}) for the coefficients $a_{n,r_1,\ldots, r_{N+1}}$
and we get 
\begin{align*}
&|P_{r_{N+1}+1,n+1}^{Z,n}P_{r_{N}+1,r_{N+1}}^{G,n}...P_{r_{1}+1,r_{2}}^{G,n}P_{1,r_{1}}^{G,n}Q_{r_{1},...,r_{N+1}}^{(N+1)}\partial _{\gamma }f(x)|
\\
\leq &C\frac 1{n^{\frac {4N+3}2}}\times \widehat {\mathrm{C}}_\ast(f,Y) \Big[%
L_{q+(N+1)(N+3)}(f)\,e^{-\frac{{\mathfrak{m}}_r^2}{128}\,n} +L_0(f) \Big] %
(1+|x|)^{l_{q+(N+1)(N+3)}(f)}
\end{align*}
where $C$ depends only on $q$ and $N$ and $\widehat{\mathrm{C}}_\ast(f,Y)$
is given by 
\begin{equation*}
\widehat{\mathrm{C}}_\ast(f,Y)=\big({\mathrm{C}}_2^{\lfloor N/2\rfloor +1}(Y)%
\big)^{N+1} {\mathrm{D}}_{q+(N+1)(N+3)}(f,Y) .
\end{equation*}
Since $\widehat{\mathrm{C}}_\ast(f,Y)=\widehat{\mathrm{C}}%
_{q+(N+1)(N+3),N}(f,Y)$, (\ref{s26}) is proved in Case 1.

We deal now with Case 2, that is, we assume (\ref{s25}). We write%
\begin{align*}
&P_{r_{N+1}+1,n+1}^{Z,n}P_{r_{N}+1,r_{N+1}}^{G,n}...P_{r_{1}+1,r_{2}}^{G,n}P_{1,r_{1}}^{G,n}Q_{N,r_{1},...,r_{N+1}}^{(N+1)}\partial _{\gamma }f(x)
\\
&=\frac 1{n^{\frac {4N+3}2}}\sum_{3\leq \left\vert \alpha \right\vert \leq
(N+1)^{2}}a_{n, r_{1},\ldots ,r_{N+1}}(\alpha )\theta _{n,r_{1},\ldots
,r_{N+1}}^{\alpha }P_{r_{N+1}+1,n+1}^{Z,n}P_{r_{N}+1,r_{N+1}}^{G,n}\cdots
P_{r_{1}+1,r_{2}}^{G,n}P_{1,r_{1}}^{G,n}\partial _{\gamma }\partial _{\alpha
}f(x).
\end{align*}%
Notice that%
\begin{equation*}
P_{r_{N}+1,r_{N+1}}^{G,n}\cdots
P_{r_{1}+1,r_{2}}^{G,n}P_{1,r_{1}}^{G,n}\partial _{\gamma }\partial _{\alpha
}f(x)={\mathbb{E}}(\partial _{\gamma }\partial _{\alpha }f(x+G))
\end{equation*}%
where $G$ is a centered Gaussian random variable of variance $\frac
1n\sum_{i=1}^{r_{N+1}}\sigma _{n,i}-\frac 1n\sum_{i=1}^{N+1}\sigma
_{n,r_{i}}\geq \frac{1}{4}\mathrm{Id}_d$, as it follows by using also (\ref%
{s24'}). So standard integration by parts yields 
\begin{equation*}
|P_{r_{N}+1,r_{N+1}}^{G,n}...P_{r_{1}+1,r_{2}}^{G,n}P_{1,r_{1}}^{G,n}%
\partial _{\gamma }\partial _{\alpha }f(x)|\leq CL_{0}(f)(1+|x|)^{l_0(f)}.
\end{equation*}%
Now the proof follows as in the previous case. So (\ref{s26}) is proved in
Case 2 as well.

Therefore, by summing over $r_{1}<r_{2}<\cdots<r_{N+1}\leq n$ (giving a
contribution of order $n^{N+1}$), inequality (\ref{s26}) gives 
\begin{align*}
&|I_{2}(\partial _{\gamma }f)(x)| \leq \frac C{n^{\frac {N+1}2}}\times 
\widehat{\mathrm{C}}_\ast(f,Y) \Big[L_{q+(N+1)(N+3)}(f)\,e^{-\frac{{%
\mathfrak{m}}_r^2}{128}\,n} +L_0(f) \Big] (1+|x|)^{l_{q+(N+1)(N+3)}(f)}
\end{align*}%
Exactly as in Case 2 presented above (using standard integration by parts
with respect to the law of Gaussian random variables) we obtain 
\begin{equation*}
|I_{1}(\partial _{\gamma }f)(x)|+|I_{3}(\partial _{\gamma }f)(x)| \leq \frac
C{n^{\frac {N+1}2}}\times \widehat{\mathrm{C}}_\ast(f,Y)
L_{0}(f)(1+|x|)^{l_0(f)}.
\end{equation*}
So, recalling that $\widehat{\mathrm{C}}_\ast(f,Y)=\widehat{\mathrm{C}}%
_{q+(N+1)(N+3),N}(f,Y)$, (\ref{s20bis}) is proved.

\smallskip

\textbf{Step 2.} We now come back and we replace $L_{q+(N+1)(N+3)}(f)$ by $%
L_{q}(f)$ in (\ref{s20bis}). We will use the regularization lemma. So we fix 
$\delta >0$ (to be chosen in a moment) and we write%
\begin{equation*}
\left\vert {\mathbb{E}}(\partial _{\gamma }f(S_{n}(Y)))-{\mathbb{E}}%
(\partial _{\gamma }f(W)\Phi _{N}(W))\right\vert \leq A_{\delta
}(f)+A_{\delta }^{\prime }(f)+A_{\delta }^{\prime \prime }(f)
\end{equation*}%
with 
\begin{align*}
&A_{\delta }(f)=\left\vert {\mathbb{E}}\big(\partial _{\gamma }f_{\delta
}(S_{n}(Y))\big)-{\mathbb{E}}\big(\partial _{\gamma }f_{\delta }(W)\Phi
_{N}(W)\big)\right\vert \\
&A_{\delta }^{\prime }(f) =\left\vert {\mathbb{E}}\big(\partial _{\gamma
}f(S_{n}(Y))\big)-{\mathbb{E}}\big(\partial _{\gamma }f_{\delta }(S_{n}(Y))%
\big)\right\vert , \\
&A_{\delta }^{\prime \prime }(f) =\left\vert {\mathbb{E}}\big(\partial
_{\gamma }f(W)\Phi _{N}(W)\big)-{\mathbb{E}}\big(\partial _{\gamma
}f_{\delta }(W)\Phi _{N}(W)\big)\right\vert .
\end{align*}

We will use (\ref{s20bis}) for $f_{\delta }.$ Notice that $L_p(f_\delta)\leq
\hat c_{p,l_0(f)}L_0(f)\delta^{-p}$, with $\hat c_{p,l}=1\vee \max_{0\leq
|\alpha|\leq p}\int(1+|x|)^l|\partial_\alpha\phi(x)|dx$, and $%
l_p(f_\delta)=l_0(f)$. So, 
\begin{equation*}
A_{\delta }(f)\leq \frac C{n^{\frac {N+1}2}}\times {\mathrm{H}}_{q,N}(f,Y)
L_0(f)\Big[\frac 1{\delta^{q+(N+1)(N+3)}}\,e^{-\frac{{\mathfrak{m}}_r^2}{128}%
\times n} +1\Big],
\end{equation*}
where 
\begin{equation*}
\begin{array}{l}
\displaystyle {\mathrm{H}}_{q,N}(f,Y) ={\mathrm{C}}_2^{(\lfloor N/2\rfloor
+1)(N+1)}(Y) 2^{(N+3)l_{0}(f)}{\mathrm{C}}_{2(N+3)}^{(N+1)/2}(Y)\big(1+{%
\mathrm{C}}_{2l_{0}(f)}^{l_{0}(f)\vee (N+1)}(Y)\big)\times\smallskip \\ 
\quad\displaystyle \times\Big(1\vee\frac{8}{{\mathfrak{m}}_r }\Big)%
^{2d(q+(N+1)(N+3))}\,\frac{{\mathrm{C}}%
_{16d(q+(N+1)(N+3))}^{(4d+1)(q+(N+1)(N+3))}(Y)} {%
r^{(q+(N+1)(N+3))(q+(N+1)(N+3)+1)}} \, \hat c_{q+(N+1)(N+3),l_0(f)}.%
\end{array}%
\end{equation*}
We use now (\ref{s16}) with $x=0$ and with some $h$ to be chosen in a
moment. We then obtain 
\begin{equation*}
A_{\delta }^{\prime }(f)\leq C{\mathrm{C}}_{2l_{0}(f)}^{1/2}(Y) {\mathrm{Q}}%
_{h,q}(Y) \Big(L_{q}(f) e^{-\frac{{\mathfrak{m}}^2_r }{32}\times n}
+L_{0}(f)\delta ^{h}\Big)
\end{equation*}
with ${\mathrm{Q}}_{h,q}(Y)$ given in (\ref{s16'}). And we also have $%
A_{\delta }^{\prime \prime }(f)\leq CL_{0}(f)\delta ^{h}$ (the proof is
identical to the one of (\ref{1b}) but one employs usual integration by
parts with respect to the Gaussian law)$.$ We put all this together and we
obtain 
\begin{align*}
&\left\vert {\mathbb{E}}(\partial _{\gamma }f(S_{n}(Y)))-{\mathbb{E}}%
(\partial _{\gamma }f(W)\Phi _{N}(W))\right\vert \leq\frac C{n^{\frac {N+1}%
2}}\,{\mathrm{H}}_{q,N}(f,Y) L_0(f)\Big[\frac 1{\delta^{q+(N+1)(N+3)}}\,e^{-%
\frac{{\mathfrak{m}}_r^2\,n}{128}} +1\Big] \\
&\qquad +C{\mathrm{C}}_{2l_{0}(f)}^{1/2}(Y) {\mathrm{Q}}_{h,q}(Y) L_{q}(f)
e^{-\frac{{\mathfrak{m}}^2_r\,n }{32}} +CL_{0}(f)\delta ^{h}
\end{align*}
We take now $\delta $ such that%
\begin{equation*}
\delta ^{h}=\frac{1}{\delta ^{q+(N+1)(N+3)}}e^{-\frac{{\mathfrak{m}}^2_r n}{%
128}}
\end{equation*}%
and $h=q+(N+1)(N+3)$, so that 
\begin{equation*}
\delta ^{h}=e^{-\frac{{\mathfrak{m}}^2_r n}{128}\times \frac{h}{%
h+q+(N+1)(N+3)}}=e^{-\frac{{\mathfrak{m}}^2_r n}{256}}.
\end{equation*}%
With this choice of $h$ and $\delta $ we get%
\begin{align*}
&\left\vert {\mathbb{E}}(\partial _{\gamma }f(S_{n}(Y)))-{\mathbb{E}}%
(\partial _{\gamma }f(W)\Phi _{N}(W))\right\vert \leq C\,{\mathrm{H}}%
_{q,N}(f,Y) L_0(f)\Big(\frac 1{n^{\frac {N+1}2}}+ e^{-\frac{{\mathfrak{m}}%
_r^2\,n}{256}}\Big) \\
&\qquad +C{\mathrm{C}}_{2l_{0}(f)}^{1/2}(Y) {\mathrm{Q}}_{q+(N+1)(N+3),q}(Y)
L_{q}(f) e^{-\frac{{\mathfrak{m}}^2_r\,n }{32}}
\end{align*}
We take now $n$ sufficiently large in order to have 
\begin{equation*}
n^{\frac{1}{2}(N+1)}e^{-\frac{{\mathfrak{m}}^2_r n}{256}}\leq 1.
\end{equation*}%
The statement now follows by observing that, with ${\mathrm{C}}_*(Y)$ given
in (\ref{s20'}), 
\begin{equation*}
{\mathrm{C}}_*(Y)\geq {\mathrm{H}}_{q,N}(f,Y)\quad\mbox{and}\quad {\mathrm{C}%
}_*(Y)\geq {\mathrm{C}}_{2l_{0}(f)}^{1/2}(Y) {\mathrm{Q}}%
_{q+(N+1)(N+3),q}(Y).
\end{equation*}
$\square$

\subsubsection{Proof of Corollary \ref{coeff}}
\label{coefficients}

We first explicitly write the expression of the polynomials $%
H_{\Gamma_{n,k}}(x)$ for $k=1,2,3$. Recall formula (\ref{M2}) for the $k$th
operator $\Gamma_{n,k}$ and recall formula (\ref{M1}) for the set $%
\Lambda_{m,k}$ appearing in (\ref{M2}). Recall also formula (\ref{coeff-1})
for $c_n(\alpha)$ and $d_n(\alpha,\beta)$.

\smallskip

\textbf{Case} $k=1.$ Then $m=1$ and $l=3,l^{\prime }=0.$ So the first order
operator is given by%
\begin{equation*}
\Gamma_{n,1}=\frac 1n\sum_{r=1}^{n}\frac 1{6}D_{n,r}^{(3)} =\frac
1{6n}\sum_{r=1}^{n}\sum_{\left\vert \alpha \right\vert
=3}\Delta_{n,r}(\alpha)\partial_{\alpha },
\end{equation*}%
so that, with $c_n(\alpha)$ given in (\ref{coeff-1}), 
\begin{equation}  \label{Gamma1-1}
H_{\Gamma_{n,1}}(x) =\frac 1{6}\sum_{\left\vert \alpha \right\vert
=3}c_n(\alpha)H _{\alpha }(x) =\mathcal{H} _{n,1}(x)
\end{equation}
and formula (\ref{H1}) holds.

\smallskip

\textbf{Case} $k=2.$ Then $m=1$ or $m=2$, and we call $\Gamma^{\prime
}_{n,2} $ and $\Gamma^{\prime \prime }_{n,2}$ the corresponding operator.
Suppose first that $m=1.$ Then we need that $l+2l^{\prime }=k+2m=4.$ This
means that we have $l=4,l^{\prime }=0.$ Then 
\begin{equation*}
\Gamma^{\prime }_{n,2}=\frac 1n\sum_{r=1}^n\frac 1{24}D_{n,r}^{(4)} =\frac
1{24\,n}\sum_{r=1}^n\sum_{|\alpha|=4}\Delta_{n,r}(\alpha)\partial_\alpha
=\frac 1{24}\sum_{|\alpha|=4}c_n(\alpha)\partial_\alpha.
\end{equation*}%
%
%
%

Suppose now that $m=2.$ Then we need that $l_{1}+l_{2}+2(l_{1}^{\prime
}+l_{2}^{\prime })=k+2m=6.$ The only possibility is $l_{1}=l_{2}=3$, $%
l_{1}^{\prime }=l_{2}^{\prime }=0$ and the corresponding term is%
\begin{align*}
\Gamma^{\prime \prime }_{n,2} &=\frac 1{n^2}\sum_{0\leq r_{1}<r_{2}\leq n}%
\frac{1}{36}D_{n,r_{1}}^{(3)}D_{n,r_{2}}^{(3)} =\frac{1}{36\,n^2}\sum_{1\leq
r_{1}<r_{2}\leq n}\sum_{\left\vert \alpha \right\vert =3}\sum_{\left\vert
\beta \right\vert =3}\Delta _{n,r_{1}}(\alpha )\Delta _{n,r_{2}}
(\beta)\partial_{\alpha}\partial_\beta \\
&=\frac{1}{72\,n^2}\sum_{1\leq r_{1}\neq r_{2}\leq n}\sum_{\left\vert \alpha
\right\vert =3}\sum_{\left\vert \beta \right\vert =3}\Delta
_{n,r_{1}}(\alpha )\Delta _{n,r_{2}} (\beta)\partial_{\alpha}\partial_\beta.
\end{align*}%
We notice that, for $|\alpha|=|\beta|=3$, 
\begin{equation*}
\frac 1{n^2}\sum_{1\leq r_{1}\neq r_{2}\leq n}\Delta _{n,r_{1}}(\alpha
)\Delta _{n,r_{2}} (\beta) =c_n(\alpha)c_n(\beta)-\frac 1nd_n(\alpha,\beta)
\end{equation*}
with 
\begin{equation}  \label{sup-1}
\sup_n|d_n(\alpha,\beta)|\leq 4{\mathrm{C}}_{3}^2(Y),\quad
|\alpha|=|\beta|=3.
\end{equation}
So, by inserting, 
\begin{equation*}
\Gamma^{\prime \prime }_{n,2}= \frac 1{72}\sum_{\left\vert \alpha
\right\vert =3}\sum_{\left\vert \beta \right\vert
=3}c_n(\alpha)c_n(\beta)\partial_\alpha\partial_\beta -\frac
1{72\,n}\sum_{\left\vert \alpha \right\vert =3}\sum_{\left\vert \beta
\right\vert =3}d_n(\alpha,\beta)\partial_\alpha\partial_\beta.
\end{equation*}
We conclude that%
\begin{align}
H_{\Gamma_{n,2}}(x) &=H_{\Gamma^{\prime }_{n,2}}(x)+H_{\Gamma^{\prime \prime
}_{n,2}}(x) =\mathcal{H} _{n,2}(x)- \frac 1{72\,n}\sum_{\left\vert \alpha
\right\vert =3}\sum_{\left\vert \beta \right\vert
=3}d_n(\alpha,\beta)H_{(\alpha,\beta)}(x),  \label{Gamma2-1}
\end{align}
$\mathcal{H} _{n,2}(x)$ being given in (\ref{H2}).

\smallskip

\textbf{Case} $k=3.$ $m=1.$ We need that $l+2l^{\prime }=k+2m=5.$ So $%
l=3,l^{\prime }=1$ or $l=5,l^{\prime }=0.$ The operator term corresponding
to $l=3,l^{\prime }=1$ is 
\begin{align*}
\Gamma^1_{n,3} &=-\frac{1}{12\,n}\sum_{r=1}^{n}D_{n,r}^{(3)}L_{\sigma
_{n,r}}^{1} =-\frac{1}{12}\sum_{\left\vert \alpha \right\vert
=3}\sum_{i,j=1}^{d}\overline{c}_n(\alpha,i,j)\partial _{\alpha }\partial
_{i}\partial _{j},
\end{align*}
$\overline{c}_n(\alpha,i,j)$ being given in (\ref{coeff-1}). The term
corresponding to $l=5,l^{\prime }=0$ is 
\begin{equation*}
\Gamma^2_{n,3}=\frac 1n\sum_{r=1}^n\frac 1{5!} D^{(5)}_{n,r} =\frac{1}{5!\,n}%
\sum_{r=1}^{n}\sum_{\left\vert \alpha \right\vert =5}\Delta _{n,r
}(\alpha)\partial _{\alpha } =\frac{1}{5!}\sum_{\left\vert \alpha
\right\vert =5}c_n(\alpha)\partial _{\alpha }.
\end{equation*}

$m=2.$ We need $l_{1}+l_{2}+2(l_{1}^{\prime }+l_{2}^{\prime })=k+2m=7.$ The
only possibility is $l_{1}=3$, $l_{2}=4$, $l_{1}^{\prime }=l_{2}^{\prime }=0$
and $l_{1}=4$, $l_{2}=3$, $l_{1}^{\prime }=l_{2}^{\prime }=0$. The
corresponding term is 
\begin{align*}
\Gamma^3_{n,3} &=\frac 2{n^2}\sum_{1\leq r_{1}<r_{2}\leq n}\frac
1{3!}D_{n,r_{1}}^{(3)}\frac 1{4!}D_{n,r_{2}}^{(4)} =\frac{1}{3! 4!}%
\sum_{\left\vert \alpha \right\vert =3}\sum_{\left\vert \beta \right\vert =4}%
\Big[c_n(\alpha)c_n(\beta)-\frac 1n d_n(\alpha,\beta)\Big]\partial _{\alpha
}\partial _{\beta },
\end{align*}
with 
\begin{equation}  \label{sup-2}
\sup_n|d_n(\alpha,\beta)|\leq C\,{\mathrm{C}}_3(Y){\mathrm{C}}_4(Y)\leq C\,{%
\mathrm{C}}^2_4(Y),\quad |\alpha|=3,|\beta|=4.
\end{equation}

$m=3.$ We need $l_{1}+l_{2}+l_{3}+2(l_{1}^{\prime }+l_{2}^{\prime
}+l_{3}^{\prime })=k+2m=3+6=9.$ The only possibility is $%
l_{1}=l_{2}=l_{3}=3,l_{1}^{\prime }=l_{2}^{\prime }=l_{3}^{\prime }=0$ and
the corresponding term is 
\begin{align*}
\Gamma^4_{n,3} &=\frac{1}{6^{3}\,n^3}\!\!\!\!\sum_{1\leq
r_{1}<r_{2}<r_{3}\leq
n}\!\!\!\!D_{n,r_{1}}^{(3)}D_{n,r_{2}}^{(3)}D_{n,r_{3}}^{(3)} \\
&=\frac{1}{6^{3}\,n^3\,3!}\!\!\sum_{1\leq r_{1}\neq r_{2}\neq r_{3}\leq n}
\sum_{\left\vert \alpha \right\vert =3}\sum_{\left\vert \beta \right\vert
=3}\sum_{\left\vert \gamma \right\vert =3}\Delta _{n,r_1}(\alpha)\Delta
_{n,r_2}(\beta)\Delta _{n,r_3 }(\gamma)\partial _{\alpha }\partial _{\beta
}\partial _{\gamma },
\end{align*}%
where the notation $r_{1}\neq r_{2}\neq r_{3}$ means that $r_1,r_2,r_3$ are
all different. Now, straightforward computations give 
\begin{equation*}
\frac 1{n^3}\sum_{1\leq r_{1}\neq r_{2}\neq r_3\leq n}\Delta
_{n,r_{1}}(\alpha )\Delta _{n,r_{2}} (\beta)\Delta _{n,r_{2}} (\gamma)
=c_n(\alpha)c_n(\beta)c_n(\gamma)-\frac 1ne_n(\alpha,\beta,\gamma),
\end{equation*}
with 
\begin{equation}  \label{sup-3}
\sup_n|e_n(\alpha,\beta,\gamma)|\leq C\times {\mathrm{C}}_3^3(Y),\quad
|\alpha|=|\beta|=|\gamma|=3.
\end{equation}
So, we obtain 
\begin{align*}
\Gamma^4_{n,3} &=\frac{1}{6^{3}\,n^3\,3!}\sum_{\left\vert \alpha \right\vert
=3}\sum_{\left\vert \beta \right\vert =3}\sum_{\left\vert \gamma \right\vert
=3} c_n(\alpha)c_n(\beta)c_n(\gamma)\partial _{\alpha }\partial _{\beta
}\partial _{\gamma } -\frac 1n \sum_{\left\vert \alpha \right\vert
=3}\sum_{\left\vert \beta \right\vert =3}\sum_{\left\vert \gamma \right\vert
=3} e_n(\alpha,\beta,\gamma)\partial_\alpha\partial_\beta\partial_\gamma.
\end{align*}%
We conclude that 
\begin{equation}  \label{Gamma3-1}
H_{\Gamma_{n,3}}(x) =\sum_{i=1}^4H_{\Gamma^i_{n,3}}(x) =\mathcal{H}_{n,3}(x)
+\frac 1n\sum_{|\alpha|\leq 9}f_n(\alpha)H_{\alpha}(x),
\end{equation}
with $\mathcal{H}_{n,3}(x)$ as in (\ref{H3}) and with 
\begin{equation}  \label{sup-4}
\sup_n|f_n(\alpha)|\leq C{\mathrm{C}}_4^3(Y).
\end{equation}
By resuming, we get 
\begin{equation*}
\Phi_{n,N}(x) =1+\sum_{k=1}^3 \frac 1{n^{k/2}}\mathcal{H}_{n,k}(x) +\frac
1{n^2}\mathcal{P}_n(x),
\end{equation*}
where, taking into account (\ref{sup-1}), (\ref{sup-2}), (\ref{sup-3}) and (%
\ref{sup-4}), 
\begin{equation*}
\mathcal{P}_n(x) =\sum_{|\alpha|\leq 9}g_n(\alpha)H_{\alpha}(x)\quad%
\mbox{with}\quad \sup_n|g_n(\alpha)|\leq C\times {\mathrm{C}}_4^3(Y),
\end{equation*}
where $C$ a universal constant. Therefore 
\begin{align*}
&\Big\vert {\mathbb{E}}\Big(\partial _{\gamma }f(S_{n}(Y))\Big)-{\mathbb{E}}%
\Big(\partial _{\gamma }f(W)\Big(1+\sum_{k=1}^3\frac 1{n^{k/2}}\mathcal{H}%
_{n,k}(W)\Big)\Big)\Big\vert \\
&\leq \left\vert {\mathbb{E}}\big(\partial _{\gamma }f(S_{n}(Y))\big)-{%
\mathbb{E}}\big(\partial _{\gamma }f(W)\Phi _{n,N}(W)\big)\right\vert +
\frac 1{n^2}\left\vert {\mathbb{E}}\big(\partial _{\gamma }f(W)\mathcal{P}%
_n(W)\big)\right\vert \\
&\leq C\times {\mathrm{C}}_{\ast }(Y)\Big(\frac{L_{0}(f)}{n^{\frac{1}{2}%
(N+1)}}+L_{q}(f)e^{-\frac{{\mathfrak{m}}^2_r }{32}\times n}\Big) +\frac
1{n^2}\left\vert {\mathbb{E}}\big(\partial _{\gamma }f(W)\mathcal{P}_n(W)%
\big)\right\vert
\end{align*}%
in which we have used (\ref{s20}). Now, by using the standard integration by
parts for the Gaussian law, we have 
\begin{align*}
\left\vert {\mathbb{E}}\big(\partial _{\gamma }f(W)\mathcal{P}_n(W)\big)%
\right\vert =\left\vert {\mathbb{E}}\big(f(W)G_\gamma(W,\mathcal{P}_n(W))%
\big)\right\vert \leq \|f(W)\|_2\|G_\gamma(W,\mathcal{P}_n(W))\|_2,
\end{align*}
where $G_\gamma(W,\mathcal{P}_n(W))$ denote the weight from the integration
by parts formula. Since $\mathcal{P}_n$ is a linear combination of Hermite
polynomials with bounded coefficients, $\|G_\gamma(W,\mathcal{P}%
_n(W))\|_2\leq C$, $C$ depending on $q$. Moreover, $|f(W)|\leq
L_0(f)(1+|W|^{l_0(f)})$, so 
\begin{align*}
\left\vert {\mathbb{E}}\big(\partial _{\gamma }f(W)\mathcal{P}_n(W)\big)%
\right\vert \leq C\times L_0(f){\mathrm{C}}_{2l_0(f)}^{l_0(f)}(Y).
\end{align*}
The statement now follows. $\square$

\subsubsection{Proof of Proposition \ref{W}}

\label{sect-proof-W}

The idea is that, since $\sum_{k=1}^{n-m}\sigma _{k}\geq \frac{1}{2}I,$ the
random variables $Y_{k},k\leq n-m$ contain sufficient noise in order to give
the regularization effect.

We show the main changes in the estimate of $I_{2}(f)$ (for $I_{1}(f),$ $%
I_{3}(f)$ the proof is analogues). We split $%
P_{r_{N+1}+1,n}^{Z,n}=P_{r_{N+1}+1,n-m}^{Z,n}P_{n-m,n}^{Z,n}$ and we need to
have sufficient noise in order that $P_{r_{N+1}+1,n-m}^{Z,n}$ gives the
regularization effect. Then, the two cases described in (\ref{s24}) and (\ref%
{s25}) are replaced now by $\sum_{i=r_{N+1}+1}^{n-m}\sigma _{n,i}\geq \frac{1%
}{4}I$ and $\sum_{i=1}^{r_{N+1}+1}\sigma _{n,i}\geq \frac{1}{4}I$
respectively. And the condition (\ref{s24'}) becomes $\sum_{i=1}^{N+1}\sigma
_{n,r_{i}}\leq \frac{1}{8}I.$ Then the proof follows exactly the same line. $%
\square $

\appendix

\section{Norms}

\label{app:norms}

The aim of this section is to prove Lemma \ref{L1}. For $F=(F_1,\ldots,F_d)$
We work with the norms 
\begin{eqnarray*}
\left\vert F\right\vert _{1,k} &=&\sum_{j=1}^d\sum_{1=1}^{k}\left\vert
D^{i}F_j\right\vert _{\mathcal{H}^{\otimes i}},\quad \quad \left\vert
F\right\vert _{k}=\left\vert F\right\vert +\left\vert F\right\vert _{1,k} \\
\left\Vert F\right\Vert _{1,k,p} &=&\Vert \left\vert F\right\vert
_{1,k}\Vert _{p},\quad \quad \left\Vert F\right\Vert _{k,p}=\left\Vert
F\right\Vert _{p}+\left\Vert F\right\Vert _{1,k,p}.
\end{eqnarray*}

To begin we give several easy computational rules:%
\begin{eqnarray}
\left\vert FG\right\vert _{k} &\leq &C\sum_{k_{1}+k_{2}=k}\left\vert
F\right\vert _{k_{1}}\left\vert G\right\vert _{k_{2}}  \label{Nn1} \\
\left\vert \left\langle DF,DG\right\rangle \right\vert _{k} &\leq
&C\sum_{k_{1}+k_{2}=k}\left\vert F\right\vert _{1,k_{1}+1}\left\vert
G\right\vert _{1,k_{2}+1}\quad and\quad  \label{Nn2} \\
\left\vert \frac{1}{G}\right\vert _{k} &\leq &\frac{C}{\left\vert
G\right\vert }\sum_{l=0}^{k}\frac{\left\vert G\right\vert _{k}^{l}}{%
\left\vert G\right\vert ^{l}}.  \label{Nn3}
\end{eqnarray}%
Now, for $F=(F_{1},\ldots ,F_{d})$ we consider the Malliavin covariance
matrix $\sigma _{F}^{i,j}=\left\langle DF^{i},DF^{j}\right\rangle $ and, if $%
\det \sigma _{F}\neq 0,$ we denote $\gamma _{F}=\sigma _{F}^{-1}.$ We write%
\begin{equation*}
\gamma _{F}^{i,j}=\frac{\widehat{\sigma }_{F}^{i,j}}{\det \sigma _{F}}
\end{equation*}%
where $\widehat{\sigma }_{F}^{i,j}$ is the algebraic complement . Then,
using (\ref{Nn1}) 
\begin{equation*}
\big\vert \gamma _{F}^{i,j}\big\vert _{k}\leq C\sum_{k_{1}+k_{2}=k}\big\vert 
\widehat{\sigma }_{F}^{i,j}\big\vert _{k_{1}}\left\vert \frac{1}{\det \sigma
_{F}}\right\vert _{k_{2}}.
\end{equation*}%
By (\ref{Nn1}) and (\ref{Nn2}), $\left\vert \widehat{\sigma }%
_{F}^{i,j}\right\vert _{k_{1}}\leq C\left\vert F\right\vert
_{1,k_{1}+1}^{2(d-1)}$ and $\left\vert \det \sigma _{F}\right\vert
_{k_{2}}\leq C\left\vert F\right\vert _{1,k_{2}+1}^{2d}.$ Then, using (\ref%
{Nn3}) 
\begin{equation*}
\left\vert \frac{1}{\det \sigma _{F}}\right\vert _{k_{2}}\leq \frac{C}{%
\left\vert \det \sigma _{F}\right\vert }\sum_{l=0}^{k_{2}}\frac{\left\vert
\det \sigma _{F}\right\vert _{k_{2}}^{l}}{\left\vert \det \sigma
_{F}\right\vert ^{l}}\leq \frac{C}{\left\vert \det \sigma _{F}\right\vert }%
\sum_{l=0}^{k_{2}}\frac{\left\vert F\right\vert _{1,k_{2}+1}^{2ld}}{%
\left\vert \det \sigma _{F}\right\vert ^{l}}
\end{equation*}%
so that 
\begin{equation}
\big\vert \gamma _{F}^{i,j}\big\vert _{k}\leq C\frac{\left\vert F\right\vert
_{1,k+1}^{2(d-1)}}{\left\vert \det \sigma _{F}\right\vert }\sum_{l=0}^{k}%
\Big(\frac{\left\vert F\right\vert _{1,k+1}^{2d}}{\left\vert \det \sigma
_{F}\right\vert }\Big)^{l}\leq C\frac{\left\vert F\right\vert
_{1,k+1}^{2(d-1)}}{\left\vert \det \sigma _{F}\right\vert }\Big(1+\frac{%
\left\vert F\right\vert _{1,k+1}^{2d}}{\left\vert \det \sigma
_{F}\right\vert }\Big)^{k}  \label{Nn4}
\end{equation}

We denote%
\begin{equation}
\alpha _{k}=\frac{\left\vert F\right\vert _{1,k+1}^{2(d-1)}(\left\vert
F\right\vert _{1,k+1}+\left\vert LF\right\vert _{k})}{\left\vert \det \sigma
_{F}\right\vert },\quad \beta _{k}=\frac{\left\vert F\right\vert
_{1,k+1}^{2d}}{\left\vert \det \sigma _{F}\right\vert }  \label{Nn4'}
\end{equation}%
and%
\begin{equation}
{\mathcal{K}}_{n,k}(F)=(\left\vert F\right\vert _{1,k+n+1}+\left\vert
LF\right\vert _{k+n})^{n}(1+\left\vert F\right\vert _{1,k+n+1})^{2d(2n+k)}.
\label{Nn4''}
\end{equation}%
We also recall that for $\eta >0,$ we consider a function $\Psi _{\eta }\in
C^{\infty }({\mathbb{R}})$ such that $1_{(0,\eta )}\leq \Psi _{\eta }\leq
1_{(0,2\eta )} $ and $\Vert \Psi _{\eta }^{(k)}\Vert _{\infty }\leq
C_{k}\eta ^{-k},\forall k\in {\mathbb{N}}.$ Then we take $\Phi _{\eta
}=1-\Psi _{\eta }.$

\begin{lemma}
\textbf{A}. For every $k,n\in {\mathbb{N}}$ there exists a universal
constant $C$ (depending on $k$ and $n)$ such that, for $\omega $ such that $%
\det \sigma _{F}(\omega )>0,$ 
\begin{equation}
\left\vert H_{\rho }^{(n)}(F,G)\right\vert _{k}\leq C\alpha
_{k+n}^{n}\sum_{p_{1}+p_{2}=k+n}\left\vert G\right\vert _{p_{2}}(1+\beta
_{k+n})^{p_{1}}.  \label{Nn5}
\end{equation}%
\textbf{B}. For every $\eta >0$ 
\begin{equation}
\left\vert H_{\rho }^{(n)}(F,\Phi _{\eta }(\det \sigma _{F})G)\right\vert
_{k}\leq \frac{C}{\eta ^{2n+k}}\times {\mathcal{K}}_{n,k}(F)\times
\left\vert G\right\vert _{k+n}.  \label{Nn6}
\end{equation}
\end{lemma}

\bigskip

\textbf{Proof A}. We first prove (\ref{Nn5}) for $n=1$. We have 
\begin{equation*}
H_{i}^{(1)}(F,G)=-\sum_{j=1}^{m}G\gamma _{F}^{i,j}LF^{j}+G\langle D\gamma
_{F}^{i,j},DF^{j}\rangle +\gamma _{F}^{i,j}\langle DG,DF^{j}\rangle.
\end{equation*}%
Using (\ref{Nn1}) 
\begin{align*}
&\big\vert H_{i}^{(1)}(F,G)\big\vert _{k} \\
&\leq C\sum_{k_{1}+k_{2}+k_{3}=k}\left( \left\vert \gamma _{F}\right\vert
_{k_{1}}\left\vert LF\right\vert _{k_{2}}\left\vert G\right\vert
_{k_{3}}+\left\vert \gamma _{F}\right\vert _{k_{1}+1}\left\vert F\right\vert
_{1,k_{2}+1}\left\vert G\right\vert _{k_{3}}+\left\vert \gamma
_{F}\right\vert _{k_{1}}\left\vert F\right\vert _{1,k_{2}+1}\left\vert
G\right\vert _{k_{3}+1}\right) \\
&\leq C(\left\vert F\right\vert _{k+1}+\left\vert LF\right\vert
_{k})\sum_{p_{1}+p_{2}\leq k}\left( \left\vert \gamma _{F}\right\vert
_{p_{1}+1}\left\vert G\right\vert _{p_{2}}+\left\vert \gamma _{F}\right\vert
_{p_{1}}\left\vert G\right\vert _{p_{2}+1}\right) .
\end{align*}%
For $n>1$, we use recurrence and we obtain 
\begin{equation*}
\left\vert H_{\gamma }^{(n)}(F,G)\right\vert _{k}\leq C(\left\vert
F\right\vert _{k+n+1}+\left\vert LF\right\vert _{k+n})^{n}\sum_{p_{1}+\ldots
+p_{n+1}\leq k+n-1}\prod_{i=1}^{n}\left\vert \gamma _{F}\right\vert
_{p_{i}}\times \left\vert G\right\vert _{p_{n+1}}.
\end{equation*}%
Then, using (\ref{Nn1}) first and (\ref{Nn4}) secondly, (\ref{Nn5}) follows.

\medskip

\textbf{B}. Let $G_{\eta }=\Phi _{\eta }(\det \sigma _{F})G).$ For every $%
p\in {\mathbb{N}}$ one has $\left\vert G_{\eta }\right\vert _{p}\leq C\eta
^{-p}\left\vert G\right\vert _{p}\left\vert F\right\vert _{1,p+1}^{d}.$
Moreover one has $H_{\rho }^{(n)}(F,G_{\eta })=1_{\{\det \sigma _{\Phi
}>\eta /2\}}H_{\rho }^{(n)}(F,G_{\eta }).$ So (\ref{Nn5}) implies (\ref{Nn6}%
). $\square $

\end{document}